\documentclass[11pt]{amsart}
\usepackage{amssymb}
\usepackage{citesort}

\newcounter{mtheorem}
\newtheorem{mtheorem}[mtheorem]{Theorem}

\newtheorem{theorem}{Theorem}[section]
\newtheorem{lemma}[theorem]{Lemma}
\newtheorem{prop}[theorem]{Proposition}
\newtheorem{corollary}[theorem]{Corollary}
\newtheorem{conjecture}[theorem]{Conjecture}

\theoremstyle{definition}

\newtheorem{example}[theorem]{Example}

\theoremstyle{remark}
\newtheorem{remark}[theorem]{Remark}

\numberwithin{equation}{section}
\newcommand{\beq}{\begin{equation}}
\newcommand{\eeq}{\end{equation}}
\newcommand{\bea}{\begin{eqnarray}}
\newcommand{\eea}{\end{eqnarray}}
\newcommand{\cy}{{C}alabi-{Y}au\ }
\newcommand{\ka}{{K}\"ahler\ }
\newcommand{\ke}{{K}\"ahler-Einstein\ }

\newcommand{\slg}{SLG\ }
\newcommand{\Slg}{{S}pecial {L}agrangian\ }

\newcommand{\tslg}{$\theta$-SLG\ }
\newcommand{\tsleg}{$\theta$-special {L}egendrian\ }
\newcommand{\tslegg}{$\theta$-special {L}egendrian}
\newcommand{\ac}{AC\ }
\newcommand{\acslg}{ACSLG\ }
\newcommand{\dsig}{\mathcal{D}_{\Sigma}}

\newcommand{\ix}{\mathcal{I}_{X'}}
\newcommand{\ox}{\mathcal{O}_{X'}}
\newcommand{\mx}{\mathcal{M}_X}
\newcommand{\ml}{\mathcal{M}_L}
\newcommand{\bdyml}{\partial\mathcal{M}_L}
\newcommand{\mlb}{\overline{\mathcal{M}}_L}

\newcommand{\C}{\mathbb{C}}
\newcommand{\CP}{\mathbb{CP}}
\newcommand{\R}{\mathbb{R}}
\newcommand{\Hyp}{\mathbb{H}}
\newcommand{\Q}{\mathbb{Q}}
\newcommand{\Z}{\mathbb{Z}}
\newcommand{\N}{\mathbb{N}}
\newcommand{\ra}{\rightarrow}

\newcommand{\no}{\noindent}
\newcommand{\sn}{\operatorname{sn}}

\newcommand{\cn}{\operatorname{cn}}
\newcommand{\dn}{\operatorname{dn}}
\newcommand{\diag}{\operatorname{diag}}

\newcommand{\vol}{\operatorname{Vol}}

\newcommand{\lind}{\operatorname{l-ind}}
\newcommand{\sind}{\operatorname{s-ind}}
\newcommand{\Real}{\operatorname{Re}}
\newcommand{\Imag}{\operatorname{Im}}
\newcommand{\ind}{\operatorname{ind}}
\newcommand{\Area}{\operatorname{Area}}
\newcommand{\range}{\operatorname{Range}}

\newcommand{\Trace}{\operatorname{Trace}}
\newcommand{\hcf}{\operatorname{hcf}}

\begin{document}

\title{The geometric complexity of special Lagrangian $T^2$-cones}

\author{Mark Haskins}
\address{Department of Mathematics, Johns Hopkins University, Baltimore MD 21218}
\curraddr{Institut des Hautes \'Etudes Scientifiques, Bures-sur-Yvette, France}
\email{haskins@ihes.fr}
\thanks{The author would like to thank Brown University Mathematics Department
for its hospitality in the final stages of writing this paper.\\}

\date{12 December, 2003.}



\begin{abstract}
We prove a number of results on the geometric complexity of special Lagrangian (SLG) $T^2$-cones in $\C^3$.

Every SLG $T^2$-cone has a fundamental integer invariant, its spectral curve genus.
We prove that the spectral curve genus of an \slg $T^2$-cone gives a lower bound for its geometric complexity,
\textit{i.e.} the area, the stability index and the Legendrian index of any SLG $T^2$-cone
are all bounded below by explicit linearly growing functions of the spectral curve genus.

We prove that the cone on the Clifford torus (which has spectral curve genus zero) in $S^5$ is the
unique \slg $T^2$-cone with the smallest possible Legendrian index  and hence that it is the unique stable SLG $T^2$-cone.
This leads to a classification of all rigid ``index 1'' SLG cone types in dimension three.
For cones with spectral curve genus two
we give refined lower bounds for the
area, the Legendrian index and the stability index.
One consequence of these bounds is that there exist
$S^1$-invariant \slg torus cones of arbitrarily large area, Legendrian
and stability indices.

We explain some consequences of our results for the programme (due to Joyce) to understand the
``most common'' three-dimensional isolated singularities of generic families of \slg submanifolds in almost Calabi-Yau manifolds.
\end{abstract}

\maketitle


\section{Introduction}

Let $M$ be a \cy manifold of complex dimension $n$ with \ka form
$\omega$ and non-zero parallel holomorphic $n$-form $\Omega$
satisfying a normalization condition.
Then $\Real{\Omega}$ is a calibrated form whose calibrated
submanifolds are called special Lagrangian (SLG) submanifolds
\cite{HL:1982}.
Moduli spaces of \slg submanifolds have appeared
recently in string theory \cite{becker,SYZ:1996}.
On physical grounds, Strominger, Yau and Zaslow argued
that a \cy manifold $M$ with a mirror partner $\hat{M}$ admits a (singular)
fibration by \slg tori, and that $\hat{M}$ should be obtained by compactifying the dual fibration \cite{SYZ:1996}.
To make this idea rigorous one needs to have control over the singularities and compactness
properties of families of \slg submanifolds. In dimensions three and higher these properties are
not well understood. Throughout this paper we will assume $n\ge 3$.

The local models for isolated singularities of \slg varieties are
SLG cones in $\C^n$. They arise as possible tangent cones to \slg currents at singular points.
The simplest \slg singularities to understand are  isolated singularities modelled on \slg cones with isolated singularities.
Recently a number of authors have used integrable systems techniques to prove existence of large numbers of such \slg cones
\cite{haskins:slgcones,joyce:intsys,mcintosh:slg,carberry:mcintosh}. Even in dimension
three there are too many \slg cones to expect a complete classification.

On the other hand, recent work of Joyce  suggests that to make
progress on two important open questions in \slg geometry one does
not need to understand all singularity types, only the ``more
common ones'' \cite{joyce:slgcones5}. The primary goal of this paper is to shed light on
which \slg $T^2$-cones are the ``common'' ones.

When properly understood this amounts to studying \slg cones for
which certain natural indices (the Legendrian index and the stability index)
constructed out of the spectral geometry of the cone $C$ are small.
We find several methods of proving lower bounds for these indices
and use them to show that many \slg cones one can construct have large indices
and hence are ``uncommon''. In certain cases, we can also prove upper bounds
for these indices and conclude that certain singularities should be relatively ``common''.

We now give a description of the main results of this paper.\\

The \textit{Legendrian index} of a \slg cone $C$ in $\C^n$ is the number of eigenvalues of the Laplacian acting on functions of the
link, $\Sigma = C \cap S^{2n-1}$, that are contained in the open interval $(0,2n)$. The Legendrian index has a simple geometric interpretation:
it is the Morse index of the link $\Sigma$ as a critical point of volume, but restricted to Legendrian variations.
For any \slg cone in $\C^n$ with an isolated singularity the Legendrian index is at least $2n$ (Corollary \ref{lb:index}).
A \slg cone is called \textit{(Legendrian) stable} if its Legendrian index is equal to $2n$.

The \textit{stability index} of a \slg cone in $\C^n$ counts the \textit{excess} number of eigenvalues of $\Delta$ in the closed interval $[0,2n]$
\textit{i.e.} the number of eigenvalues in this interval minus those eigenvalues which must occur for simple
geometric reasons (see \S \ref{harm:fn:cone} for a complete description).
A \slg cone is called \textit{strictly stable} if it has stability index equal to $0$.
Equivalently, a \slg cone is strictly stable if it is both (Legendrian) stable and
\textit{rigid}, \textit{i.e.} the eigenvalue $\lambda=2n$ occurs with the minimum geometrically
allowable multiplicity.
More generally, the spectrum of $\Delta$ on the link $\Sigma$ of the cone $C$
plays a fundamental role in many analytic and geometric problems arising in \slg geometry,
as we try to explain in Appendix \ref{analysis:spectral:geometry}.

For a large class of isolated singularities of compact singular \slg submanifolds one can define
the \textit{index of the singularity}. This index controls how commonly that singularity type
occurs. Hence singularities of small index are the most important ones to understand. Moreover, for some outstanding problems in global \slg geometry
(\textit{e.g.} counting \slg homology spheres \cite{joyce:counting} and the SYZ conjecture \cite{joyce:lectures})
it should be sufficient to understand only singularities of index up to some fixed number $k$.
In particular, $k=1$ for the counting \slg homology spheres problem.

A formula for this index due to Joyce \cite[Thm 8.10]{joyce:slgcones5},
given in (\ref{index:transverse}), shows that the index of any singularity modelled on a \slg cone $C$ with
large stability index is also large, and therefore does not occur often.
For any \slg cone in $\C^n$ with an isolated singularity
the stability index is nonnegative (see \S \ref{harm:fn:cone}).
Our first result characterizes the $T^2$-cones with stability index equal to zero.
\begin{mtheorem}
The cone on the Clifford torus in $S^5$ is (up to special unitary equivalence) the unique strictly stable \slg $T^2$-cone.
\end{mtheorem}
\noindent
In other words, the cone on the Clifford torus is ``the most common'' \slg $T^2$-cone.
The Clifford torus $T^{n-1}$ in $S^{2n-1}$ is a flat minimal Legendrian torus in $S^{2n-1}$ which
is invariant under the diagonal subgroup of $\text{SU(n)}$. It is described in detail in Example \ref{clifford:torus}.
Strictly stable singularities have other particularly nice properties (see Appendix \ref{conical:sing:section} for a brief discussion).
Unfortunately, the analogue of Theorem \ref{ct:is:only:stable:t2} is false for all higher dimensions,
because the Clifford torus fails to be strictly stable (Prop \ref{ct:stability}).
In fact, no examples of stable \slg singularities in dimension greater than three appear to be currently known.

Theorem \ref{ct:is:only:stable:t2} combined with an argument of Joyce \cite[\S 10.3]{joyce:lectures}
leads immediately to Corollary \ref{index:1:classification} which classifies \textit{all} rigid \slg cones
which can occur in an index $1$ singularity in three dimensions.

There is an analogue of Theorem \ref{ct:is:only:stable:t2} for minimal Lagrangian tori in $\CP^2$
due to Urbano.
\begin{mtheorem}[\mbox{Urbano, \cite[Cor. 2]{urbano:index}}]
The Clifford torus in $\CP^2$ is (up to unitary equivalence) the unique
Hamiltonian stable minimal Lagrangian torus in $\CP^2$.
\end{mtheorem}
\noindent
Although Theorems A and B are similar there does not appear to be a way to prove one from the other.
We give a new proof of Theorem B which is considerably simpler than Urbano's original proof.

Dimension three is also special because \slg $T^2$-cones can be studied using methods of
integrable systems/algebraic curves.
Since the fundamental work of Hitchin on harmonic $2$-tori in $S^3$ \cite{hitchin}, it has become
clear that harmonic maps of $2$-tori into many compact Lie groups or symmetric spaces
can be completely classified in terms of  algebro-geometric data, the so-called
\textit{spectral data} of the harmonic torus.
In particular, \slg $T^2$-cones are in one-to-one correspondence with certain spectral data
(see \S \ref{int:systems}). Part of this spectral data consists of a compact Riemann surface of finite even genus $2d$.

In theory, the spectral data completely determines the geometry of the harmonic $2$-torus.
A fundamental problem in this field originally posed by Hitchin \cite{hitchin} is to
recover differential geometric properties of the torus from its spectral data.
For example, for CMC tori, numerical work suggested that the spectral curve genus provides
a lower bound for the geometric complexity of the corresponding torus.
For CMC tori in $\R^3$ and minimal tori in $S^3$ a breakthrough in this problem was achieved very recently  by Ferus \textit{et. al.} \cite{pedit:pinkall}.
Using methods of ``quaternionic holomorphic'' geometry they proved, for example, that the area of a minimal torus in $S^3$
increases with the genus of its spectral curve \cite[Thm 6.7]{pedit:pinkall}.
Two of our main theorems prove analogous results for \slg $T^2$-cones, using very different methods to those of \cite{pedit:pinkall}.
\begin{mtheorem}[$\lind$ and $\sind$ grow at least linearly with spectral curve genus]
Let $C$ be a \slg $T^2$-cone in $\C^3$ with an isolated singularity at $0$,
and let $g=2d$ denote the spectral curve genus of the associated minimal Lagrangian
torus $\phi: T^2 \ra \CP^2$. Then for any $d\ge 3$ we have
\begin{equation*}
\lind(C) \ge \max\left\{\left[\tfrac{1}{2}d\right]^+,7\right\}
\end{equation*}
and
\begin{equation*}
\sind(C) \ge \max\left\{\left[\tfrac{3}{2}d-8\right]^+,d-1 \right\},
\end{equation*}
where  $[x]^+$ denotes the smallest integer not less than $x\in \R$, \textit{i.e.} $[x]^+:=\min{\{n\in \Z \, |\, n \ge x \}}$.
\end{mtheorem}

\begin{mtheorem}[Area grows at least linearly with spectral curve genus]
Let $C$ be a \slg $T^2$-cone in $\C^3$ with an isolated singularity at $0$,
and let $g=2d$ denote the spectral curve genus of the associated minimal Lagrangian
torus $\phi: T^2 \ra \CP^2$. Let $\Sigma$ be the link of the cone $C$. Then for any $d\ge 3$ we have
\begin{equation*}
\Area(\Sigma) > \tfrac{1}{3} d\pi.
\end{equation*}
\end{mtheorem}
\noindent
We conjecture that $\lind(\Sigma)$, $\sind(\Sigma)$ and $\Area(\Sigma)$ in fact grow quadratically with $d$
(see Remarks \ref{quadratic:index:rmk} and \ref{quadratic:area:growth}).

Special Lagrangian $T^2$-cones with an $S^1$-symmetry can be described very explicitly using elliptic
functions and integrals. In \S \ref{s1:invariant:tori} and \S \ref{areas} we use these explicit descriptions to
refine our lower bounds for $\lind$, $\sind$ and $\Area$ in the $S^1$-invariant case.
Two of the main results from \S \ref{s1:invariant:tori} are the following.
\begin{mtheorem}[$\lind$ and $\sind$ for $T_{m,n}$ increase linearly in $m$ and $n$]
\noindent
Let $T_{m,n}$ be one of the countably infinite family of minimal Legendrian $2$-tori described in Theorem \ref{u:alpha:0}. For $\lambda >0$, denote by $N_{m,n}(\lambda)$
the number of eigenvalues of $\Delta$ on $T_{m,n}$ in the range $[0,\lambda(\lambda+1)]$.
Then
\begin{equation*}
 N_{n,m}(1)   \ge \left\{ \begin{array}{ll}
                             4n+5  \ge 13 & \text{if mn is even}\\
                             2n+5  \ge 11 & \text{if mn is odd}.\\
                             \end{array}
                             \right.
\end{equation*}
\begin{equation*}
 \sind(T_{m,n}) \ge \left\{ \begin{array}{ll}
                                 8n + 4m -8  \ge 12 & \text{if mn is even}\\
                                 4n + 2m -8 \ge 6 & \text{if mn is odd}.\\
                                 \end{array}
                                 \right.
\end{equation*}
\begin{equation*}
 \lind(T_{m,n})  \ge \left\{ \begin{array}{ll}
                                 8n+4m-2  \ge 18 & \text{if mn is even}\\
                                 4n+2m-2 \ge 12 & \text{if mn is odd}.\\
                                 \end{array}
                                 \right.
\end{equation*}
Further, we also have the following upper bounds
\begin{equation*}
N_{n,m}(2) < \left\{
    \begin{array}{ll}
    36 \frac{\pi}{\sqrt{3}}n -1 & \text{if $mn$ is even,}\\
    18 \frac{\pi}{\sqrt{3}}n -1 & \text{if $mn$ is odd.}
\end{array}
\right.
\end{equation*}
\end{mtheorem}
\noindent
See Theorem \ref{u:alpha:0} for a description of the tori $T_{m,n}$. They are the simplest family
of $S^1$-invariant minimal Legendrian $2$-tori in $S^5$.

\begin{mtheorem}[bounds on $\lind$ for tori $u_{0,J}$]
Let $\Sigma$ be any of the countably infinite family of minimal Legendrian $2$-tori constructed in
Theorem \ref{u:zero:j}. Then $\frac{1}{2\pi}\theta_2(T_J) = \frac{m}{n}$ for relatively prime
positive integers $m\ge 4$ and $n\ge 7$, and
\begin{eqnarray*}
\sind(\Sigma) & \ge & 2n-8 \ge 6.\\
\lind(\Sigma) & \ge & 2n-2 \ge 12.\\
N_{\Sigma}(1) & \ge & 2m +5 \ge 13.\\
N_{\Sigma}(2) & \ge & 2n +6.
\end{eqnarray*}
\end{mtheorem}
\noindent
See Theorem \ref{u:zero:j} for a description of the tori $u_{0,J}$.

In \S \ref{areas} we use expressions for the area of $S^1$-invariant minimal Legendrian $2$-tori in terms
of elliptic integrals and some properties of elliptic integrals to obtain results
about the area of these tori. One of the main results from \S \ref{areas} is the following.
\begin{mtheorem}[Area bounds for minimal Legendrian tori $T_{m,n}$]
Let $T_{m,n}$ one of the countably infinite family of minimal Legendrian $2$-tori constructed in Theorem \ref{u:alpha:0}.
If $mn$ is even then  $2n < \tfrac{1}{4\pi}\Area(T_{m,n}) < \tfrac{2\pi}{\sqrt{3}}n$,
otherwise $ n < \frac{1}{4\pi}\Area(T_{m,n}) < \frac{\pi}{\sqrt{3}}n$.
\end{mtheorem}
\noindent
An immediate consequence of Theorem \ref{area:tmn} is Corollary \ref{first:eval:bounds} which
gives an upper bound for $\lambda_1(T_{m,n})$ in terms of $m$ and $n$.
Theorem \ref{area:tmn} is also used to prove the upper bounds for $N_{n,m}(2)$
stated in Theorem \ref{lindex:tmn}.\\


The rest of the paper is organized as follows.
In \S \ref{slg:geometry:basics} we recall some basic properties of \slg submanifolds.
In \S \ref{slg:cones:section} we describe the basic
differential and spectral geometry of \slg cones in $\C^n$ and their
links. In particular, the Legendrian index and the stability index
of a \slg cone are defined.

Beginning in \S \ref{dimension:three} we consider the important special case
of \slg $T^2$-cones in $\C^3$. We prove Theorem \ref{ct:is:only:stable:t2}
and its Hamiltonian analogue Theorem \ref{ct:only:hstable:t2} using results from \S \ref{spectral:geometry} and \S \ref{conformal:area:section}.
We review the spectral curve construction of SLG $T^2$-cones in $\C^3$ due to McIntosh \cite{mcintosh:slg}.
This construction associates to any SLG $T^2$-cone a basic integer invariant, its spectral curve genus.
Theorem \ref{specgenus:bds:index} gives lower bounds for the Legendrian  or stability index of any \slg $T^2$-cone in terms of its spectral curve genus.
Similarly, Theorem \ref{area:spectral:genus} gives lower bounds for the area of the link of any \slg $T^2$-cone again
in terms of its spectral curve genus.
In \S \ref{s1:invariant:tori} we recall from \cite{haskins:slgcones} the author's explicit construction
of SLG $T^2$-cones with an $S^1$-symmetry. These cones all have spectral curve genus $2$.
Using the explicit descriptions of these cones, in Theorems \ref{lindex:tmn} and \ref{lindex:u0j}
we prove refined lower bounds for their Legendrian index and stability index.
These results show the impossibility of finding upper bounds for the Legendrian index
or stability index of a \slg $T^2$-cone solely in terms of its spectral curve genus.
In \S \ref{areas} we give analytical and numerical results on the area of the link $\Sigma$
of these $S^1$-invariant torus cones $C$. These results can then be used to prove upper bounds
for certain eigenvalue-related data, \textit{e.g.} $\lambda_1(\Sigma)$ or $\sind{(C)}$.

There are two appendices. Appendix \ref{appendix:elliptic:fn} recalls some basic properties
of Jacobi elliptic functions and integrals that are used in \S \ref{s1:invariant:tori} and \S \ref{areas}.
In Appendix \ref{analysis:spectral:geometry} we describe the implications
of our results on \slg cones for several other problems in \slg geometry.

The author would like to thank the referee for their careful reading of the paper
and for suggesting a number of improvements. The author would also like to thank H. Iriyeh, Y.G. Oh and
F. Urbano for bringing reference \cite{urbano:index} to his attention.

\setcounter{mtheorem}{0}

\section{\Slg geometry in \cy and almost \cy manifolds}
\label{slg:geometry:basics}
In this section we recall some basic facts about
\slg geometry: first in $\C^n$, then generalizing to \cy manifolds
and almost \cy manifolds.

\subsection{Special Lagrangian geometry in $\C^n$}
Special Lagrangian geometry is an example of a \textit{calibrated
geometry}. We begin by reviewing some elementary
facts about calibrations \cite{HL:1982}.

Each calibrated geometry is a distinguished class of minimal
submanifolds of a Riemannian manifold $(M,g)$ associated with a
closed differential $p$-form $\phi$ of comass one.

For each $m\in M$, the \textit{comass} of $\phi$ is defined to be
\begin{displaymath}
\|\phi\|_m^* = \sup\{ <\phi_m,\xi_m> : \xi_m \ \textrm{is a unit
simple } p\textrm{-vector at} \ m\}
\end{displaymath}
\no where a $p$-vector $\xi \in \Lambda^pV$ is \textit{simple} if
$\xi = v_1 \wedge v_2 \wedge \ldots \wedge v_p$ for some
$v_1,\ldots, v_p \in V$. In other words, $\|\phi\|_m^{*}$ is the
supremum of $\phi$ restricted to the Grassmannian of oriented
$p$-dimensional planes $G(p,T_m M)$, regarded as a subset of
$\Lambda^{p}T_mM$.

To any form of comass one there is a natural subset of $G(p,TM)$
\begin{displaymath}
G_m(\phi)  = \{\xi_m \in G(p,T_mM) : <\phi_m,\xi_m> = 1 \} ,
\end{displaymath}
that is, the collection of oriented $p$-planes on which $\phi$
assumes its maximum. These planes are the planes
\textit{calibrated by $\phi$}. An oriented $p$-dimensional
submanifold of $(M,g)$ is \textit{calibrated by $\phi$} if its
tangent plane at each point is calibrated.
The key property of calibrated submanifolds is that they are
\textit{homologically volume minimizing} \cite[\S I]{HL:1982}.
In particular, they are minimal submanifolds \textit{i.e.}
they have zero mean curvature.

 Let $z_1, \ldots, z_n$ denote standard complex coordinates on
$\C^n$. For any $\theta \in [0,2\pi)$ the real $n$-form
$\alpha_{\theta} = \Real(e^{i\theta}dz^1 \wedge \ldots \wedge dz^n)$
is a calibrated form, called the \textit{\tslg calibration} on
$\C^n$.
A \textit{$\theta$-special Lagrangian plane} ($\theta$-SLG) is an oriented $n$-plane
calibrated by the form $\alpha_{\theta}$. A useful
characterization of the \tslg planes is the following.
\begin{lemma}[\mbox{\cite[Cor III.1.11]{HL:1982}}]
\label{slg:alternate}
An oriented $n$-plane $\xi$ in $\C^n$ is \tslg (for the correct
choice of orientation) if and only if
\begin{enumerate}
\item
$\xi$ is Lagrangian with respect to the standard symplectic form
$\omega = \sum{dx^i \wedge dy^i}$,
\item
$\beta_{\theta} := \Imag(e^{i\theta}dz^1 \wedge \ldots
\wedge dz^n)$ restricts to zero on $\xi$.
\end{enumerate}
\end{lemma}
One reason for considering the whole $S^1$-family of \slg
calibrations is the following result.
\begin{prop}[\mbox{\cite[Prop 2.17]{HL:1982}}]
\label{mlag} A connected oriented Lagrangian submanifold $S
\subset \C^n$ is minimal (\textit{i.e.} its mean curvature $H$ vanishes) if and only if $S$ is \tslg
for some $\theta$.
\end{prop}

\smallskip

\subsection{\Slg geometry in \cy manifolds}
\label{cy}
Let $n\ge 3$. A \textit{Calabi-Yau $n$-fold} is a quadruple
$(M,g,J,\Omega)$ consisting of an $n$-dimensional K\"ahler
manifold $(M,g,J)$ with holonomy group $Hol(g)\subseteq \text{SU(n)}$,
and a nonzero parallel holomorphic $(n,0)$-form $\Omega$
satisfying
\begin{equation}
\label{normalize} \omega^n/ n! = (-1)^{n(n-1)/2}\left(\tfrac{1}{2}i\right)^{n}\Omega
\wedge \bar{\Omega}
\end{equation}
 where $\omega$ is the K\"ahler form of $g$. With
the above choice of normalization $\alpha=\Real(\Omega)$ has comass
one with respect to $g$. Hence $(M,g,\alpha)$ is a calibrated
geometry, whose calibrated submanifolds we also call \textit{\slg
submanifolds}.
If we take the quadruple to be $\C^n$ endowed with its Euclidean
metric, standard complex structure and usual holomorphic
$(n,0)$-form $dz^1 \wedge \ldots \wedge dz^n$ then we recover the
definition of $0$-\slg from the previous section.

Another common definition of a \cy $n$-fold is as a \ka manifold
$(M,g,J)$ with holonomy contained in $\text{SU(n)}$. In this case one can
show that there is, up to a choice of complex phase $e^{i\theta}$,
a unique $(n,0)$-form $\Omega$ such that $(M,g,J,\Omega)$ is a
Calabi-Yau manifold as defined previously. With this alternative
definition, each \cy manifold naturally comes with a one-parameter
family of \slg geometries.
In this case the obvious analogue  of Propositions \ref{slg:alternate}
and \ref{mlag} hold.
However, since compact simply connected Calabi-Yau manifolds typically have finite
isometry groups there is generally no obvious relation between \slg submanifolds
with different phases. This is in contrast to $\C^n$
where one can always use a $\text{U(1)}$ motion to transform a \tslg submanifold
into a $0$-\slg submanifold.

\subsection{\Slg geometry in almost \cy manifolds}
We now extend the notion of \slg submanifolds from \cy manifolds
to almost \cy manifolds. In this generalization, \slg submanifolds
are no longer calibrated, nor even minimal submanifolds (at least
with respect to the initially chosen metric).
Nevertheless, certain nice features of \slg submanifolds in \cy
manifolds persist even in this case, \textit{e.g.} the deformation
theory for compact \slg submanifolds remains unobstructed.

The key advantage of the extension to almost \cy manifolds is that while \cy manifolds
come in finite-dimensional families, almost \cy manifolds come in
infinite-dimensional families. Hence if one allows a perturbation
of the original almost \cy structure to a nearby ``generic'' structure
then many pathologies will disappear. A similar idea, namely extending
from integrable complex structures to almost complex structures
has  proven very powerful in the context of symplectic geometry and
topology. The idea of extending special
Lagrangian geometry to the almost \cy setting has appeared in the
work of various authors \cite{bryant:slg,goldstein:calibrated,joyce:lectures}.
Below we follow the presentation in \cite{joyce:lectures}.

Let $n\ge 3$. An \textit{almost Calabi-Yau $n$-fold} is a
quadruple $(M,g,J,\Omega)$ consisting of an $n$-dimensional
K\"ahler manifold $(M,g,J)$, and a nonzero holomorphic $(n,0)$-form $\Omega$.
Given an almost \cy $n$-fold $(M,g,J,\Omega)$ in place of the
relation (\ref{normalize}) there is a unique smooth positive
function $f$ satisfying
\begin{equation}
\label{acy_normalize} \frac{f^{2n}\omega^n}{n!} =
(-1)^{n(n-1)/2}\left(\frac{i}{2}\right)^{n}\Omega \wedge \bar{\Omega}.
\end{equation}
When $f$ is not identically 1, $(M,g,J,\Omega)$ is an almost \cy
manifold but \textit{not} a \cy manifold.

An $n$-dimensional submanifold $L$ of an almost \cy $n$-fold is
\textit{special Lagrangian} if
$$\omega|_L =0, \ \Imag(\Omega)|_L =0.$$
Since a \slg submanifold of a \cy manifold is also a \slg
submanifold in this generalized sense, this definition extends the
one given in \S \ref{cy}. From the condition that
$\Imag(\Omega)|_L =0$, it follows that $\Real(\Omega)$ is
a nowhere vanishing $n$-form on $L$. Hence each component of a
\slg submanifold $L$ is orientable, with a unique orientation in
which $\Real(\Omega)$ is positive.

For an almost \cy $n$-fold in which the function $f$ defined by
(\ref{acy_normalize}) is not identically 1, \slg submanifolds in
$M$ are neither calibrated nor minimal with respect to the
original metric $g$. Let $\tilde{g}$ be the conformally equivalent
metric $f^2 g$. Then $\Real(\Omega)$ \textit{is} a calibration on the
Riemannian manifold $(M,\tilde{g})$ and \slg submanifolds in
$(M,g,J,\Omega)$ are calibrated by this form.

\section{The differential and spectral geometry of SLG cones}
\label{slg:cones:section}

\subsection{\Slg cones in $\C^n$}
For any compact oriented embedded (but not necessarily connected)
submanifold $\Sigma
\subset S^{n-1}(1)\subset \R^n$ define the \textit{cone on
$\Sigma$},
$$ C(\Sigma) = \{ tx: t\in \R^{\ge 0}, x \in \Sigma \}.$$
A cone $C$ in $\R^n$ is $\textit{regular}$ if there exists
$\Sigma$ as above so that $C=C(\Sigma)$, in which case we call
$\Sigma$ the \textit{link} of the cone $C$. $C(\Sigma) - \{0\}$ is an
embedded smooth submanifold, but $C(\Sigma)$ has an isolated
singularity at $0$ unless $\Sigma$ is a totally geodesic sphere.

To characterize the links of regular \slg cones we need to
introduce some geometric structures on the unit sphere $S^{2n-1}$
in $\C^n$. As a convex hypersurface in a K\"ahler  manifold
$S^{2n-1}$ inherits a \textit{contact form}, that
is, a $1$-form $\gamma$ so that
\begin{equation}
\label{contact} \gamma \wedge d\gamma^{n-1} \neq 0.
\end{equation}
Let $X$ denote the Euler vector field $x\cdot\partial /\partial x$
on $\C^n$ and $\omega$ denote the standard symplectic form on
$\C^n$. Then the contact form on $S^{2n-1}$ is
$ \gamma = \iota_X\omega | _{S^{2n-1}}.$
Associated with $\gamma$ is the \textit{contact distribution}, the
hyperplane field $\ker{\gamma} \subset TS^{2n-1}$. The condition
(\ref{contact}) on $\gamma$ ensures that the distribution
$\ker{\gamma}$ is not integrable. The maximal dimensional integral
submanifolds (\textit{i.e.} submanifolds on which $\gamma$ restricts to
zero) of the distribution are $(n-1)$-dimensional and are called
\textit{Legendrian submanifolds}.
The relevance of Legendrian submanifolds of the sphere is seen in the following lemma.
\begin{lemma}
\label{link:cone}
Let $\Sigma$ be an $(n-1)$-dimensional submanifold of
$S^{2n-1}$. Then $C(\Sigma)-\{0\}$ is Lagrangian if and only if
$\Sigma$ is Legendrian.
\end{lemma}

For any $p$-form $\phi$ on $\R^n$ define the \textit{normal part}
of $\phi$ by
$ \phi_N = \iota_X \phi,$
where $X$ again denotes the Euler vector field on $\C^n$. In
particular, $\alpha_{\theta,N}$ denotes the normal part of the
\tslg calibration $\alpha_{\theta}$.
An oriented $(n-1)$-dimensional submanifold $\Sigma$ of
$S^{2n-1}$ is a \textit{\tsleg submanifold} if at each point of
$\Sigma$, $\alpha_{\theta,N}$ restricts to the volume form on
$\Sigma$.

\begin{prop}[\mbox{\cite[Prop 2.5]{haskins:slgcones}}]
\label{tsleg} A regular cone $C=C(\Sigma)$ in $\C^n$ is \tslg if
and only if $\Sigma$ is \tslegg.
\end{prop}

We also have the following Legendrian analogue of Proposition \ref{mlag}.
\begin{prop}[\mbox{\cite[Prop 2.6]{haskins:slgcones}}]
\label{minlag_tsleg} A connected oriented Legendrian submanifold
of $S^{2n-1}$ is minimal if and only if it is \tsleg for some
$\theta$.
\end{prop}


Let $\pi: S^{2n-1} \ra \CP^{n-1}$ be the Hopf fibration. When
$\CP^{n-1}$ is given the Fubini-Study metric $g_{\text{FS}}$ with
constant holomorphic sectional curvature 4, $\pi$ is a Riemannian
submersion with totally geodesic fibres. Hence if $M \subset
S^{2n-1}$ is minimal in $S^{2n-1}$ then $\pi(M)$ is minimal in
$(\CP^{n-1},g_{\text{FS}})$. Moreover, the standard contact structure
$\gamma$ on $S^{2n-1}$ is a connection compatible with the natural
metric on the Hopf bundle, and whose curvature is the Fubini-Study
\ka form $\omega$. Hence the associated horizontal distribution is
exactly the standard contact distribution on the sphere. In other
words, any Legendrian submanifold of $S^{2n-1}$ is horizontal for
$\pi$. Therefore, every minimal Legendrian immersion $f:M \ra
S^{2n-1}$ induces a minimal immersion
$\phi = \pi \circ f: M \ra \CP^{n-1},$
which in fact is Lagrangian \textit{i.e.} $\phi^*
\omega=0$. Conversely, given any minimal Lagrangian immersion
$\phi: M \ra \CP^{n-1}$ there is a horizontal lift of its
universal cover $\tilde{M}$ to the sphere, which is minimal
Legendrian. In fact, a given minimal Lagrangian immersion $\phi$
has an $S^1$ family of different possible lifts to the sphere.
Different choices of lift correspond to the lift being \tsleg for
different $\theta$.

\subsection{Symmetries of \slg $n$-folds and moment maps}

In symplectic geometry continuous
group actions which preserve $\omega$ are often induced by a collection of functions,
the \textit{moment map} of the action \cite[\S 5.2]{mcduff:salamon}. We describe the
situation in $\C^n$, where the situation is particularly simple, following
the presentations in \cite[\S 3]{joyce:slgcones2} and \cite[\S 4]{joyce:symmetries}.

Let $(M,J,g,\omega)$ be a \ka manifold, and let $\text{G}$ be a Lie group
acting smoothly on $M$ preserving both $J$ and $g$ (and hence $\omega$).
Then there is a natural linear map $\phi$ from the Lie algebra $\mathfrak{g}$ of $\text{G}$ to
the vector fields on $M$. Given $x\in \mathfrak{g}$, let $v=\phi(x)$ denote the
corresponding vector field on $M$. Since $\mathcal{L}_v\omega = 0$, it follows
from Cartan's formula that $\iota_v\omega$ is a closed $1$-form on $M$. If
$H^1(M,\R)=0$ then there exists a smooth function $\mu^x$ on $M$, unique up to a constant, such that
$d\mu^x = \iota_v\omega$. $\mu^x$ is called a moment map for $x$, or a Hamiltonian function
for the Hamiltonian vector field $v$.
We can collect all these functions $\mu^x$, $x\in \mathfrak{g}$, together to make a moment map for the whole action.
A smooth map $\mu:M \ra \mathfrak{g}^{*}$ is called a \textit{moment map} for the action of $\text{G}$ on $M$ if:
\begin{enumerate}
\item[(i)]
$\iota_{\phi(x)}\omega = <x,d\mu>$ for all $x\in \mathfrak{g}^{*}$ where $<,>$ is the natural pairing of $\mathfrak{g}$ and $\mathfrak{g}^*$
\item[(ii)]
$\mu:M \ra \mathfrak{g}^*$ is equivariant with respect to the $\text{G}$-action on $M$ and the coadjoint $\text{G}$-action on $\mathfrak{g}^*$.
\end{enumerate}
In general there are obstructions
to a symplectic $\text{G}$-action admitting a moment map.
The subsets $\mu^{-1}(c)$ are the \textit{level sets} of the moment map. The \textit{centre} $Z(\mathfrak{g}^*)$ is the subspace
of $\mathfrak{g}^*$ fixed by the coadjoint action of $\text{G}$. Property (ii) of $\mu$ implies
that a level set $\mu^{-1}(c)$ is $\text{G}$-invariant if and only if $c\in Z(\mathfrak{g}^*)$.

The group of automorphisms of $\C^n$ preserving $g$, $\omega$ and $\Omega$
is $\text{SU(n)} \ltimes \C^n$, where $\C^n$ acts by translations. The Lie algebra
$\mathfrak{su}(n) \ltimes \C^n$ acts on $\C^n$ by vector fields. Let $v$ be a vector field in
$\mathfrak{su}(n) \ltimes \C^n$. Since the action of $v$ preserves $\omega$ and $\C^n$ is simply connected,
there is a moment map for the vector field $v$, \textit{i.e.}
a function $\mu:\C^n \ra \R$, unique up to a constant, for which $\iota_v \omega = d\mu$.
Since the action of any $v \in \mathfrak{su}(n) \ltimes \C^n$ preserves the condition to be special Lagrangian,
the restriction to \slg submanifolds of the moment map of any $\mathfrak{su(n)} \ltimes \C^n$ vector field
enjoys some special properties.
\begin{lemma}[Joyce, \mbox{\cite[Lem 3.4]{joyce:slgcones2}}]
\label{harmonic:momentmap}
Let $\mu: \C^n \ra \R$ be a moment map for a vector field $v$ in $\mathfrak{su}(n) \ltimes \C^n$.
Then the restriction of $\mu$ to any \slg $n$-fold in $\C^n$ is a harmonic function
on $L$ with respect to the metric induced on $L$ by $\C^n$.
\end{lemma}
\noindent
That is any \slg $n$-fold $L$ in $\C^n$ automatically has certain distinguished harmonic functions.
This fact, was proven in a different way by Fu.
\begin{lemma}[\mbox{Fu, \cite[Thm 3.2]{fu}}]
\label{fu}
A function $f$ on $\C^n$ is harmonic on every special Lagrangian
submanifold in $\C^n$ if and only if
\begin{equation}
\label{fu:eqn}
d( \iota_{X_f}\Imag(\Omega)) = 0
\end{equation}
where $X_f = J\nabla f$ is the Hamiltonian vector field associated with $f$.
For $n\ge3$, $f$ satisfies (\ref{fu:eqn}) if and only if $f$ is a harmonic Hermitian quadratic.
\end{lemma}
A function is said to be a \textit{harmonic Hermitian quadratic} if it is of the form
\begin{equation}
\label{hermitian:harmonic}
f = c + \sum_{i=1}^n{(b_i z_i + \bar{b_i}\bar{z_i})} + \sum_{i,j=1}^n{a_{ij}z_i\bar{z_j}}
\end{equation}
for some $c\in \R, \ b_i, \ a_{ij} \in \C \textrm{\ with \ } a_{ij}=\bar{a}_{ji} \textrm{\ and\ } \sum_{i=1}^n{a_{ii}}=0.$
A harmonic Hermitian quadratic with $c=0$, $a_{ij}=0$ for all $i$ and $j$ corresponds to the moment
map of a translation.
A harmonic Hermitian quadratic with $c=0$, $b_i=0$ for all $i$ corresponds to the moment map
of the element $A=\sqrt{-1}\,(a_{ij})\in \mathfrak{su}(n)$.
Hence Lemma \ref{fu} implies Lemma \ref{harmonic:momentmap} and also gives a
gives a converse for $n\ge 3$. While Lemma \ref{harmonic:momentmap} gives a geometric
interpretation of Lemma \ref{fu}.


We can use Lemma \ref{harmonic:momentmap} or \ref{fu} together with the maximum
principle to conclude that a compact \slg submanifold with boundary
inherits all the symmetries of the boundary.
\begin{prop}
\label{slg:bdy:symmetry}
Let $L^n$ be a compact connected \slg submanifold of $\C^n$ with boundary $\Sigma$.
Let $\text{G}$ be the identity component of the subgroup of $\text{SU(n)} \ltimes \C^n$ which
preserves $\Sigma$. Then $\text{G}$ admits a moment map $\mu:\C^n\ra \mathfrak{g}^*$,
both $\Sigma$ and $L$ are contained in $\mu^{-1}(c)$ for some $c\in Z(\mathfrak{g}^*)$
and $\text{G}$ also preserves $L$.
\end{prop}
\begin{proof}
Since $L$ is Lagrangian, $\Sigma = \partial L$ is isotropic, \textit{i.e.} $\omega|_{\Sigma}=0$.
Let $\mathcal{O}$ be any orbit of $\text{G}$ in $\Sigma$. Then $\omega|_\mathcal{O}=0$, since $\mathcal{O} \subset \Sigma$.
It follows from results in \cite[\S 4]{joyce:symmetries} that $\text{G}$ admits a moment map
$\mu:\C^n \ra \mathfrak{g}^{*}$ and that $\Sigma \subset \mu^{-1}(c)$ for some $c\in Z(\mathfrak{g}^*)$.
Let $v$ be any vector field in $\mathfrak{g}$, then by Lemma \ref{harmonic:momentmap} the moment map of $v$
restricts to be a harmonic function on $L$. Since $\Sigma \subset \mu^{-1}(c)$, the moment map
of any vector field $v \in \mathfrak{g}$ must be constant on $\Sigma$ and hence by the Maximum Principle
also constant on $L$. Hence $L \subset \mu^{-1}(c)$ also.

For any $v \in \mathfrak{g}$, its moment map $\mu(v)$ is constant on $L$. Hence $\nabla \mu(v)$ is normal
to $L$ for any $v\in \mathfrak{g}$. Since $L$ is Lagrangian, $J$ maps $NL$ isomorphically to $TL$. Hence
$J \nabla \mu(v)$ is tangent to $L$ for any $v\in \mathfrak{g}$. It follows that $L$ is $\text{G}$-invariant.
\end{proof}
An asymptotically conical analogue of this result (Prop \ref{acslg:symmetry}) is given in Appendix
\ref{acslg:section} where we define asymptotically conical \slg submanifolds precisely.

\subsection{Harmonic functions on \slg cones and the stability index}
\label{harm:fn:cone}

When $L$ is a \slg cone $C$, there is an intimate relationship between homogeneous harmonic functions
on $C$ and eigenfunctions of $\Delta$ on the link $\Sigma$.
Given a regular \slg cone $C$ in $\C^n$ with link $\Sigma$, denote by $C'$, the incomplete but nonsingular manifold $C - \{0\}$.
Identify $C'$ with $(0,\infty) \times \Sigma$ via the map $\iota: (0,\infty) \times \Sigma \ra \C^n$ given by
$\iota(r,\sigma) = r\sigma$. Under this identification, the cone metric $dr^2 + r^2 g_{\Sigma}$, on $(0,\infty)\times \Sigma$,
is the metric $C'$ inherits from $\C^n$.
A function $u: C' \ra \R$ is \textit{homogeneous of order} $\alpha \in \R$, if $u(r,\sigma) \equiv r^{\alpha}v(\sigma)$
for some function $v:\Sigma \ra \R$ on the link.
A straightforward calculation using the definition of the cone metric proves the following.
\begin{lemma}
\label{harmonic:cone}
Let $u(r,\sigma) = r^{\alpha}v(\sigma)$ be a homogeneous function of order $\alpha$ on the cone
$C' \cong (0,\infty) \times \Sigma$, for some $v\in C^2(\Sigma)$. Then
$$ \Delta u (r,\sigma) = r^{\alpha-2}\left( \Delta_{\Sigma}v- \alpha(\alpha+n-2)v\right),$$
where $\Delta$ and $\Delta_{\Sigma}$ are the Laplacians on $C'$ and $\Sigma$ respectively. Hence,
$u$ is harmonic on $C'$ if and only if $v$ is an eigenfunction of $\Delta_{\Sigma}$
with eigenvalue $\alpha(\alpha+n-2)$.
\end{lemma}
\begin{remark}
\label{sign:convention}
We take our Laplacian $\Delta_{\Sigma}$ to be a nonnegative operator, and hence an eigenfunction $f$ of eigenvalue
$\lambda$ satisfies $\Delta_{\Sigma}f = \lambda f$ for some $\lambda\ge0$.
\end{remark}
\begin{remark}
\label{obvious:eigenvalues}
From Lemma \ref{harmonic:momentmap} we see that \textit{any} \slg cone $C'$  has nontrivial homogeneous
harmonic functions of order $1$ and $2$.
Using Lemma \ref{harmonic:cone} we see this is equivalent to $\lambda=n-1$
and $\lambda=2n$ being eigenvalues of $\Delta_{\Sigma}$.
\end{remark}
\begin{remark}
For any minimal Legendrian submanifold $\Sigma$ of $S^{2n-1}$, Remark \ref{obvious:eigenvalues} implies that
\begin{equation}
\label{first:eval:upbd}
\lambda_1(\Sigma) \le n-1.
\end{equation}
It is well-known that in fact (\ref{first:eval:upbd}) holds for any $(n-1)$-dimensional minimal submanifold of a unit sphere.
%
In \S \ref{legn:index:section} we will define the notion of a (Legendrian) stable minimal Legendrian submanifold of $S^{2n-1}$.
Any (Legendrian) stable minimal Legendrian submanifold $\Sigma$ will attain the equality in (\ref{first:eval:upbd}).
We would like to classify all (Legendrian) stable minimal Legendrian submanifolds. Theorem \ref{ct:is:only:stable:t2}
shows that there is a unique minimal Legendrian $2$-torus $\Sigma$ for which the equality $\lambda_1(\Sigma)=2$ holds.
\end{remark}

Let $C$ be a regular \slg cone in $\C^n$ ($n\ge 3$) with link $\Sigma$.
We need to introduce some notation for objects associated with $\text{Spec}(\Sigma)$.
We will follow the notation used in Joyce \cite{joyce:slgcones1,joyce:slgcones2,joyce:slgcones3,joyce:slgcones4,joyce:slgcones5}.
First define
\begin{equation}
\label{d:sigma}
\mathcal{D}_{\Sigma} = \left\{ \alpha \in \R: \alpha(\alpha+n-2) \textrm{\ is an eigenvalue of } \Delta_{\Sigma}\right\}.
\end{equation}
Lemma \ref{harmonic:cone} implies that $\alpha\in \dsig$ if and only if there exists a nonzero homogeneous
harmonic function of order $\alpha$ on $C'$.
The spectrum of $\Delta_{\Sigma}$ consists of a countably infinite discrete set of eigenvalues
$$0=\lambda_0< \lambda_1 \le \lambda_2\le \ldots \le \lambda_i \le \ldots$$
with corresponding orthonormal eigenfunctions
$$ \phi_0 = 1, \phi_1, \phi_2, \ldots , \phi_i, \ldots$$
Every eigenvalue $\lambda_i$ of $\Delta_{\Sigma}$
gives rise to exactly two points $\alpha_i^{\pm}$ in $\dsig$
$$\alpha_i^{\pm} = -\frac{1}{2}(n-2) \pm \sqrt{\frac{1}{4}(n-2)^2 + \lambda_i}.$$

For $\alpha \in \dsig$, define $m_{\Sigma}(\alpha)$
to be the \textit{dimension of the space of homogeneous harmonic functions of order $\alpha$ on $C$}.
Since $\lambda\ge 0$, we see that $\dsig \cap (2-n,0)$ is always empty.
Clearly $m_{\Sigma}(0)=m_{\Sigma}(2-n)=b^0(\Sigma)$, the number of components of $\Sigma$.

Let $N_{\Sigma}: \R \ra \Z$ be the unique function with the following properties: monotone increasing, upper semicontinuous
with a jump discontinuity of size $m_{\Sigma}(\alpha)$ at each $\alpha \in \dsig$, with value $0$
on the interval $(2-n,0)$. $N_{\Sigma}$ may be written in terms of $m_{\Sigma}$ as
\begin{equation}
\label{nsig:plus}
N_{\Sigma}(\beta) = \sum_{\alpha \in \dsig \cap [0,\beta]}{m_{\Sigma}(\alpha)} \ \textrm{\ if \ } \beta\ge 0
\end{equation}
or
\begin{equation}
\label{nsig:minus}
N_{\Sigma}(\beta) = -\sum_{\alpha \in \dsig \cap (\beta,0)}{m_{\Sigma}(\alpha)} \ \textrm{\ if \ } \beta< 0.
\end{equation}
Equivalently for $\beta\ge0$, $N_{\Sigma}(\beta)$ is equal to the number of eigenvalues (counting multiplicities)
of $\Delta_{\Sigma}$ in the range $[0,\beta(\beta+n-2)]$. Both $\dsig$ and $N_{\Sigma}$
play important roles in governing the analytic properties of $\Delta$ not just on cones
but also on spaces which are in some sense asymptotic to one or more \slg cone.

The \textit{stability index} $\sind(C)$ is defined in \cite[Defn 3.6]{joyce:slgcones2} to be
\begin{equation}
\label{s-index}
\sind(C) = N_{\Sigma}(2) - b^0(\Sigma) - 2n - (\dim{\text{SU(n)}} - \dim{\text{G}}),
\end{equation}
where $\text{G}$ is the subgroup of $\text{SU(n)}$ which leaves $C$ invariant.
Equivalently, the stability index of a \slg cone $C$ is the number of
eigenvalues that $\Delta_{\Sigma}$ has in the range $[0,2n]$ in \textit{excess}
of those eigenvalues which by Lemmas \ref{harmonic:momentmap} and \ref{harmonic:cone}
must be present for any \slg cone.
For any regular \slg cone with a singularity it is the case that
$$ m_{\Sigma}(0)= b^0(\Sigma), \quad m_{\Sigma}(1) \ge 2n \quad \textrm{and} \quad m_{\Sigma}(2) \ge \dim{\text{SU(n)}}-\dim{\text{G}}.$$
The first inequality will be proven in \S \ref{lindex}, and is equivalent to the geometric statement
that a minimal Legendrian submanifold in $S^{2n-1}$ is \textit{linearly full}, \textit{i.e.} is not
contained in any geodesic hypersphere.
For the second inequality see \cite[Prop 3.5]{joyce:slgcones2}. Together these facts
show that
\begin{equation}
\label{sind:nonnegative}
\textrm{s-ind(C)}\ge 0
\end{equation}
for any regular \slg cone $C$ which is not totally geodesic.

We will call a \slg cone $C$ \textit{strictly stable} if
\begin{equation}
\label{strictly:stable}
\sind(C)=0.
\end{equation}
Joyce calls such a cone stable. We will define the weaker notion of a
stable \slg cone in \S \ref{lindex}, where we explain the reason for this choice of terminology.
A strictly stable cone will then of course also be stable.

The importance of $\sind(C)$ is explained in Appendix \ref{conical:sing:section} and \ref{generic:singularities}.
One reason is that it appears as the dimension of an obstruction space to deforming a \slg submanifold $X$
with a conical singularity modelled on $C$. Hence the deformation theory simplifies
considerably if all the singularities of $X$ are modelled on strictly stable cones.
Another reason is that the most common isolated singularities should be modelled
on \slg cones with low stability index.


\subsection{Legendrian variations and the Legendrian index}
\label{legn:index:section}
In this section we recall some basic facts about small
Legendrian deformations of a Legendrian submanifold $\Sigma^n$ of any contact manifold $M^{2n+1}$
and define the Legendrian index of a compact minimal Legendrian submanifold of $S^{2n-1}$.

\subsubsection{The Legendrian neighbourhood theorem and Legendrian variations}
In symplectic geometry there is a well-known result due to Weinstein \cite{weinstein}, the Lagrangian neighbourhood theorem,
stating that some neighbourhood of a compact Lagrangian submanifold $L^n$ in any symplectic
manifold $(M^{2n},\omega)$ is symplectically equivalent to a certain standard model.
For compact Legendrian submanifolds $L^n$ of a contact manifold $(M^{2n+1},\gamma)$ there is an analogue of this result,
the Legendrian neighbourhood theorem, which seems to be less well-known.
First, we describe the relevant standard model, \textit{i.e.} the Legendrian analogue of the cotangent bundle.

The $m$-th jet space $J^m(L)$ of the manifold $L^n$ is roughly speaking the manifold which parametrizes
functions on $L$, together with all their derivatives up to order $m$. The $1$st jet space $J^1(L)$ is particularly
simple to describe. Topologically, it is equivalent to the rank $n+1$ vector bundle $\R \times T^*L$ over $L$, whose
total space is $(2n+1)$-dimensional.

There is  a canonical contact $1$-form $\gamma_0$ on $J^1(L)\cong \R \times T^*L$
\begin{equation}
\label{contact:jet}
\gamma_0 = dt - \pi^*\alpha
\end{equation}
where $t$ is the coordinate on $\R$, $\pi: J^1(L) \ra T^*L$ is the natural projection and $\alpha$ is
the tautological $1$-form on $T^*L$. Moreover, $L$ sits as a Legendrian submanifold of $(J^1(L),\gamma_0)$.
More generally, given any smooth function $f$ on $L$ there is an associated section of $J^1(L)$
$$ f \rightsquigarrow (f(x), x, df(x))$$
whose image in $J^1(L)$ is a Legendrian submanifold of $(J^1(L), \gamma_0)$.

\begin{prop}[The Legendrian Neighbourhood Theorem]
\label{legn:nhd}
Let $L^n$ be a compact Legendrian submanifold of any contact manifold $(M^{2n+1},\gamma)$.
Then there exists an open neighbourhood of $L$ which is contact equivalent
to an open neighbourhood of the zero section in $(J^1(L),\gamma_0)$.
Moreover, a section of $J^1(L)$ is Legendrian if and only if it is the 1-jet associated to a function on $L$.
\end{prop}
For further information on basic notions in contact geometry see \cite[\S 3.4]{mcduff:salamon}. Note that unlike the Lagrangian case,
every nearby Legendrian submanifold has an associated (globally defined) function.
This makes it very plausible that the geometry of Legendrian submanifolds is
closely related to function theory on these submanifolds.

\subsubsection{Legendrian index: definition and basic properties}
\label{lindex}
Let $\Sigma^{n-1}$ be a compact minimal Legendrian submanifold of $S^{2n-1}$.
Define the \textit{Legendrian index of} $\Sigma$, denoted $\lind(\Sigma)$, to be the number of eigenvalues (counted with multiplicity) of $\Delta: C^{\infty}(\Sigma) \ra C^{\infty}(\Sigma)$
in the interval $(0,2n)$.
The Legendrian index of a regular \slg cone $C$ means the Legendrian index of its link $\Sigma$.

\begin{remark}There is a geometric interpretation of the Legendrian index. It is the index of the operator corresponding to the second
variation of volume of $\Sigma$ when restricted to
Legendrian variations of $\Sigma$, whence the name. Since we do not use this fact we shall not give the proof.
Using the fact that $M$ is Sasakian-Einstein if and only if the metric cone $C(M)$ is Ricci-flat K\"ahler \cite[\S 1.1]{boyer:galicki},
the same geometric interpretation of the Legendrian index stills holds
if we replace $S^{2n-1}$ with any other compact Sasakian-Einstein manifold $M$.
\end{remark}
We want to prove that the Legendrian index of a \slg cone with isolated
singularity in $\C^n$ is at least $2n$. This will follow from a slightly
more general result, Lemma \ref{coiso:index}.
First, recall some basic definitions from symplectic
geometry. The symplectic complement $W_{\omega}^{\perp}$
of a subspace $W$ of a symplectic vector space $(V^{2n},\omega)$,
is defined to be
$ \{ v\in V |\  \omega(v,w)=0 \ \forall\  w \in W \}.$
Since $\omega$ is nondegenerate $\dim{W} + \dim{W_{\omega}^{\perp}}=\dim{V}$.
$W$ is \textit{coisotropic} if $W \supseteq W_{\omega}^{\perp}$, and
\textit{isotropic} if $W \subseteq W_{\omega}^{\perp}$. If $W$ is
isotropic then $W_{\omega}^{\perp}$ is coisotropic.
If $W$ is coisotropic and $\dim{W}=n$, then
$W=W_{\omega}^{\perp}$ and $W$ is a \textit{Lagrangian} subspace of $V$.

If $V=\C^n$ with the standard complex, symplectic and Euclidean structures
$J,\ \omega,\  <,>$ we can interpret the symplectic complement in terms
of the usual orthogonal complement and the action of $J$. In fact,  $W_{\omega}^{\perp} = (JW)^{\perp}$.
Hence $W$ is coisotropic if and only if $(JW)^{\perp} \subseteq W$,
or equivalently $W\supseteq J(W^{\perp})$. More generally,
if $X$ is a coisotropic submanifold of a \ka manifold $(M,\omega,g,J)$
this last inclusion implies that $J$
maps sections of $NX$ into sections of $TX$.
If $X$ is an isotropic submanifold of a \ka manifold $(M,\omega,g,J)$, then
instead $J$ maps sections of $TX$ into sections of $NX$.
In the case that $X$ is Lagrangian, $J$ is an isometry between
$TX$ and $NX$.

\begin{lemma}
\label{coiso:index}
Let $C$ be a coisotropic cone of real dimension $m$ (for some $m\ge n$) with link $\Sigma$ in $\C^n$
which has an isolated singularity at $0$. Then $C$ is linearly full,
\textit{i.e.} $C$ is not contained in any hyperplane. In particular, if $C$ is also minimal then
$\lambda=m-1$ is an eigenvalue  of $\Delta_{\Sigma}$ of multiplicity at least $2n$.
\end{lemma}

\begin{proof}
Define $C^{'}=C\setminus\{0\}$.
Suppose that $C$ is contained in some hyperplane.
Then there exists some $a\in \C^n$ so that the linear function $l_a(x)=<x,a>$
restricts to zero on $C$. It follows that $a$ is normal to $C^{'}$ at any
point of $C^{'}$. Since $C^{'}$ is a coisotropic submanifold this implies that $Ja$ is tangent
to $C^{'}$ at any point in  $C^{'}$. It follows that $C^{'}$
is invariant under translations in the direction $Ja$.
This contradicts our assumption that $C$ has an isolated singularity at $0$.

For any minimal $m$-dimensional cone $C=C(\Sigma)$ in $\C^n$, it is well-known that
the restriction to the link $\Sigma$ of any linear function on $\C^n$
is an eigenfunction of $\Delta_{\Sigma}$ with eigenvalue equal to
$\dim{\Sigma}=m-1$. The result now follows from the first part of the lemma.
\end{proof}
\noindent
Since any hypersurface of $\C^n$ is automatically coisotropic, Corollary \ref{coiso:index}
implies that any minimal hypersurface of $S^{2n-1}$ which is not a round hypersphere is automatically linearly full.

\begin{corollary}[\mbox{\cite[Prop 3.5]{joyce:slgcones2}}]
\label{lb:index}
The Legendrian index of any \slg cone in $\C^n$ with an isolated singularity
is at least $2n$.
\end{corollary}

\begin{proof}
Clearly a \slg cone $C$ is a minimal coisotropic cone in $\C^n$, and
hence $C$ is linearly full by the previous corollary.
Hence $n-1$ is an eigenvalue of $\Delta_{\Sigma}$ with multiplicity at least $2n$.
The result is now immediate from the definition of Legendrian index.
\end{proof}


We call a \slg cone $C$ in $\C^n$ with link $\Sigma$, \textit{Legendrian stable} (or often just stable) if
$$\lind(C)=2n,$$
\textit{Jacobi integrable} if every eigenfunction of $\Delta_{\Sigma}$ with eigenvalue $\lambda=2n$
arises as the variation vector field of some $1$-parameter family of minimal Legendrian submanifolds
and \textit{rigid} if
$$\dim{\ker{(\Delta-2n)}}=\dim{\text{SU(n)}}-\dim{\text{G}}$$
where $\text{G}$ is the subgroup of $\text{SU(n)}$ which leaves $C$ invariant.
A rigid cone is automatically Jacobi integrable.

We call $C$ \textit{(Legendrian) strictly stable} if it is both Legendrian stable and rigid.
This is equivalent to the definition given in \S \ref{harm:fn:cone}
that $C$ is (Legendrian) strictly stable if $\sind(C)=0$.
This terminology seems preferable to Joyce's use of stable for $\sind(C)=0$ since it is more
consistent with the use of stable and strictly stable in other minimal submanifold settings.

In terms of eigenvalues of the Laplacian on functions on $\Sigma$,
(Legendrian) stability implies
$$\lambda_1=n-1 \textrm{\ (with multiplicity $2n$), \quad and\quad }\lambda_{2n+1}=2n.$$
Strict stability implies further
that $\lambda_{2n+1}=2n$ has the minimum allowable multiplicity.
In \S \ref{dimension:three} we prove that if the link $\Sigma$
is a $2$-torus, and the cone $C(\Sigma)$ is (Legendrian) stable, then $\Sigma$ must be
the Clifford torus described below in Example \ref{clifford:torus}.


There are relatively few Riemannian manifolds $(M,g)$ for which one can compute $\textrm{Spec}(\Delta)$
explicitly. Typically $(M,g)$ must be very symmetric for this to be possible.  Below
we compute $\textrm{Spec}(\Delta)$ and hence $\lind$ and $\sind$ for the simplest minimal Legendrian submanifolds in $\C^n$:
the totally geodesic $S^{n-1}$ and the Clifford torus $T^{n-1}$.
\begin{example}[The round sphere]
Let $\Sigma$ be the totally geodesic minimal Legendrian sphere $S^{n-1} \subset S^{2n-1}$ defined by
$S^{n-1}=\{ z\in S^{2n-1}: \Imag(z)=0\}$.
In this case it is well-known that all eigenfunctions of $\Delta_{\Sigma}$ arise from the
restrictions of harmonic homogeneous polynomials on $\R^n$ \cite{gallot}.
Since $\R^n$ with the Euclidean metric is the metric cone over $S^{n-1}$, Lemma \ref{harmonic:cone}
shows that the restriction of a harmonic homogeneous polynomial of degree $k$ on $\R^n$ is an eigenfunction
of $\Delta_{\Sigma}$ with eigenvalue $\lambda=k(k+n-2)$.
It follows immediately from (\ref{d:sigma}) that
$$\dsig = \Z -\{-1, -2, \ldots , 3-n\}.$$
The dimension of the space of harmonic homogeneous polynomials of degree $k$ is also known,
\textit{e.g.} \cite[p. 1059]{cheng:li:yau}.
It is also immediate that
$$\lind(\Sigma)=n,\ \sind(\Sigma)=-n \textrm{\ and \ }
m_{\Sigma}(2)=\frac{1}{2}(n-1)(n+2)= \dim{\text{SU(n)}}-\dim{\text{SO(n)}}.$$
The third equality implies that $\Sigma$ is rigid and hence also Jacobi integrable.
The first two equalities imply that $\Sigma$ violates inequality
(\ref{sind:nonnegative}) and the inequality in Corollary \ref{lb:index}. But since
$C(\Sigma) = \R^n$ is smooth everywhere, $\Sigma$ does not
satisfy the condition that $C=C(\Sigma)$ has a singularity at the origin.
\end{example}

\begin{example}[The Clifford torus]
\label{clifford:torus}
Let $\Sigma$ be the torus $T^{n-1}\subset S^{2n-1}$ defined by
$$ T^{n-1}= \left\{ z\in S^{2n-1}: |z_1|^2= \ldots = |z_n|^2=\frac{1}{n}, \ \textrm{arg}(z_1\ldots z_n)=0\right\}.$$
Then $\Sigma$ is a minimal Legendrian submanifold invariant under the group $T^{n-1}$ of diagonal matrices in $\text{SU(n)}$.
The metric on $\Sigma$ induced by its embedding in $S^{2n-1}$ is flat. For any flat torus one can explicitly compute
$\textrm{Spec}(\Delta)$ as follows (see \cite[p. 200]{gallot} for a proof).

Any flat $n$-torus can be realized as the quotient $\R^n/\Gamma$
of Euclidean $n$-space by some lattice $\Gamma$. Define the dual lattice $\Gamma^{*}$ by
$$ \Gamma^* = \{ x\in \R^n : \ <x,y> \in \Z \textrm{\ \ for all \ } y\in \Gamma\}.$$
Then the spectrum of $\Delta$ on $\R^n/\Gamma$ is given by
$$ \textrm{Spec}(\R^n/\Gamma) = \{ 4\pi^2|x|^2: \ x\in \Gamma^*\}.$$
Hence to compute $\textrm{Spec}(\Delta)$ it is enough to find the dual lattice $\Gamma^*$ corresponding to $\Sigma$.

$\Sigma$ is a subset of an $n$-torus $T^{n}=\left\{ z\in S^{2n-1}: |z_1|^2= \ldots = |z_n|^2=\frac{1}{n}\right\}$, which
is itself a flat minimal torus in $S^{2n-1}$. Since $T^n$ is the product of $n$ circles of radius $\frac{1}{\sqrt{n}}$,
the induced flat metric on $T^n$ corresponds to the lattice $\Gamma_0$ spanned by $E_1=\frac{2\pi}{\sqrt{n}}e_1, \ldots , E_n=\frac{2\pi}{\sqrt{n}}e_n$ where
$e_1, \ldots, e_n$ is the standard orthonormal basis of $\R^n$.
The dual lattice $\Gamma^*_0$ is spanned by $F_1=\frac{\sqrt{n}}{2\pi}e_1, \ldots , F_n=\frac{\sqrt{n}}{2\pi}e_n$.
$\Sigma$ now corresponds to the sublattice $\Gamma$ of $\Gamma_0$
defined by
$$\Gamma= \left\{\sum_{i=1}^{n}{k_i E_i}\in \Gamma_0: k_i \in \Z, \ \sum_{i=1}^{n}{k_i}=0 \right\}.$$
A basis for $\Gamma^*$ is given by $V_1, \ldots, V_{n-1}$ where
$$V_i = \frac{1}{2\pi\sqrt{n}}\left(ne_i -\sum_{j=1}^{n-1}{e_j}\right).$$
The $V_i$ satisfy
$$<V_i, V_i>=\frac{n-1}{4\pi^2} \quad \textrm{and\ } <V_i,V_j> = -\frac{1}{4\pi^2} \quad \textrm{if\ } i\neq j,$$
and so
$$\left|\sum_{i=1}^{n-1}{k_i V_i}\right|^2 = \frac{1}{4\pi^2}\left( n\sum_{i=1}^{n-1}{k_i^2} - \sum_{i,j=1}^{n-1}{k_ik_j}\right).$$
Hence we have
\begin{equation}
\label{ct:spectrum}
\textrm{Spec}({\Delta_{T^{n-1}}}) = \left\{n\sum_{i=1}^{n-1}{k_i^2} - \sum_{i,j=1}^{n-1}{k_ik_j} : \ (k_1, \ldots , k_{n-1})\in \Z^{n-1} \right\}.
\end{equation}
In principle one can simply use (\ref{ct:spectrum}) and a computer to find any finite part of the spectrum (including the multiplicities
of the eigenvalues). Since the number of lattice points close to the origin in an $n$-dimensional lattice grows exponentially with $n$,
a naive approach to computing $\textrm{Spec}(\Delta_{T^{n-1}})$ will fail to be effective when $n$ is large.
However, since we are most interested in eigenvalues in the range $[0,2n]$ one can easily do the computations for up to $n=12$.
From these computations one finds the following results (which in fact one can show hold for all $n$).
\begin{prop}[\mbox{\cite[\S 3.2]{joyce:slgcones2}}]
\label{ct:stability}
For all $n$, $\lambda_1(T^{n-1}) = n-1$.\\
For $n=3$, the Clifford torus $T^2$ is strictly stable, and hence both Legendrian stable and rigid.\\
For $n>3$, the Clifford torus $T^{n-1}$ is not Legendrian stable and hence not strictly stable.\\
For $n>3$, the Clifford torus $T^{n-1}$ is rigid except when $n=8$ or $n=9$.
\end{prop}
By Lemma \ref{harmonic:cone} the previous proposition implies: there are no extra sublinear homogeneous harmonic functions on $C(T^{n-1})$ for any $n$,
there are extra subquadratic homogeneous harmonic functions on $C(T^{n-1})$ for $n>3$, and there are extra quadratic harmonic functions when $n=8$ or $9$.

For a table of $N_{\Sigma}(2)$, $m_{\Sigma}(2)$ and $\sind(\Sigma)$ for $n\le 12$ see  \cite[\S 3.2, Table 1]{joyce:slgcones2}.
In fact, the table shows that for $n\le 9$ the behaviour of the eigenvalues in the range $[0,2n]$ is somewhat erratic, but settles down for $n\ge 10$.
Below we give a proof that for $n>3$ the Clifford torus is not Legendrian stable and
offer some explanation why the behaviour for $n\le 9$ is different from higher values of $n$ and
why the Clifford torus fails to be rigid for $n=8,9$.

For $m\in \N$ and $c\in \R$, define $c_\mathbf{m}$ to be the vector $(c, \ldots , c)\in \R^m$.
Given $p$, $q\in \N$ with $p+q \le n-1$ define $\mathbf{k}_{p,q} = (1_{\mathbf{p}}, -1_{\mathbf{q}}, 0_{\mathbf{n-p-q-1}})\in \Z^{n-1}$.
Denote by $\lambda_{p,q}= (p+q) n - (p-q)^2$ the eigenvalue of $\Delta_{T^{n-1}}$ corresponding to $\mathbf{k}_{p,q}$.
The erratic behaviour of the eigenvalues in $[0,2n]$ for $n\le9$ can be understood by looking at which of the eigenvalues
$$\lambda_{p,0}=p(n-p) \textrm{\quad for\ } p\in \N \cap \left[1,\frac{n}{2}\right]$$
(by the symmetry $p\ra n-p$ it is enough to consider only $p\le \frac{n}{2}$) are in the range $[0,2n]$.

$\lambda_{1,0}=n-1$ and $\lambda_{2,0}=2(n-2)$
are in the interval $[0,2n]$ for any $n$ and are distinct provided $n\neq3$.
Since for $n>3$ we have an eigenvalue $2n-4$ which lies in the interval $(n-1,2n)$
it follows immediately that we have $\lind(T^{n-1}) > 2n$,
\textit{i.e.} for $n>3$ the Clifford torus is not Legendrian stable, and hence also not strictly stable.

For other values of $p$, whether $\lambda_{p,0}\in [0,2n]$ depends on both $p$ and the dimension $n$.
For $n\le 8$ the maximum value of $p(n-p)$ is less than or equal to $2n$.
That is, for $n\le 8$, $\lambda_{p,0} \in [0,2n]$ for any $p\in \N \cap [1,\frac{n}{2}]$.
For $n>8$, by comparing $p(n-p)$ to $p(8-p)$ it follows that $\lambda_{p,0} \in [0,2n]$ implies $p< 4$.
Since we already know $\lambda_{1,0}, \lambda_{2,0}\in [0,2n]$ for any $n$, we need only deal with $p=3$.
But $\lambda_{3,0}\le 2n$ implies $n\le 9$ with equality if and only if $n=9$.
Hence the eigenvalues in the range $[0,2n]$ of the form $\lambda_{p,0}$ are:\\
\begin{center}
$\begin{array}{rr}
\lambda_{1,0} &  \textrm{for \ }n=3, \\
\lambda_{1,0}<\lambda_{2,0} & \textrm{for \ }n=4,5,\ n\ge 10 \\
\lambda_{1,0}<\lambda_{2,0}<\lambda_{3,0} & \textrm{for\ } n=6,7 \textrm{\ and\ } 9\\
\lambda_{1,0}<\lambda_{2,0}<\lambda_{3,0}<\lambda_{4,0} &  \textrm{for\ } n=8.
\end{array}
$
\end{center}
Note that for $n=8$, $\lambda_{4,0}$ coincides with $2n$, and for $n=9$, $\lambda_{3,0}$ coincides with $2n$.
These account for the extra $\lambda=2n$ eigenfunctions which in dimensions $8$ and $9$ cause the Clifford torus not to be rigid.
\end{example}
In \S \ref{dimension:three} and \S \ref{s1:invariant:tori} we will find methods to estimate
$\lind(C)$ and $\sind(C)$ for more complicated \slg cones $C$ which do
not require us to compute the whole spectrum of $\Delta$.

\subsection{Comparison with Hamiltonian index}

Oh gave the formula for the second variation operator $\mathcal{J}$
of a compact minimal Lagrangian submanifold $L$ of a \ke manifold $M$ thought of as an operator
on $1$-forms on $L$ \cite{oh}.
When $M=\CP^{n-1}$ with the metric induced by the Hopf fibration of the unit sphere $S^{2n-1}$,
his formula gives $\mathcal{J}=\Delta_L - 2n$.
Since we took $\Delta$ to be a nonnegative operator (Rmk \ref{sign:convention}),
$\mathcal{J}$ is bounded below but in general has negative eigenvalues.
Accordingly, for $\CP^{n-1}$ the index of $L$ as a minimal submanifold, denoted $\textrm{ind}(L)$, is
equal to the number of eigenvalues (counting multiplicity) in the
range $[0,2n)$ of the Laplacian acting on $1$-forms. In
particular, Hodge theory implies that $\textrm{ind(L)}\ge b^1(L)$.

We may also consider $\mathcal{J}$ restricted to
Lagrangian variations, and use this to define the \textit{Lagrangian
index}, denoted $\textrm{Lag-ind}$.
By the Lagrangian
neighbourhood theorem we need to look at variations corresponding
to closed $1$-forms. Clearly harmonic $1$-forms on $L$
still contribute to the Lagrangian index. By the Hodge
decomposition, any $1$-form $\alpha$ orthogonal to the harmonic $1$-forms and satisfying
$\Delta \alpha =\lambda \alpha$, for nonzero $\lambda$ is a sum of an exact and a coexact piece.
The coexact $1$-forms do not correspond to Lagrangian variations.
When $\alpha$ is exact, \textit{i.e.} $\alpha=df$ for some function $f$ on $L$, the
corresponding variation of $L$ is called a \textit{Hamiltonian
variation}. We can accordingly define a notion of the
\textit{Hamiltonian index} of $L$, denoted $\textrm{h-ind}(L)$.
$L$ is said to be \textit{Hamiltonian stable} if its Hamiltonian index is equal to
zero.

Since $ \Delta df= \lambda df \iff \Delta f = \lambda f$,
it follows that the Hamiltonian index of $L$ in $\CP^{n-1}$ is equal
to the number of eigenvalues (with multiplicity) of $\Delta: C^\infty(L) \ra  C^{\infty}(L)$
contained in the open interval $(0,2n)$. Hence we have $\textrm{Lag-ind}(L)=b^1(L)+ \textrm{h-ind}(L)$.
In the special case where $L$ is a
surface then the coexact $1$-forms also come from eigenfunctions
of $\Delta$ acting on $C^{\infty}(L)$. Hence when $L$ is a surface we have $\textrm{ind}(L)=b^1(L) +2 \ \textrm{h-ind}(L).$

Since on its horizontal space the Hopf projection $\pi: S^{2n-1}\ra \CP^{n-1}$ is an isometry,
it follows that the Legendrian index of a minimal Legendrian submanifold
$\Sigma$
in $S^{2n-1}$ is equal to the Hamiltonian index of the
corresponding minimal Lagrangian submanifold $\pi(\Sigma)$ in $\CP^{n-1}$.
By Remark \ref{obvious:eigenvalues}, $\lambda=n-1$ is always an eigenvalue of any minimal Legendrian
submanifold $\Sigma$ of $S^{2n-1}$. Hence the corresponding minimal Lagrangian
submanifold $\pi(\Sigma)$ of $\CP^{n-1}$ is \textit{never} Hamiltonian stable.
In particular,
the image in $\CP^{n-1}$ of the Clifford torus in $S^{2n-1}$ is not
Hamiltonian stable. On the other hand Oh \cite[Thm 5.4]{oh} showed that what he called the Clifford
torus $T^{n-1}$ in $\CP^{n-1}$ is Hamiltonian stable for all $n$.

To reconcile these two facts we notice that the metric induced on
the standard Clifford torus in $S^{2n-1}$ is a $\Z_{n}$-cover of
the metric on the standard Clifford torus in $\CP^{n-1}$, \textit{i.e.}
${T_{\CP^{n-1}}} = T_{S^{2n-1}}/{\Z_n}$.
$\Z_{n}$ arises because it is the centre of $\text{SU(n)}$.
Hence the
eigenfunctions on the complex projective Clifford torus are
exactly the eigenfunctions on the spherical Clifford torus which
are $\Z_{n}$-invariant.
In \S \ref{dimension:three} we will characterize
the Clifford tori both in $\CP^2$ and in $S^5$ as the only Hamiltonian stable/Legendrian stable
minimal Lagrangian/Legendrian $2$-torus.

\section{Special Lagrangian cones in dimension three}
\label{dimension:three}
The theory of \slg cones in $\C^3$ has a distinctive flavour.
There are two principal reasons for this.
First, the link $\Sigma$ of such a cone is a surface and the geometry of $\Delta$ on
a surface has special features.
Second, if $\Sigma$ is homeomorphic to a $2$-torus then powerful methods from
integrable systems can be applied, \textit{e.g.}  \cite{mcintosh:slg,mcintosh,mcintosh:nonisotropic,mcintosh:remarks}.
In principle, these methods allow the construction
of \textit{all} \slg $T^2$-cones. In practice, it has proven difficult to extract
useful differential geometric/analytic information from the integrable systems constructions.
In \S \ref{int:systems} we prove two results in this direction.
The main results of this section are the following.

First, in Theorem \ref{ct:is:only:stable:t2} we prove that the cone on the Clifford torus
in $S^5$ is the unique \textit{stable} (and hence also strictly stable) \slg $T^2$-cone.
This answers a problem posed by  Joyce
\cite[\S 10.3]{joyce:lectures},
and proves that the Clifford torus is the only \slg $T^2$-cone
which can appear in an index $1$ singularity (Corollary \ref{index:1:classification}).
In Theorem \ref{ct:only:hstable:t2} we prove the analogous result for the Clifford torus in $\CP^2$, replacing the notion
of (Legendrian) stability by Hamiltonian stability.

Second, in Theorem \ref{specgenus:bds:index} for any minimal Legendrian $2$-torus we give
an explicit lower bound for the \textit{Legendrian index} and for the \textit{stability index}, in terms of a
fundamental invariant of the integrable systems side, the \textit{spectral curve genus}. These bounds imply
that minimal Legendrian $2$-tori with large spectral curve genus have large Legendrian and stability indices.

Third, in Theorem \ref{area:spectral:genus} we give explicit lower bounds for the area of a minimal Legendrian $2$-torus in terms of the Legendrian index,
the stability index or the spectral curve genus. This result gives a special Lagrangian analogue of
results proved in \cite[Thm 6.7]{pedit:pinkall} for minimal tori in $S^3$. \cite{pedit:pinkall} makes essential
use of ``quaternionic holomorphic''  methods which do not appear applicable to \slg geometry.
Our results employ completely different methods, namely results from the spectral theory of $\Delta$
on a compact surface.

In \S \ref{spectral:geometry} we review results from the spectral theory of $\Delta$
on a compact Riemann surface. In \S \ref{conformal:area:section} we use a result
on minimal immersions by first eigenfunctions to prove Theorem \ref{ct:is:only:stable:t2}.
The Hamiltonian analogue of this result, Theorem \ref{ct:only:hstable:t2}, is proved using results from \S \ref{spectral:geometry}.
In \S \ref{int:systems} we review the construction of \slg $T^2$-cones in terms of spectral curve data.
By combining this construction with results of \S \ref{spectral:geometry} we prove Theorems \ref{specgenus:bds:index} and \ref{area:spectral:genus}.

\subsection{The geometry of the Laplacian on surfaces}
\label{spectral:geometry}
We recall some facts from the spectral geometry of surfaces.
In a number of respects the geometry of the Laplacian $\Delta: C^{\infty}(\Sigma) \ra C^{\infty}(\Sigma)$
on a compact surface $\Sigma$ is more rigid than in higher dimensions.

Let $M$ be a compact manifold endowed with any Riemannian metric $h$.
Let $\Delta_h: C^{\infty}(M) \ra C^{\infty}(M)$ be the Laplace-Beltrami operator
with respect to the metric $h$. $\Delta_h$ is a second order self-adjoint elliptic operator
when acting on the Sobolev space $L^2_1(M)$. Its spectrum consists of a countably infinite set
of discrete eigenvalues
$\textrm{Spec}(\Delta) = \{0=\lambda_0 < \lambda_1 \le \ldots \le \lambda_i \le \ldots \}$
each of finite multiplicity.

When $M=\Sigma^2$ an orientable surface of genus $g$, Yang-Yau \cite{yang:yau} proved that
\begin{equation}
\label{yang:yau}
\lambda_1 \Area(\Sigma) \le 8\pi (1+g)
\end{equation}
holds for any metric $h$ on $\Sigma$. The natural higher-dimensional analogue of (\ref{yang:yau}) turns out to
be false \cite{urakawa}.
Several authors have proved
upper bounds for the multiplicity of the $i$-th eigenvalue in terms of $i$ and the genus $g$ of $\Sigma$.

\begin{theorem}[\cite{nadir}]
\label{besson:thm}
Let $(\Sigma^2,h)$ be a compact orientable Riemannian surface of genus $g$, and let $m_i$ be the multiplicity of the $i$-th
eigenvalue $\lambda_i$ of $\Delta_{\Sigma}$ acting on functions.
If $g=0$, then $m_i \le 2i+1$.
If $g=1$, then $m_i \le 2i+4$.
If $g>1$, then $m_i \le 4g + 2i -1$.
\end{theorem}

\begin{remark}
The first result of this kind was proved by Cheng \cite{cheng}, and was subsequently improved by
Besson \cite{besson} and then Nadirashvili \cite{nadir}.
\end{remark}
\begin{remark}
\label{mult:growth:rmk}
The bound in Theorem \ref{besson:thm} is linear in the index $i$ and the genus $g$. However, the multiplicity
grows like at most $\sqrt{i}$ and $\sqrt{g}$ in all known examples. It is an open problem to
determine the correct rate of growth of $m_i$ with respect to $i$ and $g$.
\end{remark}
\begin{remark}
In contrast to Theorem \ref{besson:thm}, Colin de Verdiere  \cite{colindv} proved that on any compact $n$-manifold
with $n\ge 3$, one can prescribe any finite part of the spectrum of $\Delta_h$ by an appropriate choice of metric $h$.
\end{remark}

Let $f$ be an eigenfunction of $\Delta$, and define the
\textit{nodal set of $f$} to be $f^{-1}(0)$.
The structure of the nodal set of any eigenfunction $f$ of $\Delta$ on a surface $\Sigma$ is well understood.
It consists of the union of smooth simple curves (the \textit{nodal lines}) and isolated singular points
where finitely many nodal lines intersect. Near each singular point $f^{-1}(0)$ is homeomorphic
to the zero set of a harmonic homogeneous polynomial on $\R^2$. Hence at each singular point $p$
the curves form an ``equiangular configuration'', and the number of curves that intersect at $p$
is determined by the vanishing order of $f$ at $p$.
In higher dimensions $f^{-1}(0)$ may have a more complicated structure \cite{hardt:simon}.

The set $M - f^{-1}(0)$ is in general not connected. Its connected components are called the
\textit{nodal domains} of the eigenfunction $f$. Because $f$ has
the same sign on each nodal domain and vanishes on its boundary,
the restriction of $f$ to each nodal domain is a first
eigenfunction for the Dirichlet problem on that nodal domain.

\begin{theorem}[Courant Nodal Domain Theorem, \cite{berard:meyer,cheng}]
\label{courant:nd}
Let $(M,h)$ be a compact Riemannian manifold
with $\partial M = \emptyset$, and $\lambda_i$ be the $i$th
eigenvalue of $\Delta_h$ acting on $C^\infty(M)$. Then
\begin{enumerate}
\item[(i)] $\lambda_1 > 0$ and any eigenfunction with eigenvalue $\lambda_1$ has exactly $2$ nodal domains.
\item[(ii)]
 the number of nodal domains of any eigenfunction with eigenvalue $\lambda_i$ is at most $i+1$.
\end{enumerate}
\end{theorem}
\no
That is,
the number of nodal domains of an
eigenfunction gives us  a \textit{lower bound} on the index
of the corresponding eigenvalue. In \S \ref{s1inv:indexbounds}
we will use Theorem \ref{courant:nd} to prove
lower bounds for the Legendrian and stability indices
of $S^1$-invariant minimal Legendrian tori in $S^5$.

Another powerful tool to study the spectral geometry of $\Delta$ on any Riemannian manifold $(M,g)$ is the \textit{heat kernel}.
The heat kernel $H(x,y,t) \in C^{\infty}(M\times M\times \R^{+})$ is the fundamental solution to the heat equation
$$ \frac{\partial u}{\partial t} = \Delta_g u.$$
If $(M,g)$ is complete, then the heat kernel $H(x,y,t)$ always exists \cite[p94]{schoen:yau}.
If $M$ is compact then \cite[p103]{schoen:yau}
$$H(x,y,t) = \sum_{i=0}^{\infty}{e^{-\lambda_i t}\phi_i(x)\phi_i(y)}$$
where $\{\phi_i\}_{i=0}^{\infty}$ is an orthonormal basis of $L^2(M)$ consisting of eigenfunctions
of $\Delta_g$ with eigenvalues $\{\lambda_i\}_{i=0}^{\infty}$. In particular, the trace of the heat kernel satisfies
$$ \int_M{ H(x,x,t) dx} = \sum_{i=0}^{\infty}{e^{-\lambda_i t}}.$$

Cheng-Li-Yau \cite{cheng:li:yau} proved comparison results for the heat kernel of a minimal submanifold
immersed in any of the simply connected space forms $\Hyp^n$, $\R^n$ or $S^n$. From these results one can prove the following.
\begin{theorem}[\mbox{\cite[Thms 5,6]{cheng:li:yau}}]
\label{cheng:li:yau:thm}
Let $M$ be a compact $n$-manifold without boundary minimally immersed in some sphere $S^{n+l}$.
Let $\lambda_i$ and $\mu_i$ denote the eigenvalues of $\Delta$ on $M$ and $S^n$ respectively.
Define $\Theta(M)$ to be the ratio
$$\Theta(M) = \frac{\vol(M^n)}{\vol(S^n)}.$$
Then the trace of the heat kernel on $M$ and on $S^n$ satisfy
\begin{equation}
\label{heat:kernel:comparison}
\sum_{i=0}^{\infty}{e^{-\lambda_i t}} \le \Theta(M) \sum_{i=0}^{\infty}{e^{-\mu_i t}}
\end{equation}
for $t> 0$.
\end{theorem}
$\Theta(M)$ is the $(n+1)$-density of the singular point $\mathbf{0}$ of the minimal cone $C(M)$ over $M$.
By taking the limit as $t\ra \infty$ in (\ref{heat:kernel:comparison}) we see immediately that $\Theta(M) \ge 1$.
We will use Theorem \ref{cheng:li:yau:thm} in \S \ref{int:systems} to prove that the area of a minimal Legendrian
$2$-torus increases at least linearly with its spectral curve genus.

\subsection{Stable SLG cones and minimal immersions by 1st eigenfunctions}
\label{conformal:area:section}

Recall from Corollary \ref{lb:index} that a non-planar regular
\slg cone $C$ in $\C^n$ has Legendrian index at least $2n$.
$C$ is called \textit{Legendrian stable}
if  $\lind(C)= 2n$.
In terms of eigenvalues of $\Delta_{\Sigma}$ on the link $\Sigma$,
Legendrian stability is equivalent to
$$\lambda_1=n-1 \textrm{\ (with multiplicity precisely $2n$),\ } \quad \lambda_{2n+1}=2n.$$
When $n=3$ this puts considerable
restrictions on the minimal Legendrian surface $\Sigma$.
The main result of this section, Theorem \ref{ct:is:only:stable:t2}, shows there is essentially a unique such surface
if $\Sigma$ is a $2$-torus. To prove this theorem we need a result (Theorem \ref{elsoufi:ilias}) on
isometric minimal immersions induced by first eigenfunctions.


Let $\phi: M \ra S^{n}$ be a minimal immersion of a compact $m$-manifold $M$ into
any unit sphere. Since the restrictions of linear functions on $\R^{n+1}$ to $M$ are always
eigenfunctions for $\Delta_{M}$ with eigenvalue $m$, we always have the inequality $\lambda_1(M)\le m$.
Minimal immersions of $m$-manifolds into spheres with $\lambda_1=m$ are very special. They appear naturally
in other geometric problems and
naturally give rise to isometric minimal immersions induced by first eigenfunctions
of the Laplacian \cite{li:yau,montiel:ros}. In their work extending results of Li-Yau \cite{li:yau}
and Montiel-Ros \cite{montiel:ros} to higher dimensions, El Soufi-Ilias proved the following useful result.

\begin{theorem}[El Soufi-Ilias \cite{elsoufi:ilias1}]
\label{elsoufi:ilias}
Let $(M,g_0)$ be an $n$-dimensional compact Riemannian homogeneous manifold which is not
homothetic to $S^n$, and let $g=fg_0$ be any metric conformal to $g_0$. Suppose
that $(M,g)$ admits a minimal isometric immersion by first eigenfunctions into some unit sphere,
then $f$ is constant.
\end{theorem}
A simple corollary of the previous result is the following.
\begin{corollary}
\label{lam1min:is:flat}
Let $(T^2,g)$ be any metric on a $2$-torus. Suppose that $(T^2,g)$ admits a minimal
isometric immersion by first eigenfunctions into some unit sphere. Then $g$ must
be flat.
\end{corollary}
\begin{proof}
Any metric $g$ on a $2$-torus is conformal to some flat metric, and
any flat metric on $T^2$ is isometric to the quotient of $\R^2$ with its Euclidean metric
by the action of some lattice $\Gamma$. Clearly, $\R^2/{\Gamma}$ is homogeneous.
The result now follows immediately from Theorem \ref{elsoufi:ilias}.
\end{proof}
Using the previous corollary we can now give a very simple proof of the following result.
\begin{mtheorem}
\label{ct:is:only:stable:t2}
The Clifford torus $T^2$ is the unique (up to unitary equivalence) minimal Legendrian
$2$-torus in $S^5$ which is Legendrian stable, \textit{i.e.} has $\lind(T^2)=6$.
\end{mtheorem}
\begin{proof}
Let $\Sigma^2$ be any Legendrian stable minimal Legendrian surface in $S^5$,
and let $g$ denote the metric induced on $\Sigma$ from $S^5$. Then
$\lind(\Sigma)=6 \Rightarrow \lambda_1(\Sigma)=2$.
It follows that the first eigenfunctions of $(\Sigma,g)$ (the restrictions to $\Sigma$ of the linear functions on $\C^3$)
give a minimal isometric immersion into $S^5$.
%
%
If $\Sigma$ is a $2$-torus, Corollary \ref{lam1min:is:flat} implies that it must be
flat. But (up to finite covers) the Clifford torus is the unique  flat minimal
Legendrian torus in $S^5$ (up to unitary equivalence) \cite[\S 2.3]{haskins:slgcones}.
In fact, one can prove that no nontrivial finite cover of the Clifford torus
will still have the property that $\lambda_1=2$,
\textit{e.g.} see the final paragraph in the proof of Thm \ref{ct:only:hstable:t2}.
\end{proof}

There is a natural analogue of Theorem \ref{ct:is:only:stable:t2} for $\CP^2$ due to Urbano.
Although Theorems \ref{ct:is:only:stable:t2} and \ref{ct:only:hstable:t2} are closely related
there does not seem to be a way to deduce one from the other, especially since in higher dimensions the Clifford
torus is Hamiltonian stable but no longer Legendrian stable. Below we give a new proof of Theorem B using Theorem \ref{besson:thm}.
Our proof is considerably shorter than Urbano's original proof.
\begin{mtheorem}[\mbox{\cite[Cor. 2]{urbano:index}}]
\label{ct:only:hstable:t2}
The Clifford torus $T^2$ is the unique (up to unitary equivalence) Hamiltonian
stable minimal Lagrangian torus in $\CP^2$.
\end{mtheorem}
\begin{proof}
Let $\Sigma$ be a minimal Lagrangian torus in $\CP^2$. Then
$\Sigma$ is Hamiltonian stable if and only $\lambda_1=6$.
Any element
of $\text{SU(3)}$ acts in a Hamiltonian fashion on $\Sigma$ giving a variation through minimal
Lagrangian tori. Differentiating such a variation gives a Jacobi
field $\phi$ on $\Sigma$, \textit{i.e.} a solution of $\mathcal{J}\phi=0$,
and hence an eigenfunction satisfying $\Delta \phi = 6\phi$.
 Let $\text{G}$ be the identity component of the
subgroup of $\text{SU(3)}$ preserving $\Sigma$.
Then
$$ \dim{\ker(\Delta_{\Sigma} - 6)} \ge \dim{\text{SU(3)}} - \dim{\text{G}}.$$
For a $2$-torus, it follows from the integrable systems machinery described in \S \ref{int:systems}
that the possible choices for $\text{G}$ are $\{e\}, S^1, T^2$.
Hence we have
$\dim\ker(\Delta - 6) \ge 6$, with equality if and only if the symmetry
group of $\Sigma$ is $T^2$.

On the other hand, for any metric $g$ on a $2$-torus, Theorem \ref{besson:thm}
shows that the first nonzero eigenvalue of
$\Delta_g$ on functions has multiplicity less than or equal to 6.
Hence if $\Sigma$ is Hamiltonian stable then $\Sigma$ is invariant
under a $T^2$-subgroup of $\text{SU(3)}$ which makes $\Sigma$ into a homogeneous space. But homogeneity
on a $2$-torus implies that the metric is flat. But the only flat minimal Lagrangian torus in $\CP^2$ is up to a
finite Riemannian cover the Clifford torus \cite[\S 2.3]{haskins:slgcones}. The Clifford torus is already known to be Hamiltonian stable
by computing its spectrum via the dual period lattice as in Example \ref{clifford:torus}.

Clearly not all finite covers of the Clifford torus can be Hamiltonian stable \textit{e.g.}
the $\Z_3$ cover of the Clifford torus given by the projection of the Clifford torus in $S^5$
to $\CP^2$ is not Hamiltonian stable. In fact, no nontrivial finite cover of the Clifford torus is Hamiltonian stable.
First a short computation shows that the area of the Clifford torus in $S^5$ is
${4\pi^2}/\sqrt{3}$. Hence the area of the Clifford torus in $\CP^2$ is ${4\pi^2}/{3\sqrt{3}}$.
Hence if $\Sigma_k$ is any $k$-fold Riemannian cover of the Clifford torus in $\CP^2$  then
$$\Area(\Sigma_k) = \frac{4k\pi^2}{3\sqrt{3}}.$$
But if $\Sigma_k$ is Hamiltonian stable then $\lambda_1(\Sigma_k)=6$. Hence by (\ref{yang:yau}) we have
$$  \frac{24k\pi^2}{3\sqrt{3}} \le 16\pi \quad \implies k \le \frac{2\sqrt{3}}{\pi} \approx 1.103 < 2.$$
\end{proof}
\noindent
Several parts of the above proof fail in higher dimensions. Most crucially we have no higher dimensional analogue of Thm \ref{besson:thm}.

\subsection{Integrable systems and spectral curves of \slg $T^2$-cones}
\label{int:systems}
In the general setting of harmonic maps from a $2$-torus into a Lie group or symmetric
space, powerful integrable systems and loop group methods have been developed by various authors
\textit{e.g.} \cite{burstall:pedit,hitchin,pinkall:sterling,mcintosh,mcintosh:slg}.
In many cases these methods allow one to prove
that harmonic $2$-tori up to natural geometric equivalence are in
one-to-one correspondence with certain algebro-geometric data (its so-called \textit{spectral data}).

While in theory this yields all harmonic tori of certain types, in practice it is not obviously so useful.
There are two principal difficulties. The first is to recognize when some given spectral data actually corresponds
to a harmonic torus. This is problematic because of certain period conditions which must hold if the construction
is to yield a torus. A second more fundamental issue is that
it has proven very difficult to derive useful differential-geometric information directly from the algebro-geometric data.

For example, one suspects that harmonic tori with spectral curves of high genus must be complicated geometrically.
So one would like to prove that there exists a lower bound for the energy (area if it is also a conformal map)
of any harmonic torus with spectral curve of a given genus $g$
and that this bound increases with $g$.
A breakthrough in this problem was achieved only very recently.
For minimal tori in $S^3$ and constant mean curvature tori in $\R^3$ (the Gauss map of which is a non-conformal
harmonic map to $S^2$) Thm 6.7 of \cite{pedit:pinkall} proves a lower bound for the energy that grows quadratically with
the spectral curve genus. The method of proof by Ferus \textit{et. al.} uses the theory of ``quaternionic holomorphic geometry''
(and the quaternionic Pl\"ucker formula in particular)
in an essential way. It does not appear possible to fit minimal Legendrian tori into the same framework.

Nevertheless, we can prove analogous results in the minimal Legendrian setting. We exhibit lower bounds not only for
the area, but also the Legendrian index and the stability index of any minimal Legendrian $2$-torus in terms
of its spectral curve genus. We use completely different methods from \cite{pedit:pinkall} and use
in an essential way the results of \S \ref{spectral:geometry} on the geometry of eigenvalues of $\Delta$ on a surface.

We begin with results due to McIntosh on the existence of a ``spectral curve'' for any \slg $T^2$-cone.
The fundamental facts that the spectral curve machinery gives us (and the only ones which
we use heavily) are the following:

\smallskip
\emph{
Any minimal Lagrangrian torus in $\CP^2$ has an associated algebraic curve of even genus $g$, its \emph{spectral curve}.
Any minimal Lagrangian torus in $\CP^2$ with spectral curve of genus $g=2d$ comes in at least a $(d-2)$-dimensional family of geometrically
distinct tori.}

\smallskip

By combining these fundamental facts with results from \S \ref{spectral:geometry} we
prove Theorem \ref{specgenus:bds:index}. It gives an explicit lower bound for the Legendrian and stability indices
of any \slg $T^2$-cone in terms of its spectral curve genus.
Theorem \ref{area:spectral:genus} uses
Theorem \ref{cheng:li:yau:thm} to prove that the area
of any minimal Legendrian $2$-torus is bounded below by an explicit function of the spectral curve genus.
Theorems \ref{specgenus:bds:index} and \ref{area:spectral:genus} provide further evidence for the general principle that harmonic tori
with high spectral curve genus are complicated differential-geometrically.

\subsubsection{Spectral curves of \slg $T^2$-cones}
Every \slg $T^2$-cone $C$ in $\C^3$ determines
both a minimal Legendrian torus $\Sigma$ in $S^5$ (its link) and a minimal Lagrangian torus $\pi(\Sigma)$
in $\CP^2$, where $\pi:S^5\ra \CP^2$ is the Hopf map.
Conversely, \cite[Prop 4]{mcintosh} proves  that up to a possible $\Z_3$-cover we can reverse this construction
and lift any minimal Lagrangian torus in $\CP^2$ to a minimal Legendrian torus in $S^5$, which
forms the link of a \slg torus cone.

Integrable systems methods can be used to construct all minimal Lagrangian
tori in $\CP^2$.
Below we review McIntosh's work on the existence of spectral data for every minimal Lagrangian torus
in $\CP^2$. For further details of this construction we refer the reader to \cite{mcintosh,mcintosh:slg,joyce:intsys}.
However, it turns out that to prove our results we will need very little detailed information about
the spectral curve construction or other integrable systems intricacies.

McIntosh \cite{mcintosh:nonisotropic,mcintosh:remarks} proved that any harmonic $2$-torus
in $\CP^n$ which is not (complex) isotropic has associated to it unique spectral data and that this spectral data determines the torus up to congruence.
In particular, this proves that any minimal Lagrangian $2$-torus in $\CP^2$ is constructed from some spectral data.
What remains is to characterize precisely which spectral curves correspond to minimal tori which are Lagrangian.

The construction of \cite{mcintosh:nonisotropic} associates to each minimal torus in $\CP^2$ its spectral data $(X,\lambda,\mathcal{L})$.
This consists of a real algebraic curve $X$, a degree $3$ function $\lambda$ on $X$ and a holomorphic line bundle $\mathcal{L}$ over $X$.
$\lambda$ and $\mathcal{L}$ both respect the real involution of $X$. If the torus is also Lagrangian then $X$ must possess an
additional holomorphic involution $\mu$ which both $\lambda$ and $\mathcal{L}$ also respect. This implies, for example, that
$\mathcal{L}$ lies in (a translate of ) the Prym variety defined by $(X,\mu)$. The main result of \cite{mcintosh:slg} is that
these additional conditions characterize the spectral data of minimal Lagrangian $2$-tori.

\begin{theorem}[\mbox{\cite[Thm 1]{mcintosh:slg}}]
There is a bijective correspondence between:\\
(i) congruence classes of minimal Lagrangian tori $\phi: T^2 \ra \CP^2$
(with base point on $T^2$), and\\
(ii) equivalence classes of spectral data $(X,\lambda,\mathcal{L},\mu)$.\\
This spectral data consists of  a compact Riemann surface $X$ of even genus $g=2d$, $\lambda:X \ra \mathbb{\hat{C}}$
a degree three cover, $\mathcal{L}$ a holomorphic line bundle over $X$ of degree $g+2$
and $\mu$ a holomorphic involution on $X$ all of which satisfy some compatibility conditions
and some double periodicity conditions.
\end{theorem}
For the precise statement of the conditions that the spectral data must satisfy we refer the reader
to \cite[\S 4.1]{mcintosh:slg}. We will not need them in the remainder of this paper.
Taking into account all the conditions that the spectral data of a minimal Lagrangian torus must satisfy
and performing a parameter count suggests that for each $d$ there should be a countable family
of data $(X,\lambda,\mu)$ producing minimal Lagrangian tori. However, the important fact for our purposes
is that any minimal Lagrangian $2$-torus whose spectral curve has genus $g=2d$
comes in a real $(d-2)$-dimensional family of minimal Lagrangian tori.
This $(d-2)$-dimensional family arises by varying the line bundle $\mathcal{L}$
in the $d$-dimensional Prym variety of $(X,\mu)$, and discarding the choice of base point on $T^2$.

Given the periodicity conditions that the spectral data must satisfy it is not
at all obvious that spectral data corresponding to minimal Lagrangian tori exist
for each even genus $g=2d$. If $X$ has genus $0$, then the associated minimal Lagrangian torus
is the Clifford torus. If $X$ has genus $2$, it follows from \cite{mcintosh}
that the associated minimal Lagrangian torus is $S^1$-equivariant. Such tori have been studied by
several authors including the author \cite{haskins:slgcones,haskins:thesis}, \cite{castro:94} and \cite{joyce:symmetries}.
The explicit construction of such tori in terms of elliptic functions and integrals will be reviewed in \S \ref{s1:invariant:tori}.
Joyce has given some special examples of minimal Lagrangian tori with spectral curve genus equal to 4 \cite{joyce:intsys}.
For $g>4$, no explicit examples are known of spectral curves which satisfy the periodicity conditions
necessary to correspond to a minimal Lagrangian torus.
However, Carberry and McIntosh recently proved the following result.
\begin{theorem}[\cite{carberry:mcintosh}]
\label{carb:mc:thm}
For each $d>2$, there are countably infinite families of spectral data $(X,\lambda,\mu)$
of spectral genus equal to $2d$ which give minimal Lagrangian tori in $\CP^2$.
\end{theorem}

\smallskip
\noindent
The facts from the discussion above that we will need in the rest of this paper are summarized below.
\begin{enumerate}
\item Every minimal Lagrangian torus $T^2$ in $\CP^2$ has associated spectral data
which determines the torus up to congruence.
\item
The spectral data includes an algebraic curve of even genus $g=2d$ called the \textit{spectral curve}.
$g$ is called the \textit{spectral curve genus} of the torus $T^2$.
\item
Any minimal Lagrangian torus of spectral curve genus $g=2d$
comes in a family of non-congruent minimal Lagrangian tori of real dimension at least $(d-2)$.
\end{enumerate}

\subsubsection{Area, Legendrian index, stability index and spectral curve genus}
\label{area:index:genus}
In this section we prove that the Legendrian index, the stability index and the
area of a \slg $T^2$-cone all increase at least linearly with the spectral curve genus.

\begin{mtheorem}
\label{specgenus:bds:index}
Let $C$ be a \slg $T^2$-cone in $\C^3$ with an isolated singularity at $0$,
and let $g=2d$ denote the spectral curve genus of the associated minimal Lagrangian
torus $\phi: T^2 \ra \CP^2$. Then for any $d\ge 3$ we have
\beq
\label{lind:specgenus}
\lind(C) \ge \max\left\{\left[\tfrac{1}{2}d\right]^+,7\right\}
\end{equation}
and
\begin{equation}
\label{sind:specgenus}
\sind(C) \ge \max\left\{\left[\tfrac{3}{2}d-8\right]^+,d-1 \right\},
\end{equation}
where  $[x]^+$ denotes the smallest integer not less than $x\in \R$, \textit{i.e.} $[x]^+:=\min{\{n\in \Z \, |\, n \ge x \}}$.
\end{mtheorem}

\begin{proof}
Given $\phi: T^2 \ra \CP^2$ a minimal Lagrangian torus, deformations of
$\phi$ through minimal Lagrangian surfaces give rise to Jacobi fields on $T^2$, \textit{i.e.} solutions of $\mathcal{J}v=0$.
Clearly, any element of $\text{SU(3)}$ determines a deformation of $\phi$
as a minimal Lagrangian surface. For $d>1$, the corresponding minimal Lagrangian torus has no
continuous symmetries and hence $\text{SU(3)}$ gives rise to an $8$-dimensional space of Jacobi fields.
For $d>2$, there is at least a $(d-2)$-dimensional family of additional Jacobi fields
arising from varying the line bundle $\mathcal{L}$ in the spectral curve description of $\phi$
within $\textrm{Prym}(X,\mu)$.

But Jacobi fields on $T^2$ arising from Lagrangian deformations
are in one-to-one correspondence with eigenfunctions $v$ satisfying $\Delta v=6v$,
where $\Delta$ is acting on functions on $T^2$.
So the multiplicity of $6$ as an eigenvalue of $\Delta$ is at least $\dim{\text{SU(3)}}+(d-2)=6+d$.
Trivially, we have the lower bound $\sind(C)\ge d-2$. We get an improvement of $1$ to this lower bound
using the fact that from Theorem \ref{ct:is:only:stable:t2}, $\lind(C)\ge7$ whenever $d>0$.

We can also prove a lower bound for $\lind(C)$ and use this to improve the lower bound for $\sind$.
Let $i\in \N$ be the unique index so that $\lambda_{i-1}<6$ and $\lambda_i=6$.
Theorem \ref{besson:thm} implies that
this index $i$ must satisfy $ 4 + 2i \ge d + 6$.
Hence there are at least $[\frac{d}{2}]^+$ eigenvalues of $\Delta$ in the range $(0,6)$.
Since by Theorem \ref{ct:is:only:stable:t2}, $\lind(C)\ge 7$ whenever $d>0$, we have
obtained the claimed lower bound for $\lind(C)$.
Finally, $\sind(C) \ge (d-2) + \frac{d}{2} - 6  = \frac{3d}{2}-8$ as claimed.
\end{proof}

\begin{remark}
\label{quadratic:index:rmk}
Theorem \ref{specgenus:bds:index} shows that $\lind$ and $\sind$ grow at least linearly with
the spectral curve genus $g=2d$, using the result (Thm \ref{besson:thm}) that the multiplicity of an eigenvalue of $\Delta$ on a surface
grows at most linearly with its index. As we already noted in Remark \ref{mult:growth:rmk}, in all known examples
the multiplicity $m_i$ of an eigenvalue $\lambda_i$ of $\Delta$ on a surface grows at most like the square root of the index.
If one could prove $m_i \le c \sqrt{i} + d$ for some constants $c>0$ and $d\in \R$, then the proof given in
Theorem \ref{specgenus:bds:index} would show that $\lind$ and $\sind$ must grow at least quadratically with the
spectral curve genus.
\end{remark}

\begin{mtheorem}
\label{area:spectral:genus}
Let $C$ be a \slg $T^2$-cone in $\C^3$ with an isolated singularity at $0$,
and let $g=2d$ denote the spectral curve genus of the associated minimal Lagrangian
torus $\phi: T^2 \ra \CP^2$. Let $\Sigma$ be the link of the cone $C$. Then for any $d\ge 3$ we have
\begin{equation}
\label{genus:bounds:area}
\Area(\Sigma) > \tfrac{1}{3} d\pi.
\end{equation}
\end{mtheorem}
\begin{proof}
We will prove this result using Theorem \ref{cheng:li:yau:thm} and Theorem \ref{specgenus:bds:index}.
To make use of Theorem \ref{cheng:li:yau:thm} we need to have an upper bound for the trace of the heat kernel $H$
on the standard $S^2$. This is not difficult to obtain. From the explicit expression for the eigenvalues $\mu_i$ of $\Delta$ on $S^2$, the trace of the heat kernel $H$ is given by
$$ \Trace_{S^2}{H}(t)= \int_{S^2}H(x,x,t)\, dx = \sum_{i=0}^{\infty}{e^{-\mu_i t}} = \sum_{i=0}^{\infty}{(2i+1)e^{-i(i+1)t}}$$
for any $t>0$.
Define the function $f(t,x)=(2x+1) e^{-x(x+1)t}$ and let $a_i(t)=f(t,i)$ for $i\in \N$. Then
$\Trace_{S^2}{H}(t) = \sum_{i=0}^{\infty}{a_i(t)}$. A calculation shows that $f$ is a decreasing function
of $x$ provided $t> \frac{2}{(2x+1)^2}$. Hence by the Integral Test, we have
$$ \Trace_{S^2}H(t)
\le \sum_{i=0}^{k}{(2i+1)e^{-i(i+1)t}} + \int_{k}^{\infty}{(2x+1)e^{-x(x+1)t}dx} = \sum_{i=0}^{k}{(2i+1)e^{-i(i+1)t}} + \tfrac{1}{t}e^{-k(k+1)t}$$
provided $t> \frac{2}{(2k+1)^2}$.

Let $\lambda_i$ be the eigenvalues of $\Delta$ on $\Sigma$, and define $N_{\Sigma}(\lambda)$ as
in \S \ref{harm:fn:cone}. Then for the trace of the heat kernel on $\Sigma$ we have
$$\Trace_{\Sigma}H(t) = \sum_{i=0}^{\infty}{e^{-\lambda_i t}} > 1 + (N_{\Sigma}(2)-1)e^{-6t}.$$
Hence given any $k\in \N$, applying Theorem \ref{cheng:li:yau:thm} implies that
\begin{equation}
\label{area:neval:t}
\Area(\Sigma)> \frac{4\pi \left(1 + (N_{\Sigma}(2)-1)e^{-6t}\right)}{\left(\sum_{i=0}^{k}{(2i+1)e^{-i(i+1)t}} + \tfrac{1}{t}e^{-k(k+1)t}\right)}
\end{equation}
holds provided $t>\frac{2}{(2k+1)^2}$. By taking $k=4$, and letting $t=\frac{4}{25}>\frac{2}{(2k+1)^2}=\frac{2}{81}$ from (\ref{area:neval:t})
we obtain
\begin{equation}
\label{area:gt:nsigma}
 \Area(\Sigma) > 4\pi \left( \frac{1}{7} + \frac{1}{18}(N_{\Sigma}(2)-1) \right).
\end{equation}
From the proof of Theorem \ref{specgenus:bds:index} we see that $N_{\Sigma}(2) -1 \ge \frac{3}{2}d -2$, where $g=2d$ denotes
the spectral curve genus of $\Sigma$. Hence we have
$$ \Area(\Sigma) > \frac{4\pi}{18} \left(2 + \tfrac{3}{2}d -2\right)  = \frac{1}{3}d\pi$$
as claimed.
\end{proof}

\begin{remark}
\label{quadratic:area:growth}
As in Remark \ref{quadratic:index:rmk},
if one could prove improve the bounds from Theorem \ref{besson:thm} for the multiplicity $m_i$ of the eigenvalue $\lambda_i$ to
$m_i \le c \sqrt{i} + d$ for some constants $c>0$ and $d\in \R$, then the previous proof
would show that $\Area(\Sigma)$ must grow at least quadratically with the
spectral curve genus of $\Sigma$. We conjecture that the area does indeed grow quadratically with spectral curve genus.
For minimal tori in $S^3$
it was proved in \cite{pedit:pinkall} using quaternionic holomorphic methods that the area does grow quadratically with
the spectral curve genus.
\end{remark}
\begin{remark}
The constants appearing in the conclusion of Theorem \ref{area:spectral:genus} are not the optimal constants that one
could achieve using this method. Rather they were chosen to keep the algebraic manipulations simple.
\end{remark}

\begin{remark}
There is a natural higher dimensional analogue of Theorem \ref{area:spectral:genus}. Since $\text{Spec}(\Delta)$
is explicitly known for $S^{n-1}$ for all $n$, it is once again possible to obtain upper bounds for the trace of the heat
kernel on $S^{n-1}$ for $n-1\ge 3$ \cite[Lemma 3]{cheng:li:yau}. Let $\Sigma^{n-1}$ be any compact minimal Legendrian submanifold
of $S^{2n-1}$, and for $\lambda>0$ define $N_\Sigma(\lambda)$ as in (\ref{nsig:plus}).
Combining the upper bound for the trace of the heat kernel on $S^{n-1}$ with Theorem \ref{cheng:li:yau:thm} leads
to the existence of some positive constant $C_{n,\lambda}$ depending on $n$ and $\lambda$ but not on $\Sigma$
so that
$$ \vol{(\Sigma^{n-1})} > C_{n,\lambda} N_{\Sigma}(\lambda) \vol(S^{n-1}).$$
In particular, by taking $\lambda=2$ we see that the volume of any minimal Legendrian submanifold $\Sigma^{n-1}$ must grow at least linearly with the
number of eigenvalues in the range $[0,2n]$. However, in higher dimensions we have no equivalent of
the spectral curve description through which we can force the Legendrian index or $N_{\Sigma}(2)$ to grow.
\end{remark}

Theorems \ref{specgenus:bds:index} and \ref{area:spectral:genus} together with Theorem \ref{carb:mc:thm} imply
\begin{corollary}
There exist \slg $T^2$-cones with arbitrarily large Legendrian index,  stability index and area.
\end{corollary}

In \S \ref{s1:invariant:tori} we will prove that the same result holds even if one restricts to $S^1$-invariant $T^2$-cones.
This shows that there is no \textit{upper} bound for the area, Legendrian index or stability index just in
terms of the spectral curve genus. The  proof is independent of the results in this section.

\section{Explicit \slg $T^2$-cones and refined lower bounds for {l-index}}
\label{s1:invariant:tori}

Explicit descriptions of minimal Legendrian tori in $S^5$ invariant
under a $1$-parameter subgroup of $\text{SU(3)}$ were given by the author in \cite{haskins:slgcones,haskins:thesis}.
In \S \ref{s1:descriptions} we recall the main features of this construction. In \S \ref{s1inv:indexbounds}
we use the explicit descriptions to obtain refined estimates
of the Legendrian and stability indices of $S^1$-invariant minimal Legendrian $2$-tori.
We use these explicit descriptions again in \S \ref{areas} to obtain analytic and numerical results
on the area of $S^1$-invariant minimal Legendrian tori.
The main results of the section are Theorems \ref{lindex:tmn} and \ref{lindex:u0j}.

\subsection{Explicit descriptions of $S^1$-invariant minimal Legendrian tori}
\label{s1:descriptions}

For each pair $(\alpha,J) \in [0,1] \times [0,\frac{1}{3\sqrt{3}}]$ we will associate an minimal
Legendrian immersion $u_{\alpha,J}: \R^2 \ra S^5$ as follows.
Given $(\alpha,J) \in [0,1] \times [0,\frac{1}{3\sqrt{3}}]$, define
$$\vec{\lambda} = (1,\alpha, -1-\alpha ), \quad \vec{\mu} = \vec{1} \times \vec{\lambda}, \quad \vec{J}=J \vec{\mu},
\quad A = \sqrt{-1} \diag{(\vec{\lambda})} \in \mathfrak{su}(3).$$

Define $\gamma_{\alpha,J}(t)$ to be the unique solution to the first order nonlinear ODE
\beq
\label{gammadot}
\frac{1}{4} \dot{\gamma}^2 + J^2 = \gamma^3 \mu_1 \mu_2 \mu_3 + \frac{\gamma^2}{3} \sum_{i\neq j}{\mu_i \mu_j} + \frac{1}{27}
\end{equation}
with the property that $ \gamma(0) = \gamma_{min}$
where $\gamma_{min}$ denotes the minimum value attained by $\gamma$.
$\gamma(t)$ can be expressed in terms of the Jacobi elliptic function $\sn$.
The definition of $\sn$ and the other Jacobi elliptic functions, together with some basic properties,
are given in Appendix \ref{appendix:elliptic:fn}.
For $(\alpha,J) \in  [0,1) \times [0,\frac{1}{3\sqrt{3}})$, there are three solutions $\Gamma_1, \Gamma_2, \Gamma_3$
to (\ref{gammadot}) with $\dot{\gamma}=0$. Let us label these solutions so that
$\Gamma_2\le 0 \le \Gamma_1 \le \Gamma_3$. We can rewrite (\ref{gammadot}) as
\begin{equation}
\label{gammadot2}
\dot{\gamma}^2 = 4 \mu_1 \mu_2 \mu_3 (\gamma - \Gamma_1 ) (\gamma - \Gamma_2 ) (\gamma - \Gamma_3 ).
\end{equation}
The following proposition  expresses $\gamma(t)$ in terms
of the Jacobi elliptic function $\sn$.

\begin{prop}[\mbox{Haskins, \cite[Prop 4.2]{haskins:slgcones}}]
\label{gamma}
The unique solution of (\ref{gammadot2}) with $\gamma(0) = \gamma_{min} = \Gamma_2$
is given by
\beq
\gamma(t) = \Gamma_2 - (\Gamma_2 - \Gamma_1) \sn^2{(rt,k)}
\end{equation}
where
\beq
\label{rk:gamma}
 r^2 = \mu_1 \mu_2 \mu_3 (\Gamma_3 - \Gamma_2), \quad k^2 = \frac{\Gamma_2 - \Gamma_1}{\Gamma_2 - \Gamma_3}
\end{equation}
and $\sn$ is the Jacobi elliptic sn-noidal function.
\end{prop}
The proof is a straightforward computation using elementary properties of the Jacobi elliptic functions.
Since the period of $\sn{(t,k)}$ is $4K(k)$ where $K(k)$ is the complete elliptic integral of the first kind defined
in (\ref{ellipK}),
Proposition \ref{gamma} implies that $\gamma$ is periodic of basic period $T_{\alpha,J}= \frac{2K(k)}{r}$.

Define the following nine functions $R_i, {\theta}_i, z_i$ of $t$ by
\beq
\label{Ri}
R_i^2(t) = \gamma(t) \mu_i + \frac{1}{3},
\end{equation}
\beq
\label{thetai}
\quad \dot{\theta_i}(t) = \frac{J\mu_i}{R_i^2}, \quad \textrm{with\ } \theta_i(0)=0,
\end{equation}
\beq
\label{zi}
z_i(t) = R_i(t) \exp{(\sqrt{-1}\theta_i(t))}
\end{equation}
for $i=1,2,3$. Define $u_{\alpha,J}$ by
\beq
\label{u:alpha:j}
u_{\alpha,J}(s,t) = e^{As}\left(z_1, z_2, z_3\right).
\end{equation}

\begin{theorem}[\mbox{Haskins, \cite[Thm D]{haskins:slgcones}}]
For any $(\alpha,J) \in [0,1] \times [0,\frac{1}{3\sqrt{3}}]$, $u_{\alpha,J}$ defined by (\ref{u:alpha:j})
is a minimal Legendrian immersion $\R^2 \ra S^5$ with the following properties:\\
(i) $u_{\alpha,J}$ is invariant under the $1$-parameter subgroup of \text{SU(3)} generated by $A \in \mathfrak{su}(3)$.\\
(ii) For $\alpha=1$ or $J= \frac{1}{3\sqrt{3}}$ these immersions all describe the Clifford torus,
but otherwise all the immersions are geometrically
distinct.
\end{theorem}

For $\alpha \in \Q \cap [0,1]$, the $1$-parameter subgroup of $\text{SU(3)}$ generated by $A$ closes up
to form a circle -- the corresponding immersions $u_{\alpha,J}$ then factor through a minimal Legendrian cylinder
which in general will have dense image in $S^5$. For special values of $J$, the immersion $u_{\alpha,J}$ will
admit a second independent period and hence $u_{\alpha,J}$ will factor through a minimal Legendrian torus.
This happens whenever the curve in $\C^3$ defined by
$$ t \ra \left(z_1(t), z_2(t), z_3(t)\right)$$
forms a closed curve. Since $\gamma$ and hence each of the $R_i$ are periodic of period $T_{\alpha,J}$, this
happens if and only if $\theta_i(T_{\alpha,J}) \in \pi \Q$ for $i=1,2,3$. However, these three rationality conditions
are not independent. From  \cite[Lem 5.1]{haskins:slgcones} the sum of the angles $\sum{\theta_i(t)}$
satisfies
\begin{equation}
\label{gammadot:anglesum}
\dot{\gamma}(t) = 2J \tan{\sum{\theta_i(t)}}.
\end{equation}
It follows from this \cite[Cor 5.2]{haskins:slgcones} that
\begin{equation}
\label{theta:sum:period:zero}
\sum{\theta_i(T)} = 0.
\end{equation}

The following proposition gives the conditions under which $u_{\alpha,J}$ admits two independent periods.

\begin{prop}[\mbox{Haskins, \cite[Prop 5.3]{haskins:slgcones}}]
\label{double:period:prop}
(a) $(\sigma, \tau)$ is a period of $u_{\alpha,J}$ $\Rightarrow$ $\tau \in \N T_{\alpha,J}$.\\
(b) $u$ admits two independent periods $\Rightarrow$ it admits a period of the form $(\sigma,0)$.\\
(c) $u$ admits a period of the form $(\sigma,0)$ if and only if $\alpha \in \Q$.\\
(d) $u$ admits two independent periods if and only if
$ \alpha$ and $\frac{1}{2\pi} (\alpha \theta_1(T) - \theta_2(T) ) \in \Q.$
\end{prop}

There are two cases of the previous proposition where it is straightforward to verify that the
rationality condition needed for double periodicity holds: when $J=0$ or when $\alpha=0$.

For the case $J=0$, it follows from (\ref{thetai}) that $\theta_i(t)\equiv0$ for $i=1,2,3$.
Hence the rationality conditions are satisfied by any $u_{\alpha,0}$ with $\alpha \in \Q$.
In this case the period lattices of the resulting minimal Legendrian tori are described by the
following proposition.
\begin{theorem}[Haskins, \mbox{\cite[Prop 5.4]{haskins:slgcones}}]
\label{u:alpha:0}
Let $\alpha=\frac{m}{n}$, where $0<m < n \in \N$ and $(m,n)=1$. Then the immersion
$u_{\alpha,0}$ is doubly periodic and gives rise to an embedded minimal Legendrian torus $T_{m,n}$ in $S^5$.
The period lattice of $u_{m/n,0}$ is:
rectangular with basis $\omega_1 = (2n\pi,0)$, $\omega_2 = (0,4K(k)/r)$ if $mn$ is even, and otherwise
has basis $\omega_1 = (2n\pi,0)$, and $\omega_2= (n\pi, 2K(k)/r)$.
\end{theorem}

For the case $\alpha =0$, the rationality conditions reduce to $\theta_2(T) \in \Q \pi$.
One can express $\theta_2(T)$ in terms of elliptic integrals of the third kind
and prove that $\theta_2(T_J)$ is a monotone continuous function of $J\in (0,\frac{1}{3\sqrt{3}})$ with
$\lim_{J\ra 0}{\theta_2(T_J)} = \pi$, and $\lim_{J \ra 1/3\sqrt{3}}{\theta_2(T_J)} = \frac{2\pi}{\sqrt{3}}$.
\begin{theorem}[Haskins, \mbox{\cite[Thm E]{haskins:slgcones}}]
\label{u:zero:j}
For a dense set of $J \in (0, \frac{1}{3\sqrt{3}})$, the immersion $u_{0,J}$ is
doubly periodic and hence gives rise to a minimal Legendrian torus.
Given any $J$ such that the immersion $u_{0,J}$ is doubly periodic then
$\theta_2(T) =  2\pi \tfrac{M}{N}$ for relatively prime positive integers $M$ and $N$. The period lattice of $u_{0,J}$ is:
rectangular with basis $\omega_1 = (2\pi,0)$, $\omega_2 = (0,NT)$ if $M$ is even, and otherwise
has basis $\omega_1= (2\pi,0)$, $\omega_2= (\pi,NT)$.
\end{theorem}

\begin{proof} The first part follows from the monotonicity of $\theta_2(T_J)$ and Proposition \ref{double:period:prop}.
The part concerning the structure of the period lattice follows
by an application of the method used below to prove results on the period lattice in Theorem \ref{general:double:period}.
%
%
\end{proof}

With more work one can prove that for each $\alpha \in \Q \cap (0,1)$
there are infinitely many values of $J$ for which $u_{\alpha,J}$ is doubly periodic
(see also \cite[Thm 8.5]{joyce:symmetries}).

\begin{theorem}
\label{general:double:period}
For each $\alpha = \frac{m}{n} \in \Q \cap (0,1)$ there is a countable dense set of $J\in (0,\frac{1}{3\sqrt{3}})$ for which
$u_{\alpha,J}$ is doubly periodic.
Whenever $u_{\alpha,J}$ is doubly periodic then $\theta_2(T)-\alpha \theta_1(T) = 2\pi\frac{M}{N}$ for relatively prime integers $M$ and $N$.
The period lattice of $u_{\alpha,J}$ is generated by $\omega_1 = 2\tilde{n}\pi$ and $\omega_2 = \tilde{N}(-\theta_1(T_J),T_J)$,
where $\tilde{n}$ and $\tilde{N}$ are defined by $\tilde{n}\hcf(n,N)=n$ and $\tilde{N}\hcf(n,N)=N$.
\end{theorem}

\begin{proof}
By Proposition \ref{double:period:prop}(d) we need to prove there is a countable dense set of
$J\in (0,\frac{1}{3\sqrt{3}})$ for which $\Theta(\alpha,J):=\theta_2(T_J) - \alpha \theta_1(T_J) \in 2\pi \Q$.
Since for fixed $\alpha$, $\Theta(\alpha,J)$ is a real analytic
function of $J\in (0,\frac{1}{3\sqrt{3}})$ this will follow if we can prove that $\Theta(\alpha,J)$ is a non-constant function of $J$.
In fact, we will prove
\begin{equation}
\label{theta:limits}
\lim_{J\ra 0}{\Theta(\alpha,J)} = (1+\alpha)\pi \ < \lim_{J\ra \frac{1}{3\sqrt{3}}}{\Theta(\alpha,J)} =
2\pi \left(\frac{1+\alpha+\alpha^2}{3}\right)^{\frac{1}{2}}
\end{equation}
holds for any $\alpha \in (0,1)$. In these limits we mean a limit from the right as $J\ra 0$ and a limit from the left as $J\ra \frac{1}{3\sqrt{3}}$.
Throughout the rest of the proof we will suppress this to reduce notation.

The proof of the second equality in (\ref{theta:limits}) is quite straightforward.
It follows directly from Proposition \ref{area:limits}(b).
Since from (\ref{theta:sum:period:zero}) we have $\sum{\theta_i}(T_J)=0$, the first equality in (\ref{theta:limits}) will follow if we show that
\begin{equation}
\label{theta23:limits}
\mathrm{(a)\ } \lim_{J\ra 0} \theta_3(T_J)=0  \quad \quad  \mathrm{\ and \quad \quad (b)\ } \lim_{J\ra 0}\theta_2(T_J)=\pi,
\end{equation}
hold for any $\alpha\in(0,1)$.

\medskip
\textit{Proof of (\ref{theta23:limits}a)}:
It follows from (\ref{gammadot}) and (\ref{rk:gamma}) that for any $\alpha \in (0,1)$
$$  k_{0,\alpha}^2:= \lim_{J\ra 0}{k^2} = \frac{1-\alpha^2}{1+2\alpha}  < 1 \textrm{\ and }  \lim_{J\ra0}{r^2} = 1+2\alpha.$$
The period $T_J$ is given by $\frac{2K(k)}{r}$, where $K(k)$ is the complete elliptic integral of the first kind
defined in (\ref{ellipK}).
$K$ is a positive strictly increasing continuous function of $k\in (0,1)$ with $\lim_{k\ra 0}{K(k)}=\frac{\pi}{2}$ and $\lim_{k\ra 1}{K(k)}=+\infty$.
Hence for any $\alpha\in (0,1)$ we have
$$ \lim_{J\ra0}{T_J} = \frac{2K(k_{0,\alpha})}{1+2\alpha} < \infty$$
(on the other hand for $\alpha=0$, we have $k_{0,0}=1$ and hence $\lim_{J\ra 0}{T_J}= +\infty$).
In particular, for any $\alpha\in (0,1)$ $T_J$ is bounded on $(0,\frac{1}{3{\sqrt{3}}})$.
Also, the minimum value of $R_3$ as a function of $t$ is an increasing function of $J$ with
$ \lim_{J\ra 0}({R_3})_\text{min} = \left(\frac{\mu_1-\mu_3}{3\mu_1}\right)^\frac{1}{2} = \left(\frac{3\alpha}{1+2\alpha}\right)^{\frac{1}{2}}.$
In particular,
$ R_3^2 > \frac{3\alpha}{1+2\alpha}$
holds for all $t$ and any $J>0$.
Combining these observations we find that for any $\alpha\in (0,1)$ there exists a constant $C(\alpha)$,
depending on $\alpha$ but independent of $J$, such that
$$|\theta_3(T_J)|  \le  \  \int_0^{T_J}|\dot{\theta}_3(t) dt |
 \le \ \int_0^{T_J}\frac{J|\mu_3|}{R_3^2} dt
\   \le  \int_0^{T_J} \frac{J|\alpha-1|}{(R_3^2)_\text{min}} dt  \le  \frac{J(1-\alpha)(1+2\alpha)}{3\alpha} T_J  \le \
 C(\alpha) J\\
$$
holds for all $J>0$. Hence $\lim_{J\ra 0}\theta_3(T_J)=0$ as claimed.

\textit{Proof of (\ref{theta23:limits}b)}:
Since $R_2$ is reflection symmetric about $t=\frac{1}{2}T_J$ it is enough to prove that
\begin{equation}
\lim_{J\ra 0} \theta_2\left(\tfrac{1}{2}T_J\right) = \frac{\pi}{2}.
\end{equation}
Let $\dot{\gamma}_+$ denote the value of $\dot{\gamma}$ at the time $t_+ \in (0, \frac{1}{2}T_J)$ when $\tan{\sum{\theta_i(t)}}$ attains its maximum.
$\dot{\gamma}_+$ and $t_+$ depend on both $\alpha$ and $J$, but to ease notation we will suppress the explicit dependence on $\alpha$ and $J$.
It follows from (\ref{gammadot:anglesum}) that the maximum value of $\tan{\sum{\theta_i(t)}}$ occurs when $\ddot{\gamma}=0$.
Differentiating (\ref{gammadot}) shows that $\gamma$ satisfies $\ddot{\gamma} = - 6|Az|^2 \gamma$. Hence
$\ddot{\gamma}=0 \Leftrightarrow \gamma=0$. Substituting for $\gamma$ in (\ref{gammadot}) it follows that
$\dot{\gamma}_+ = 2 \sqrt{\frac{1}{27}-J^2},$
and so
$\max\left({\tan{{\textstyle\sum{\theta_i}}}}\right) =  \tan{{\textstyle{\sum{\theta_i}}}}(t_+) =
\frac{\dot{\gamma}_+}{2J} = \frac{1}{J}\sqrt{\frac{1}{27}-J^2}.$
Hence we have
\begin{equation}
\label{theta:sum:limit}
\lim_{J\ra0} \tan{{\textstyle{\sum{\theta_i(t_+)}}}} = +\infty.
\end{equation}
 From (\ref{gammadot:anglesum}) it follows that
$ R_1R_2R_3 \cos{({\textstyle{\sum{\theta_i}}})} = J.$ For any $J>0$, all the $R_i$ are positive and hence it follows that $\cos{\sum{\theta_i}}>0$.
Since we defined $\theta_i(0)=0$, this implies that $\sum{\theta_i(t)} \in (-\frac{\pi}{2},\frac{\pi}{2})$ for all $J>0$.
Hence (\ref{theta:sum:limit}) implies that
\begin{equation}
\label{theta:sum:limit2}
\lim_{J\ra0}{\textstyle{\sum{\theta_i(t_+)}}} = \frac{\pi}{2}.
\end{equation}
Since $\gamma(t) \le 0$ for any $t\in [0,t_+]$, and $R_1^2 = -(1+2\alpha)\gamma + \frac{1}{3}$, then we have
$R_1^2 \ge \frac{1}{3}$ for $t\in [0,t_+]$. Hence,
$$ 0 < - \theta_1(t_+) = \int_0^{t_+} -\dot{\theta_1}(t)dt = \int_0^{t_+} \frac{J(1+2\alpha)}{R_1^2} dt \le 3J(1+2\alpha)t_+
\le \tfrac{3}{2}J(1+2\alpha) T_J \le C(\alpha)J$$
holds for some constant $C(\alpha)$ depending on $\alpha$ but independent of $J$.
So
\begin{equation}
\label{theta1:tplus}
\lim_{J\ra0}{\theta_1(t_+)} = 0,
\end{equation}
holds for any $\alpha\in (0,1)$. Also, since $0<t_+<\frac{1}{2}T_J$ the proof of (\ref{theta23:limits}a) implies that
\begin{equation}
\label{theta3:tplus}
\lim_{J\ra0}{\theta_3(t_+)}=0.
\end{equation}
Combining (\ref{theta:sum:limit2}), (\ref{theta1:tplus}) and (\ref{theta3:tplus}) implies
\begin{equation}
\label{theta2:tplus}
\lim_{J\ra0}{\theta_2(t_+)} = \frac{\pi}{2}.
\end{equation}
Since $\gamma(t) \ge 0$ for any $t\in [t_+,\frac{1}{2}T_J]$, and $R_2^2 = (2+\alpha)\gamma+ \frac{1}{3}$, then
$R_2^2 \ge \frac{1}{3}$ for $t\in [t_+, \frac{1}{2}T_J]$. Hence,
$$ 0 < \theta_2(\tfrac{1}{2}T_J) - \theta_2 (t_+) = \int_{t_0}^{\frac{1}{2}T_J} \dot{\theta}_2(t)dt \le 3J(2+\alpha)(\tfrac{1}{2}T_J - t_+) \le D(\alpha)J$$
holds for some constant $D(\alpha)$ depending on $\alpha$ but independent of $J$. Combining this with (\ref{theta2:tplus}) we have
$$ \lim_{J\ra0}{\theta_2(\tfrac{1}{2}T_J)} = \lim_{J\ra0}{ \theta_2(t_+)} = \frac{\pi}{2},$$
as required.\\

Let us finish by analyzing the possible periods of any doubly periodic $u_{\alpha,J}$. The condition that $(s,pT)$ where $p\in \N$
be a period of $u_{\alpha,J}$ is equivalent to the following three equations
$$
\exp i(s+p\theta_1(T))=1, \quad \exp i\left(\frac{m}{n}s+ p\theta_2(T)\right) =1, \quad \exp i\left(-\frac{n+m}{n}s + p\theta_3(T)\right) =1.
$$
Since $\sum{\theta_i}(T)=0$, the first two equations imply the third. By Proposition \ref{double:period:prop}(d) the existence of two periods
of $u_{\alpha,J}$ implies
$$ \Theta= \theta_2(T) - \alpha \theta_1(T) = \frac{2\pi M}{N}$$
for some integers $M$ and $N$. Given any two integers $p_1$ and $p_2$ we get a solution of the periodicity equations by setting
\begin{equation}
\label{s:p:eqn}
\text{(a)\ }s = \frac{N}{M} (\theta_2(T)p_1 - \theta_1(T)p_2) \text{\quad and \quad } \text{(b)\ } p = \frac{N}{M}\left( p_2 - \frac{m}{n}p_1 \right).
\end{equation}
However, since we require that $p\in \N$, most pairs of integers $p_1$ and $p_2$ will not give us valid solutions of the periodicity equations.
Taking $(p_1,p_2)=(0,M)$ gives $(s,p) = (-N \theta_1(T), N)$, while $(p_1,p_2)=(n,m)$ gives $(s,p)= (2n\pi,0)$.
Hence the period lattice of $u_{\alpha,J}$ always contains $\omega_1=(2\pi n,0)$ and $\omega_2= (NT,-N\theta_1(T))$.

We now determine the most general element of the period lattice of $u_{\alpha,J}$. Since $M$ and $N$ are relatively prime and
$p_2 - \frac{m}{n}p_1 \in \Q$, it follows from (\ref{s:p:eqn}b) that $p=k N_f$ , where $k\in \Z$ and $N_f$ is some factor of $N$.
Define $\tilde{N}_f$ to be the integer $\frac{N}{N_f}$. Substituting for $p$ into (\ref{s:p:eqn}) gives us
\begin{equation}
\label{p:eqn}
np_2-mp_1 = \frac{Mkn}{\tilde{N}_f},
\end{equation}
and hence $Mkn \in {\tilde{N}_f}\Z$. Since $M$ and $\tilde{N}_f$ are relatively prime this implies
$kn = l \tilde{N}_f$, for some $l\in \Z$. Define $f$ by $f=\hcf(n,\tilde{N}_f)$,
$n_f$ by $n=f n_f$ and $\tilde{N}_{ff}$ by $\tilde{N}_f=f \tilde{N}_{ff}$.
Hence,
$$k=a \tilde{N}_{ff}, \quad l=a n_f  \text{\  for some } a\in \Z.$$
Substituting for $k$ in (\ref{p:eqn}) gives us
\begin{equation}
\label{p:eqn:2}
f p_2 - \frac{m}{n_f}p_1 = aM,
\end{equation}
and hence $m p_1 \in n_f \Z$. Since $m$ and $n_f$ are relatively prime this implies
\begin{equation}
\label{p1:solution}
p_1 = bn_f \quad \text{for some \ } b\in \Z.
\end{equation}
Substituting for $p_1$ in (\ref{p:eqn:2}) gives
\begin{equation}
\label{p2:solution}
p_2 = \frac{aM + mb}{f}.
\end{equation}
Substituting for $p_1$ and $p_2$ in (\ref{s:p:eqn}) gives
\begin{equation}
\label{s:p:solutions}
s = \frac{1}{f} \left( 2bn\pi - Na \theta_1(T)\right) \quad \text{and} \quad p = \frac{1}{f}Na, \quad \text{where \ }a,b\in \Z.
\end{equation}

(\ref{s:p:solutions}) gives the most general solution of the periodicity conditions.
In fact, it is simple to check that if we choose the initial factor $N_f$ of $N$ to be $1$, so that $\tilde{N}_f=N$ and
$f=\hcf(n,N)$, then this still generates all solutions of the periodicity conditions.
It follows immediately that the period lattice is of the form claimed.
\end{proof}

In the proof of Theorem \ref{general:double:period} we saw that $\Theta = \theta_2(T_J) - \alpha \theta_1(T_J)$ had limiting
values of $(1+\alpha)\pi$ and $2\pi \sqrt{\frac{1+\alpha+\alpha^2}{3}}$ as $J\ra 0$ and $J\ra \frac{1}{3\sqrt{3}}$ respectively.
Computer experiments suggest the following:
\begin{conjecture}
\label{theta:conjecture}
For each $\alpha\in (0,1)$ $\Theta$ is an increasing function of $J\in(0,\frac{1}{3\sqrt{3}})$ with
$\range{\Theta} = \left( \pi(1+\alpha),2\pi \sqrt{\frac{1+\alpha+\alpha^2}{3}}\right)$.\\
\end{conjecture}
\noindent
Two consequences of this conjecture would be:

(1) for any $\alpha \in (0,1)$ we have $\range{\Theta} \subset (\pi,2\pi)$.

(2) as $\alpha\ra 1$, the interval $\range{\Theta}$ shrinks to the point  $\{2\pi\}$.\\

\noindent
(1) implies, for example, that the integer $N$ appearing in Theorem \ref{general:double:period} must be greater than or equal to three.
(2) implies that as $\alpha \ra 1$, the integer $N$ appearing in Theorem \ref{general:double:period} must tend to infinity.
More generally, for any fixed $\alpha$, by determining  the rational number with the smallest denominator in the
interval $( \frac{1}{2}(1+\alpha), \frac{1}{\sqrt{3}}\sqrt{{1+\alpha+\alpha^2}})$ one obtains restrictions on values of
the integer $N$ which may occur for that value of $\alpha$. For example, if $\alpha = \frac{1}{2}$ then
$\frac{1}{2\pi}\Theta\in ( \frac{3}{4}, \sqrt{\frac{7}{12}}\,)$. Since the rational number with the smallest denominator in this interval is
$\frac{14}{19}$, for $\alpha = \frac{1}{2}$ we have $N\ge 19$ in Theorem \ref{general:double:period}.
Similarly, if $\alpha = \frac{1}{3}$ then $\frac{1}{2\pi}\Theta\in ( \frac{2}{3}, \sqrt{\frac{13}{27}})$.
Since the rational number with the smallest denominator
in this interval is $\frac{9}{13}$, for $\alpha=\frac{1}{3}$ we have $N \ge 13$ in Theorem \ref{general:double:period}.


\subsection{Refined lower bounds for the Legendrian and stability indices}
\label{s1inv:indexbounds}
In this section we use the explicit descriptions of the $S^1$-invariant minimal Legendrian tori
just exhibited to prove refined lower bounds for the Legendrian index and stability index of these examples.

The basic strategy is simple. Any \slg cone in $\C^n$ has distinguished homogeneous harmonic functions
of order $1$ and $2$ coming from the moment maps of translations and $\mathfrak{su}(n)$ vector fields respectively.
This implies that the restriction of any linear function on $\C^3$ to a minimal Legendrian surface $\Sigma^2$ in $S^5$
is an eigenfunction of $\Delta_{\Sigma}$, with eigenvalue $2$. Similarly, the restriction to $\Sigma$
of certain quadratic functions on $\C^3$ give eigenfunctions of $\Delta_{\Sigma}$ with eigenvalue $6$.
Since the $S^1$-invariant examples
have such explicit descriptions we can estimate the
number of nodal domains of these special eigenfunctions. The Courant Nodal Domain Theorem (Thm \ref{courant:nd})
will then guarantee that
if one of these eigenfunctions has at least $k+1$ nodal domains, then the corresponding eigenvalue is
at least the $k$-th eigenfunction of $\Delta_{\Sigma}$. In this way we get lower bounds for the number of
eigenvalues less than $2$ and less than $6$ respectively, and hence for the Legendrian and stability indices.

First we fix a convenient choice of basis for the space of quadratic functions on $\C^3$
which are moment maps of $\mathfrak{su}(3)$ vector fields.
The quadratic functions in question are the harmonic Hermitian quadratics of (\ref{hermitian:harmonic})
which have no constant nor linear part. Equivalently, given any $A\in \mathfrak{su}(3)$ the
corresponding quadratic function $Q_A$ on $\C^3$ is given by
\begin{equation}
\label{sun:quadratics}
Q_A(z)=\omega(Az,z).
\end{equation}
Define six quadratic functions $G_{ij}$ and $H_{ij}$ on $\C^3$ by
\begin{equation}
\label{gij:hij}
G_{ij} = \Real(z_i \bar{z}_j), \quad\quad H_{ij} = \Imag(z_i \bar{z}_j), \quad\quad \textrm{for\ } i< j \in \{1,2,3\}
\end{equation}
where $z_i$ denotes the $i$-th component of $z \in \C^3$.
The six functions $G_{ij}$, $H_{ij}$ together with the two functions $F_1$ and $F_2$ on $\C^3$ defined by
\begin{equation}
\label{fi}
 F_1 = |z_1|^2 - |z_3|^2, \quad F_2 = |z_2|^2 - |z_3|^2
 \end{equation}
are all of the form (\ref{sun:quadratics}) for some $A\in \mathfrak{su}(3)$ and give a basis for all such quadratics.
In particular, the restriction of any of the eight quadratic functions $F_i$, $G_{ij}$, $H_{ij}$
to any minimal Legendrian surface $\Sigma$ gives an eigenfunction of $\Delta_{\Sigma}$
with eigenvalue $6$.
As a basis for the linear functions on $\C^3$, we choose $\Real(z_i)$ and $\Imag(z_i)$ for $i=1, 2$ or $3$.
Each of these six linear functions gives an eigenfunction of eigenvalue $2$ when restricted to any minimal Legendrian surface $\Sigma$.


%


We first prove a result for the simplest $S^1$-invariant minimal Legendrian tori: the countable
family $T_{m,n}$ constructed in Theorem \ref{u:alpha:0}.
\begin{mtheorem}
\label{lindex:tmn}
\noindent
Let $T_{m,n}$ be one of the countably infinite family of minimal Legendrian $2$-tori described in Theorem \ref{u:alpha:0}. For $\lambda >0$, denote by $N_{m,n}(\lambda)$
the number of eigenvalues of $\Delta$ on $T_{m,n}$ in the range $[0,\lambda(\lambda+1)]$.
Then
\begin{equation*}
 N_{n,m}(1)   \ge \left\{ \begin{array}{ll}
                             4n+5  \ge 13 & \text{if mn is even}\\
                             2n+5  \ge 11 & \text{if mn is odd}.\\
                             \end{array}
                             \right.
\end{equation*}
\begin{equation*}
 \sind(T_{m,n}) \ge \left\{ \begin{array}{ll}
                                 8n + 4m -8  \ge 12 & \text{if mn is even}\\
                                 4n + 2m -8 \ge 6 & \text{if mn is odd}.\\
                                 \end{array}
                                 \right.
\end{equation*}
\begin{equation*}
 \lind(T_{m,n})  \ge \left\{ \begin{array}{ll}
                                 8n+4m-2  \ge 18 & \text{if mn is even}\\
                                 4n+2m-2 \ge 12 & \text{if mn is odd}.\\
                                 \end{array}
                                 \right.
\end{equation*}
Further, we also have the following upper bounds
\begin{equation*}
N_{n,m}(2) < \left\{
    \begin{array}{ll}
    36 \frac{\pi}{\sqrt{3}}n -1 & \text{if $mn$ is even,}\\
    18 \frac{\pi}{\sqrt{3}}n -1 & \text{if $mn$ is odd.}
\end{array}
\right.
\end{equation*}
\end{mtheorem}

\begin{proof} First we obtain the lower bounds.
We will consider the two cases $mn$ is even and $mn$ is odd separately.
In both cases the immersion $u_{\frac{m}{n},0}$ is given explicitly \cite[Prop 4.3]{haskins:slgcones} by
\bea
u_{\frac{m}{n},0}(s,t)&=&e^{As} \left(R_1(t), iR_2(t), R_3(t) \right) \nonumber \\
& = & e^{As} \left(c_1\cn{(rt,k)}, ic_2\sn{(rt,k)}, c_3\dn{(rt,k)} \right) \nonumber
\eea
where
$$ c_1 = \sqrt{\frac{\mu_2-\mu_1}{3\mu_2}}, \quad c_2 = \sqrt{\frac{\mu_1-\mu_2}{3\mu_1}}, c_3 = \sqrt{\frac{\mu_2-\mu_3}{3\mu_2}}$$
and
$$ \vec{\mu} = \frac{1}{n} \left( -n-2m, 2n+m, m-n \right), \quad r^2 = \frac{n+2m}{n}, \quad k^2 = \frac{n^2-m^2}{n(n+2m)}.$$

From this explicit formula for $u_{\frac{m}{n},0}$ we see immediately
that the restrictions to $T_{m,n}$ of the six quadratic functions $G_{ij}$ and $H_{ij}$ defined by (\ref{gij:hij})  are
\bea
G_{12} &= R_1 R_2 \cos{\mu_3 s},\quad   H_{12} &= - R_1 R_2 \sin{\mu_3 s} \nonumber\\
G_{13} &= R_1 R_3 \cos{\mu_2 s}, \quad   H_{13} &=  R_1 R_3 \sin{\mu_2 s} \nonumber\\
G_{23} &= R_2 R_3 \cos{\mu_1 s}, \quad   H_{23} &= - R_2 R_3 \sin{\mu_1 s}. \nonumber
\eea
\textit{Case 1}: $mn$ is even. In this case Theorem \ref{u:alpha:0} implies that the period lattice of $T_{m,n}$ is rectangular with basis
$\sigma_1 = (2n\pi,0)$ and $\sigma_2=(0,\frac{4K(k)}{r})$.
It follows that the nodal set $N_{12}$ of $H_{12}$ within the period rectangle $\mathcal{R}:= [0,2n\pi] \times [0,4K/r]$ is the following
union of horizontal and vertical lines
$$ N_{12} = \left\{(s,t) \in \mathcal{R}\ :\  t=\frac{aK}{r}, a=0,\ldots,4 \quad s=\frac{bn\pi}{n-m},\  b=0, \ldots, 2n-2m\right\}.$$
$N_{12}$ divides $\mathcal{R}$ into $8n-8m$ identical subrectangles. Since $G_{12}$ differs from $H_{12}$ only by a translation
of $\frac{\pi}{2}$ in the $s$ variable, the nodal set of $G_{12}$ is a translate of $N_{12}$.
Similarly the nodal sets for $H_{13}$ and $H_{23}$ are the following union of horizontal and vertical lines
$$N_{13} = \left\{(s,t) \in \mathcal{R}\ :\  t=\frac{aK}{r}, a=1,3 \quad s=\frac{bn\pi}{2n+m},\  b=0, \ldots, 4n+2m\right\}$$
and
$$N_{23} = \left\{(s,t) \in \mathcal{R}\ :\  t=\frac{aK}{r}, a=0,2,4 \quad s=\frac{bn\pi}{n+2m},\  b=0, \ldots, 2n+4m\right\}$$
respectively. $N_{13}$ and $N_{23}$ divide $\mathcal{R}$ into $8n+4m$ and $4n+8m$ identical subrectangles respectively. Again
the nodal sets for $G_{13}$ and $G_{23}$ are translates of the sets $N_{13}$ and $N_{23}$.
Since $n>m>0$, $N_{13}$ gives us the most nodal domains, namely $8n+4m$ of them.

Similarly, one can examine the nodal domains of the coordinate functions $\Real(z_i)$ and $\Imag(z_i)$, since
these are eigenfunctions of $\Delta$ with eigenvalue $\lambda=2$.
The zero sets are again easy to describe since each coordinate function is the product of a Jacobi elliptic function
with either $\cos$ or $\sin$. Once again the nodal domains are unions
of horizontal and vertical lines dividing the original rectangle $\mathcal{R}$ into a number of identical subrectangles. Also
the nodal domains for $\Real(z_i)$ and $\Imag(z_i)$ again differ only by a translation. The number of nodal domains
for $\Real(z_i)$ turns out to be $4n$, $4m$ and $2m+2n$ for $i=1,2,3$ respectively (remember that $\dn$ is positive, while both
$\cn$ and $\sn$ have zeroes). Hence $\Real(z_1)$ gives us the
most nodal domains, namely $4n$ of them.

Let $I$ denote the unique index so that $\lambda_{I-1}<2$ and $\lambda_I=2$,
and $J$ denote the unique index so that $\lambda_{J-1}<6$ and $\lambda_J=6$.
Then $\phi_I = \Real(z_1)$ and $\phi_J = \Real(z_1 \bar{z}_3)$ are eigenfunctions corresponding to $\lambda_I$ and $\lambda_J$
respectively.
Since $\phi_I$ has $4n$ nodal domains and $\phi_J$ has $8n+4m$ nodal domains,
by the Courant Nodal Domain Theorem (Thm \ref{courant:nd}), we must have $I\ge 4n-1$ and $J\ge 8n+4m-1$.
Now because $\lambda_{J-1}<6$, we have $\lind(T_{m,n}) \ge 8n+4m-2$.
Since $\lambda=6$ occurs with multiplicity at least $7$ we obtain $N_{n,m}(2) \ge 1 + (8n+4m-2) + 7 = 8n+4m+6$.
Hence from definition (\ref{s-index}) we have $\sind(T_{m,n}) = N_{n,m}(2) - 1 - 6 - (8-1) \ge 8n+4m-8$.
Similarly, since $\lambda=2$ occurs with multiplicity at least $6$ we obtain $N_{n,m}(1) \ge 1 + (4n-2) + 6 = 4n+5$.

\textit{Case 2}: $mn$ is odd. Theorem \ref{u:alpha:0} implies that the period lattice of $T_{m,n}$ is
generated by $\sigma_1=(2n\pi,0)$ and $\sigma_2 = (n\pi,T)$, where $T=\frac{2K}{r}$ is the basic period of $\gamma$.
Another basis for the period lattice is $\tilde{\sigma_1} = (0,2T)$, $\tilde{\sigma_2}=\sigma_2$.
We shall take the period parallelogram $\mathcal{R}$ to be the one with vertices $0$, $\tilde{\sigma_1}$, $\tilde{\sigma_2}$
and $\tilde{\sigma_1}+\tilde{\sigma_2}$.
The nodal sets $N_{ij}$ in this case, are again composed of unions of horizontal and vertical line segments contained in $\mathcal{R}$.
The two main differences from the previous case are that (i) the union of the line segments
$rt=2K$, $rt=4K$ contained in $\mathcal{R}$, forms a single closed curve in $T_{m,n}$, rather than two parallel closed
curves as in Case 1, and (ii) we are now only interested in $s\in [0,n\pi]$ and not $s\in[0,2n\pi]$.
The same holds for the union of the line segments $rt=K$, $rt=3K$, $rt=5K$ that are contained in $\mathcal{R}$.
It is not difficult to check that the number of nodal domains for the functions $H_{ij}$ and hence also $G_{ij}$
is half that of Case 1. Hence the claimed lower bounds follow from another application of the Courant Nodal Domain Theorem.\\

It still remains to prove the upper bounds claimed. For any compact minimal Legendrian surface $\Sigma$ in $S^5$,
we obtain from (\ref{area:gt:nsigma}) the following upper bound
for $N_{\Sigma}(2)$ in terms of $\Area{(\Sigma})$
$$
N_{\Sigma}(2) <  \frac{18\Area{(\Sigma})}{4\pi} - \frac{11}{7}\,.$$
Hence to prove upper bounds for $N_{n,m}$ it suffices to find upper bounds for $\Area(T_{m,n})$.
Properties of the area of $S^1$-invariant minimal Legendrian $2$-tori are studied in \S \ref{areas}.
In particular, from Theorem \ref{area:tmn} in \S \ref{areas} we have
$$
\frac{1}{4\pi}\Area(T_{m,n}) < \left\{
        \begin{array}{ll}
            \frac{2\pi}{\sqrt{3}}n & \text{if $mn$ is even,}\\
            \frac{\pi}{\sqrt{3}}n & \text{if $mn$ is odd.}
        \end{array}
\right.
$$
The claimed upper bounds now follow immediately.

By making the obvious modifications in the proof of Theorem \ref{area:spectral:genus} one can also
find explicit constants $C_1>0$ and $C_2\in \R$ so that
$N_{\Sigma}(1) < C_1 \Area(\Sigma) + C_2$
holds for any compact minimal Legendrian surface $\Sigma$ in $S^5$.
For example, one can show that $C_1 = \frac{13}{8\pi}$ and $C_2 = -\frac{7}{6}$ work.
Hence the area estimates for $T_{m,n}$
also lead to upper bounds for $N_{n,m}(1)$ improving the one coming from $N_{n,m}(1) < N_{n,m}(2)$.
\end{proof}

\begin{remark}
One could also look at the nodal sets of the quadratic functions which
correspond to elements in the diagonal subgroup of $\mathfrak{su}(3)$,
\textit{i.e.}
$Q_{\vec{a}}(z):=\sum{a_i |z_i|^2} \textrm{\ for any \ }\vec{a}=(a_1,a_2,a_3) \textrm{\  such that \ } \vec{a}.\vec{1}=0.$
If we take $\vec{a}=\vec{\lambda}=(1,\frac{m}{n},-1-\frac{m}{n})$
to be the vector determining the $1$-parameter symmetry group of the surface,
it follows from (\ref{Ri}) that $Q_{\vec{\lambda}} \equiv 0$ on any such surface.
If instead we take $\vec{a} = \vec{\mu} = \vec{1} \times \vec{\lambda}$, then again by (\ref{Ri})
$Q_{\vec{\mu}} = \sum{\mu_i |z_i|^2} = \sum{\left(\mu_i^2 \gamma + \frac{1}{3}\mu_i\right)} = \gamma \left(\sum{\mu_i^2}\right).$
$\gamma$ is periodic of period $\frac{2K}{r}$ and has exactly one zero $t_0\in [0,\frac{K}{r}]$ and another zero
$t_1=\frac{2K}{r}-t_0$ in $[\frac{K}{r},\frac{2K}{r}]$. From this it is easy
to see that when $mn$ is even, $Q_{\vec{\mu}}$ has exactly $4$ nodal domains on the period rectangle
$\mathcal{R}$ independent of $m$ and $n$. When $mn$ is odd, the union of the line segments
$rt=rt_i$, $rt=2K+rt_i$, $rt=4K+rt_i$ forms a single closed curve in $T_{m,n}$ for either $i=0$ or $i=1$.
Hence in this case $Q_{\vec{\mu}}$ has only $2$ nodal domains independent of $m$ and $n$.
\end{remark}


\begin{remark}
\label{ac:dimension:unbounded} (At this point the reader is advised to consult Appendix \ref{acslg:section} on
asymptotically conical \slg submanifolds.)
By applying Lemma \ref{rotated:cone} to the cone $C_{m,n}=C(T_{m,n})$ and setting $t=1$ we obtain a two-ended \acslg
submanifold ${L_{m,n}}$ in $\C^3$ with rate $\lambda=-1$ asymptotic to $C:= C_{m,n}\cup e^{i\frac{\pi}{3}}C_{m,n}$.
Because $C$ is invariant under an $S^1$-subgroup of $\text{SU(3)}$, by Proposition \ref{acslg:symmetry} so is $L_{m,n}$.
Since
${L_{m,n}}$ has rate $\lambda=-1$ it also
has any rate $\lambda\in (-1,2)$. Hence as in Appendix \ref{acslg:section}, for any $\lambda \in (-1,2)$
we may consider $\mathcal{M}_{L_{m,n}}^\lambda$, the set of \acslg submanifolds in $\C^3$ with
rate $\lambda$ and cone $C$ which are isotopic to ${L_{m,n}}$ as an \ac submanifold of $\C^3$.
Since the link of $C$, $\Sigma = T_{m,n} \cup e^{i\frac{\pi}{3}} T_{m,n}$, the set of homogeneous harmonic functions on $C$
is two copies of the harmonic homogeneous functions on $C_{m,n}$. Hence $\dsig = \mathcal{D}_{T_{m,n}}$ and
$N_{\Sigma}(\lambda) = 2N_{T_{m,n}}(\lambda)$.
By Theorem \ref{acslg:deformations}, if $\lambda\in (0,2) - \dsig$ then $\mathcal{M}_{L_{m,n}}^{\lambda}$ is a manifold
with $\dim{\mathcal{M}_{L_{m,n}}^{\lambda}} = b^1({L_{m,n}})-b^0({L_{m,n}}) + N_\Sigma(\lambda)$.
In particular, for any $\lambda\in (1,2) - \dsig$ applying Theorem \ref{lindex:tmn}
we have $\dim{\mathcal{M}_{L_{m,n}}^{\lambda}} \ge b^1({L_{m,n}})-b^0({L_{m,n}}) + N_\Sigma(1) \ge b^1({L_{m,n}})-b^0({L_{m,n}})+4n+5$.

In particular, we have a countably infinite family of \acslg submanifolds $L_{m,n}$ of $\C^3$ each diffeomorphic to $T^2\times \R$
which come in moduli spaces of unbounded dimension. Hence even in this simple case the topology of an \acslg submanifold does
not come close to controlling its geometry.
\end{remark}

The next result uses the method of Theorem \ref{lindex:tmn} to
prove improved bounds for the l-index of the next simplest infinite family of
$S^1$-invariant minimal Legendrian $2$-tori: the examples constructed in Theorem \ref{u:zero:j}.
\begin{mtheorem}
\label{lindex:u0j}
Let $\Sigma$ be any of the countably infinite family of minimal Legendrian $2$-tori constructed in
Theorem \ref{u:zero:j}. Then $\frac{1}{2\pi}\theta_2(T_J) = \frac{m}{n}$ for relatively prime
integers $m\ge 4$ and $n\ge 7$, and
\begin{eqnarray*}
\sind(\Sigma) & \ge & 2n-8 \ge 6.\\
\lind(\Sigma) & \ge & 2n-2 \ge 12.\\
N_{\Sigma}(1) & \ge & 2m +5 \ge 13.\\
N_{\Sigma}(2) & \ge & 2n +6 \ge 20.
\end{eqnarray*}
\end{mtheorem}

\begin{proof}
By Prop \ref{double:period:prop}(d) the necessary and sufficient condition that $u_{0,J}$ closes to give a torus is $\theta_2(T_J)\in 2\pi \Q$.
As remarked above Theorem \ref{u:zero:j}, $\theta_2(T_J)$ is a monotone continuous function of $J \in (0,\frac{1}{3\sqrt{3}})$
with $\range{\theta_2(T_J)} = (\pi, \frac{2\pi}{\sqrt{3}})$.  Hence we are interested in rational numbers in the interval $(\frac{1}{2},\frac{1}{\sqrt{3}})$.
The rational number with the smallest denominator in this interval turns out to be $\frac{4}{7}$. Hence $m\ge 4$ and $n\ge 7$ as claimed.

Any torus constructed in Theorem \ref{u:zero:j} has the form
$$u_{0,J}(s,t) = (e^{is}R_1e^{i\theta_1}, R_2 e^{i\theta_2}, e^{-is}R_1 e^{i\theta_1})$$
for some $J\in (0,\frac{1}{3\sqrt{3}})$.
Because of the form of $u$ there are three
$\mathfrak{su}(3)\ltimes \C^3$ moment map functions which are particularly straightforward to
analyze: $\Imag(z_2)$, $Q_{\vec{\mu}}=\sum{\mu_i|z_i|^2} = (\sum{\mu_i^2})\gamma$ and $\Imag(z_1 \bar{z}_3)$.

When $m$ is even the period lattice  of $\Sigma$ is rectangular with basis
$\omega_1=(2\pi,0)$ and $\omega_2=(0,nT)$ and otherwise has basis $\omega_1=(2\pi,0)$ and $\omega_2=(\pi,nT)$.
We will take a fundamental domain to be the region $\mathcal{R}$ with vertices $0$, $\omega_1$, $\omega_2$ and $\omega_1+\omega_2$.
Consider first the case where $m$ is even and hence $\mathcal{R}$ is rectangular.
Since $\Imag(z_1 \bar{z}_3) = R_1^2(t)\sin{(2s)}$ and $R_1(t)>0$, it is follows that $\Imag(z_1 \bar{z}_3)$
vanishes in $\mathcal{R}$ if and only if $s$ is an integer multiple of $\frac{\pi}{2}$ contained in the interval $[0,2\pi]$.
Whether $m$ is even or odd, the number of nodal domains
of $\Imag(z_1 \bar{z}_3)$ is at most $4$ independent of $m$ and $n$. The other two functions give us more useful information.

To prove the bound for $N_{\Sigma}(1)$, we consider $\Imag(z_2) = R_2(t) \sin{\theta_2(t)}$, which is a function of $t$ alone. Since $R_2(t)>0$ for all $t$, its zeroes
in $\mathcal{R}$ correspond to the solutions of $\sin{\theta_2(t)}=0$ in $[0,nT]$. Because $\theta_2(t)$ is a strictly increasing function
of $t$ and $\theta_2(nT) = n \theta_2(T)= 2\pi m$, it follows by continuity of $\theta_2$ that it has exactly $2m+1$
zeroes in $[0,nT]$, two of which are $t=0$ and $t=nT$. Hence the nodal set is a union of $2m$ lines each of the form
$t=c_i$ for some constant $c_i$. It follows that there are exactly $2m$ nodal domains.
Since $\Imag(z_2)$ is an eigenfunction of $\Delta_\Sigma$ with $\lambda=2$, applying the Courant Nodal Domain Theorem
(Thm \ref{courant:nd}) and arguing as in the proof of Theorem \ref{lindex:tmn} shows that there are at least
$2m-2$ nonzero eigenvalues below $2$. Hence $N_{\Sigma}(1) \ge 1+(2m-2) + 6 = 2m+5$.

To prove the remaining three lower bounds we consider $Q_{\vec{\mu}}= (\sum{\mu_i^2})\gamma$, which again is a function of $t$ alone.
Since $\gamma$ is periodic of period $T$ and has exactly $2$ zeroes on $[0,T]$ both of which are in the interior,
$\gamma$ has exactly $2n$ zeroes on $[0,nT]$
all of which are in the interior of the interval. It follows that there are exactly $2n$ nodal domains.
Since $Q_{\vec{\mu}}$ is an eigenvalue of $\Delta_{\Sigma}$ with eigenvalue $6$, applying
the Courant Nodal Domain Theorem again shows that there are at least $2n-2$ eigenvalues below $6$, \textit{i.e.}
$\lind(\Sigma)\ge 2n-2$. The two remaining lower bounds now follow easily.

In the case when $m$ is odd and hence $\mathcal{R}$ is not rectangular, it is not difficult to check
that since both eigenfunctions we considered depended only on $t$ that we obtain the same number of nodal domains
as in the rectangular case.
\end{proof}

\begin{remark}
The first three rational numbers in the interval $(\frac{1}{2},\frac{1}{\sqrt{3}})$
written in the form $\frac{m}{n}$ for relatively prime $m$ and $n$ and
ordered by $m+n$ are $\frac{4}{7}, \ \frac{5}{9}$, and $\frac{6}{11}$.
In  Table 2 (\S \ref{areas}) we list the area and index bounds of these tori. The area is obtained by
computing numerically the values of certain elliptic integrals.
\end{remark}

For the general doubly periodic $u_{\alpha,J}$ constructed in Theorem \ref{general:double:period}
one can also prove lower bounds for the Legendrian index and the stability index in terms
of the integers $n$ and $N$ appearing in the statement of Theorem \ref{general:double:period}.
The simplest eigenfunction to analyze is $Q_{\vec{\mu}}= \sum{\mu_i|z_i|^2}=(\sum{\mu_i^2})\gamma$,
since it depends only on $t$. In this case the best results are obtained when $n$ is prime
and hence $\hcf(n,N)=1$. Establishing Conjecture \ref{theta:conjecture} would be useful here. It would
give lower bounds on the possible $N$, for a given choice of $\alpha=\frac{m}{n}$.
We leave the details of the analysis which is very similar to the previous two cases to the interested reader.

\section{$S^1$-invariant \slg cones with small area}
\label{areas}

In this section we establish a formula for the area of an $S^1$-invariant minimal Legendrian $2$-torus
in terms of elliptic integrals. Using some basic properties of elliptic integrals one can easily
evaluate certain limiting behaviour and in some cases deduce some useful monotonicity properties.
These monotonicity properties then imply upper and lower bounds for the area of certain
$S^1$-invariant tori. As well as its intrinsic geometric interest,
the area upper bound is particularly useful because it implies upper bounds for the
Legendrian index and stability index. These complement the lower bounds proved by the nodal domain
counting technique of the previous section. A lower bound for area implies via (\ref{yang:yau})
an upper bound for $\lambda_1$.

Since all the immersions $u_{\alpha,J}$ are conformal, the metric $g$ induced
on $\R^2$ by the immersion is determined by a single positive function $y(t)$
where $g = y(t) (ds^2 + dt^2)$. A calculation shows that $\gamma$ and $y$ are related
by
\beq
\label{y:gamma}
y = - \gamma \mu_1 \mu_2 \mu_3 + {\textstyle{\frac{1}{3}\sum{\lambda_i^2}}}.
\end{equation}
It follows from (\ref{gammadot}) and (\ref{y:gamma}) that $y$ satisfies
\beq
\label{ydot}
\dot{y}^2 + 4y^3 -2y^2 {\textstyle{\sum{\lambda_i^2}}} = -4(\lambda_1^2 \lambda_2^2 \lambda_3^2 + J^2 \mu_1^2 \mu_2^2 \mu_3^2).
\end{equation}
When $\dot{y}=0$, (\ref{ydot}) has three solutions $y_2 \le 0 \le y_1 \le y_3$.
$y$ can also be expressed in terms of $\sn$ as
\beq
\label{y:elliptic}
y = y_3 - (y_3 -y_1) \sn^2{(rt,k)}
\end{equation}
where
\beq
\label{rk:y}
 r^2 = y_3-y_2\, ,\quad k^2 = \frac{y_3 -y_1}{y_3 - y_2}.
\end{equation}
Define  $A(\alpha,J)$ by
\begin{equation}
\label{A:alpha:J}
A(\alpha,J) = \int_0^{T_{\alpha,J}}{y_{\alpha,J}(t) dt}
\end{equation}
\textit{i.e.} $A$ is the integral of the conformal factor $y$ over one basic period of $y$.

\begin{prop} Let $K$ and $E$ be the complete elliptic integrals of the first and second kinds respectively,
as defined in (\ref{ellipK}) and (\ref{ellipE}).
Then $A(\alpha,J)$ can be expressed in terms of $K$ and $E$ as
\beq
\label{area:period}
\tfrac{1}{2}A(\alpha,J) = \frac{y_2 K(k)}{r} + r E(k).
\end{equation}
\end{prop}

\begin{proof}
Using the explicit expression for $y$ in terms of $\sn{}$ we have
\bea
\label{area1}
A & = & \int_0^{\frac{2K(k)}{r}}{y_3 - (y_3-y_1)\sn^2{(rt,k)}}\ dt
= \frac{2y_3K}{r} - \frac{2}{r}(y_3-y_1) \int_0^{K}{\sn^2{(t,k)}\ dt}.
\eea
But using $\dn^2t=1-k^2\sn^2t$, and the definition of $E(k)$ in terms of
$\dn$, (\ref{ellipE}),  we have
$$
\int_0^K {\sn^2{t}\ dt} = \frac{1}{k^2}\int_0^K{1-\dn^2{t}\ dt} =
\frac{1}{k^2}(K-E).$$
Inserting this expression into (\ref{area1}) and simplifying we have
\begin{equation}
\label{area2}
\tfrac{1}{2} A=  \frac{y_3}{r}K - r(K-E) = \frac{y_2}{r}K + rE,
\end{equation}
where the final equality follows from $r^2=y_3-y_2$.
\end{proof}

\begin{prop}[The limiting cases $\alpha=1$, $J\ra \frac{1}{3\sqrt{3}}$ and $J=0$]
\label{area:limits}
\quad \\
(a) For $\alpha=1$ we have
$$ T(J) \equiv \frac{\pi}{\sqrt{3}}, \quad A(J) \equiv \frac{2\pi}{\sqrt{3}}, \quad \theta_3 \equiv 0.$$
(b) In the limit $J \ra  \frac{1}{3\sqrt{3}}$ we have
$$ T_{\alpha} = \frac{\pi}{\sqrt{1+\alpha+\alpha^2}}\ , \quad A(\alpha) = \frac{2\pi}{3}\sqrt{(1+\alpha+ \alpha^2)}\ , \quad
\theta_i(T_\alpha) = \frac{\pi \mu_i}{\sqrt{3(1+\alpha+\alpha^2)}}.$$
(c) In the case $J=0$ the function $A(\alpha)$ defined by (\ref{A:alpha:J}) is an increasing function of $\alpha$
with $\range(A)=[2,\frac{2\pi}{\sqrt{3}}]$.
\end{prop}

\begin{proof}
Parts (a) and (b) are simple calculations which we omit. Part (c) is the nontrivial part of the Proposition and was first observed in numerical
experiments. Its proof is also a (not so straightforward) computation, which makes some use of some standard facts from the
theory of elliptic integrals about derivatives with respect to $k$ of $E$ and $K$.

%
\no

The basic strategy is simple. First we express $y_1, y_2, y_3, k$ in terms of $r$, where
$r$ is defined in (\ref{rk:y}). Then we use
(\ref{area:period}) together with some standard expressions for $\frac{dE(k)}{dk}$ and $\frac{dK(k)}{dk}$ to prove that
$\frac{dA}{dr} > 0$. The detailed calculation is somewhat involved, and some details will be suppressed
below. Since $J=0$, it follows from (\ref{ydot}) that $y_1, y_2, y_3$ satisfy the equation
\beq
y^3 - py^2 + (p-1)^2 =0
\end{equation}
where $p:=1+\alpha+\alpha^2$. It is easy to verify that $y_1=p-1$ is always a root. If follows that
$y_2+y_3=1$ and $y_2y_3 = (1-p)$. Since from (\ref{rk:y}) we also have $r^2=y_3-y_2$
this implies that
\beq
\label{y23}
y_2 = \frac{1}{2}(1-r^2),\quad  y_3=\frac{1}{2}(1+r^2).
\end{equation}
and hence
\beq
\label{y1:r}
 y_1 = \frac{1}{4}(r^4-1)
\end{equation}
where in the last equality we use the fact that $p$ and $r$ are related by
\begin{equation}
\label{p:r}
p = \frac{1}{4}(r^4+3).
\end{equation}
It is easy to see that as $\alpha$ varies in $[0,1]$, $p$ varies in $[1,3]$ and $r^4$ varies in $[1,9]$.
$p$ is an increasing function of $\alpha$ in this range, $r$ is an increasing function of $\alpha$ \textit{etc.}
Substituting for $y_2$ from (\ref{y23}) into (\ref{area:period}) we obtain
$A(r) = \frac{1}{2}\left( \frac{1}{r} - r \right)K(k) + rE(k)$ and
\beq
\label{da:dr}
\frac{dA}{dr} = E + r \frac{dE}{dr} + \frac{1}{2r}(1-r^2)\frac{dK}{dr} - \frac{1}{2r^2}(1+r^2)K.
\end{equation}
Combining equations (\ref{rk:y}), (\ref{y23}) and (\ref{y1:r}) we obtain
\beq
\label{k:r}
k^2 = \frac{(1+r^2)(3-r^2)}{4r^2}, \quad {k'}^2 := 1-k^2 = \frac{(r^2-1)(r^2+3)}{4r^2}.
\end{equation}
Equations (\ref{k:r}), (\ref{dE:dk}) and (\ref{dK:dk}) imply that
\beq
\frac{dE}{dr} = \frac{dE}{dk}\frac{dk}{dr} = \frac{E-K}{k} \frac{-1}{4k}\left(\frac{3}{r^3}+r \right)
\end{equation}
and
\beq
\frac{dK}{dr} = \frac{dK}{dk}\frac{dk}{dr} = \frac{1}{k} \left(\frac{E}{{k'}^2} -K\right) \frac{-1}{4k}\left(\frac{3}{r^3}+r \right).
\end{equation}
After substituting these expressions for $\frac{dE}{dr}$ and $\frac{dK}{dr}$ into (\ref{da:dr})
and performing some algebraic manipulation one finds
\beq
\frac{dA}{dr} = \frac{r^2-1}{3-r^2} K + \frac{2(3-r^4)}{(3+r^2)(3-r^2)}E.
\end{equation}
Since $K$ and $E$ are both positive it follows immediately that $\frac{dA}{dr} >0$ for $r^4\in (1,3)$, but
is not clear that this still holds for $r^4\in [3,9)$. However, the previous expression for $\frac{dA}{dr}$
can be rewritten as
\beq
\label{da:dr:2}
\frac{dA}{dr} = \frac{1+r^2}{3+r^2} E + \frac{r^2-1}{3-r^2}\left(K-E \right).
\end{equation}
Since $K-E>0$, $E>1$ both hold for $k\in(0,1)$, (\ref{da:dr:2}) implies that
$$
\frac{dA}{dr} > \frac{1+r^2}{3+r^2} = 1 - \frac{2}{3+r^2} > \frac{1}{2}
$$
for $r^2 \in (1,3)$. Hence $A$ is an increasing function of $r$ for $r^4 \in (1,9)$
and hence also an increasing function of $\alpha$ for $\alpha \in (0,1)$.
All that remains is to find the values $A(0)$ and $A(1)$.
This is a short calculation that we omit.
\end{proof}

We can use part (c) of the previous proposition to estimate analytically the areas of the minimal Legendrian tori
$T_{m,n}$ from Theorem \ref{u:alpha:0}.

\begin{mtheorem}
\label{area:tmn}
Let $T_{m,n}$ one of the countably infinite family of minimal Legendrian $2$-tori constructed in Theorem \ref{u:alpha:0}.
Let $A(\alpha)$ denote the function $A(\alpha,0)$ defined in (\ref{A:alpha:J}). Then
$$\Area(T_{m,n}) = \left\{
            \begin{array}{ll}
                4n\pi A\left(\tfrac{m}{n}\right) & \quad \text{if $mn$ is even,}\\
                2n\pi A\left(\tfrac{m}{n}\right) & \quad \text{if $mn$ is odd.}
            \end{array}
            \right.
$$
Hence if $mn$ is even then  $2n < \tfrac{1}{4\pi}\Area(T_{m,n}) < \tfrac{2\pi}{\sqrt{3}}n$,
otherwise $ n < \frac{1}{4\pi}\Area(T_{m,n}) < \frac{\pi}{\sqrt{3}}n$.
\end{mtheorem}

\begin{proof}
Case 1: $mn$ is even. In this case Proposition \ref{u:alpha:0} tell us that the period lattice of $T_{m,n}$
is generated by $\sigma_1=(2n\pi,0)$ and $\sigma_2 = (0,\frac{4K(k)}{r})$. Hence the first part of the theorem follows from
$$ \Area(T_{m,n}) = \int_{s=0}^{s=2n\pi} \int_{t=0}^{t=\frac{4K(k)}{r}} y(t) ds dt = 4n\pi \int_0^{\frac{2K}{r}}y(t) dt =
4n\pi A\left(\tfrac{m}{n}\right).$$
The upper and lower bounds for the area now follow from  Proposition \ref{area:limits}(c).

Case 2: $mn$ is odd. In this case Proposition \ref{u:alpha:0} tell us that the period lattice of $T_{m,n}$
is generated by $\sigma_1=(2n\pi,0)$ and $\sigma_2 = (n\pi,T)$.
$\Area(T_{m,n})$ is equal to
the integral of $y$ over the parallelogram with vertices
$0, {\sigma_1}, {\sigma_2}$ and ${\sigma_1}+{\sigma_2}$.
Hence the first part of the theorem follows from
$$ \Area(T_{m,n}) = \int_{t=0}^{t=T} \int_{s=ct}^{s=ct+2n\pi}  y(t) ds\, dt = 2n\pi \int_{t=0}^{t=T} y(t) dt =
2n\pi A\left(\tfrac{m}{n}\right)$$
where $c=\frac{n\pi}{T}$, and $ T=\frac{2K}{r}$ is the basic period of $y$.
The upper and lower bounds for the area follow once again from Proposition \ref{area:limits}(c).
\end{proof}
\begin{remark}
Numerical evidence suggests that for any $\alpha \in (0,1)$, the function $A(\alpha,J)$ is a strictly increasing
function of $J\in (0,\frac{1}{3\sqrt{3}})$. If one can establish this fact then one will obtain
upper and lower bounds for the area of a general torus $u_{\alpha,J}$ in same way we proved Theorem \ref{area:tmn} from
Proposition \ref{area:limits}(c).
\end{remark}

Two immediate consequences of Theorem \ref{area:tmn} are the following.
\begin{corollary}[monotonicity of area for $T_{m,n}$]
Let $T_{m,n}$ one of the minimal Legendrian tori constructed in Proposition \ref{u:alpha:0}.\\
(i) If both $m_1 n$ and $m_2 n$ are even, then $m_1< m_2 \Rightarrow \Area(T_{m_1,n}) < \Area(T_{m_2,n})$.
In particular for any fixed $n$, $\Area(T_{m,2n})$ is an increasing function of $m$.\\
(ii) If both $m_1 n$ and $m_2 n$ are odd, then $m_1< m_2 \Rightarrow \Area(T_{m_1,n}) < \Area(T_{m_2,n})$.
In particular for any fixed $n$, $\Area(T_{2m+1,2n+1})$ is an increasing function of $m$.
\end{corollary}
\begin{proof}
The result follows immediately from Theorem \ref{area:tmn} and the fact that according to Proposition \ref{area:limits}(c)
$A(\alpha)$ is an increasing function of $\alpha\in [0,1]$.
\end{proof}

\begin{corollary}[first eigenvalue estimates for $T_{m,n}$]
\label{first:eval:bounds}
Let $T_{m,n}$ one of the minimal Legendrian tori constructed in Proposition \ref{u:alpha:0}.\\
(i) If $m$ and $2n$ are relatively prime then $\lambda_1(T_{m,2n}) < \frac{1}{n}$.\\
(ii) If $2m+1$ and $2n+1$ are relatively prime, then $\lambda_1(T_{2m+1,2n+1}) < \frac{4}{2n+1}$.
\end{corollary}
\begin{proof}
The result follows immediately by combining inequality (\ref{yang:yau})
with the lower bound for $\Area(T_{m,n})$ from Theorem \ref{area:tmn}.
\end{proof}

The previous theorem gives explicit expressions for $\Area(T_{m,n})$ in terms of the value of
certain elliptic integrals. It is straightforward to use any of the standard mathematics software
packages (\textit{e.g.} Mathematica or Matlab) to numerically evaluate these elliptic integrals to any
desired accuracy. In particular, it is easy to construct a list of the minimal Legendrian tori
$T_{m,n}$ with area less than some fixed area.
A list of tori $T_{m,n}$ by increasing area is given in Table \ref{table:smallareatori}.

Empirically $A(\alpha)$ is reasonably well-approximated by the straight line
$L(\alpha) = M \alpha + 2$ where  $M:=\frac{2\pi}{\sqrt{3}} - 2$
is the slope of the line segment between the points $(0,A(0))$ and $(1,A(1))$.
Hence a rough estimate of the area of the torus $T_{m,n}$ is given by
\begin{equation}
\label{area:approx}
 \frac{\Area(T_{m,n})}{4\pi} \sim  \left\{
\begin{array}{ll}
n + \frac{1}{2}Mm & \mbox{\ if $mn$ is odd}\\
2n+Mm & \mbox{\ if $mn$ is even}
\end{array}
 \right.
\end{equation}
where $M=\frac{2\pi}{\sqrt{3}} - 2 \sim 1.628$. This linear approximation
to the area of $T_{m,n}$ is also listed in Table 1.
Note that for these small values of $m$ and $n$ it produces the correct
ordering by area.

\begin{table}
\label{table:smallareatori}
\begin{tabular}{||c||c|l|l|l|l|l|l|l|l||} \hline
\mbox{\rule[-0.1in]{0cm}{0.7cm}\textbf{$(m,n)$}} & (0,1) & (1,3) & (1,2) & (1,5) & (3,5) & (1,7) & (2,3) & (3,7) & (1,4)\\ \hline
\mbox{\rule[-.4cm]{0cm}{1.15cm}$\dfrac{1}{4\pi}\Area(T_{m,n})$} & 1.81 & 3.71 & 5.49 & 5.66 & 7.29 & 7.63 & 9.10 & 9.19 & 9.36\\ \hline
\mbox{\rule[-0.1in]{0cm}{0.7cm}Linear approx using (\ref{area:approx})} & - & 3.81 & 5.63 & 5.81 & 7.44 & 7.81 & 9.26 & 9.44 & 9.63 \\ \hline
\mbox{\rule[-0.1in]{0cm}{0.7cm}$\sind$} & $0$  & $\ge 6$  & $\ge 12$ & $\ge 14$ & $\ge 18$ & $\ge 22$ & $\ge 24$ & $\ge 26$ & $\ge 28$\\ \hline
\end{tabular}\\
\medskip
\caption{Tori $T_{m,n}$ with small area.}
\end{table}

\begin{table}
\label{table:area:u:0:j}
\begin{tabular}{||c||c|c|c||} \hline
\mbox{\rule[-0.4cm]{0cm}{1.1cm}$\dfrac{\theta_2(T_J)}{2\pi}$} & $\dfrac{4}{7}$ & $\dfrac{5}{9}$ & $\dfrac{6}{11}$ \\ \hline
\mbox{\rule[-.4cm]{0cm}{1.15cm}$\dfrac{1}{4\pi}\Area$} & $7.26$ & $9.16$ & $11.12$ \\ \hline
\mbox{\rule[-0.4cm]{0cm}{1.15cm}$\dfrac{J}{3\sqrt{3}}$} & $0.831$ & $0.495$ & $0.344$ \\ \hline
\mbox{\rule[-0.2cm]{0cm}{0.7cm}$\sind$} & $\ge 12$ & $\ge 16$ & $\ge 20$\\ \hline
\mbox{\rule[-0.4cm]{0cm}{1.0cm}$\lambda_1$} & $< \dfrac{4}{7}$ & $< \dfrac{4}{9}$ & $ < \dfrac{4}{11}$\\ \hline
\end{tabular}\\
\medskip
\caption{Area, $\sind$ and $\lambda_1$ of tori $u_{0,J}$ with small areas.}
\end{table}

\section{Concluding Remarks}

We conclude with a few remarks about natural extensions of the results proved in this paper
and some related questions.

The basic strategy of \S \ref{s1:invariant:tori} for proving lower bounds for the stability index of a \slg cone
should be useful beyond the $S^1$-invariant $T^2$-cones of \S \ref{s1:invariant:tori}.
Basically, the method just needs a sufficiently explicit description of the \slg cone in question. So the results of
\S \ref{s1:invariant:tori} should be thought of as an illustration of a basic strategy. Two other natural classes
of examples to which the same kind of analysis could be applied are described below.

Explicit examples of \slg $T^2$-cones of spectral curve genus $4$ were constructed in \cite{joyce:intsys}.
These tori also have descriptions in terms of
(products of) elliptic integrals and elliptic functions.
Using the basic method of \S \ref{s1:invariant:tori} it should be possible
to prove lower bounds for the stability index and Legendrian index of these tori.
Lower bounds for the area of these tori will then follow from (\ref{area:gt:nsigma}).
As in the $S^1$-invariant case the area of these tori can
be expressed in terms of elliptic integrals. If one can use these expressions to prove upper bounds for the
area of these tori then again (\ref{area:gt:nsigma}) implies corresponding
upper bounds for the stability index. Certainly one can use the expression for the area in terms of elliptic integrals
to investigate numerically the area of these tori.

Examples of $T^{n-1}$-invariant \slg cones in $\C^n$ generalizing the $S^1$-invariant $T^2$-cones of
\S \ref{s1:invariant:tori} were constructed in \cite{joyce:symmetries}. In these examples there
is a natural analogue of the function $\gamma$ which gives us an eigenfunction of $\Delta$ of eigenvalue
$\lambda=2n$, and whose behaviour is particularly simple to analyze. So the methods of \S \ref{s1:invariant:tori}
should again lead to lower bounds for the stability index, and hence by the higher-dimensional analogue of
(\ref{area:gt:nsigma}) to lower bounds for the volume. We plan to discuss both the $T^2$-cones
of spectral curve genus $4$ and the $T^{n-1}$-invariant cones in a sequel to this paper.

The results of \S \ref{s1:invariant:tori} and \S \ref{areas} give a quite good understanding
of the behaviour of the stability index and area of $S^1$-invariant $T^2$-cones. One outstanding question
which our methods do not address is: which $S^1$-invariant \slg $T^2$-cones are rigid or Jacobi integrable?
Jacobi integrability of a cone implies, for example, extra regularity for geometric measure theoretic SLG objects
which resemble that cone.
Some more sophisticated numerical work would allow us to compute numerically at least the low end
of $\text{Spec}{(\Delta)}$ with reasonable accuracy. This should indicate which $S^1$-invariant examples
are rigid.

One could also hope to demonstrate numerically the existence of Jacobi integrable
but not rigid $T^2$-cones. To do this consider a \slg $T^2$-cone of spectral curve genus $g=2d$ for some $d>2$.
The discussion of \S \ref{area:index:genus} shows that there is at least a $6+d$-dimensional family of integrable
Jacobi fields of $T^2$, and hence that the multiplicity of the eigenvalue $\lambda=6$ is at least $6+d$.
Hence by computing numerically that the multiplicity of $\lambda=6$ is exactly $6+d$ then one establishes
integrability.
Without the aid of the spectral curve genus it is difficult to see how numerics
can help resolve whether non-rigid examples are integrable or not.
For example, it is unknown whether the Clifford torus in dimensions $8$ and $9$ is
Jacobi integrable, even though from equation (\ref{ct:spectrum}) we know the spectrum explicitly.
We know of no other general methods for proving Jacobi integrability in our setting.

We have seen that the Clifford torus $T^{n-1}$ fails to be strictly stable for $n>3$ (Prop \ref{ct:spectrum}). It is natural to
wonder what are the strictly stable SLG cones in higher dimensions. It would be natural to start by
finding the stability index of the higher-dimensional homogeneous examples.

The results of \cite{pedit:pinkall} prove that the area of a minimal torus in $S^3$ grows quadratically with spectral curve
genus. It is natural to ask if the (linear) lower bounds established in this paper for both stability index and area  can be improved to quadratic
ones.

\appendix
\section{Jacobi elliptic functions and elliptic integrals}
\label{appendix:elliptic:fn}
In this appendix we recall some basic properties of the Jacobi
elliptic functions and elliptic integrals which are needed in \S \ref{s1:invariant:tori} and \S \ref{areas}.
For a more leisurely description of elliptic functions we refer the reader to
\cite{lawden}.

We take the following definition of the Jacobi sn-noidal
elliptic function. Let $\sn{(t,k)}$, the \textit{Jacobi sn-noidal function with modulus} $k \in [0,1)$, be the unique solution of the equation
\begin{equation}
  \label{sn-noid}
  \dot{z}^2= (1 -z^2)(1-k^2 z^2)
\end{equation}
with $z(0)=0,\  \dot{z}(0)=1$. It is straightforward to see from this definition that in the case $k=0$ we have
$\sn{(t,0)}= \sin{t}$ and that for $k=1$ we get $\sn(t,1)=\tanh{t}$. By analogy with the trigonometric functions
there is also a \textit{Jacobi cn-noidal function} $\cn{(t,k)}$ which satisfies
$$ \cn^2{(t,k)} = {1 - \sn^2{(t,k)}}.$$
There is another Jacobi elliptic function $\dn{}$ satisfying
$$ \dn^2{(k,t)} = {1-k^2 \sn^2{(t,k)}}.$$
Using the definition of $\sn$ given in (\ref{sn-noid}) and the
relationships between the squares of the other Jacobi elliptic
functions  we find that
\begin{eqnarray*}
\frac{d}{dt} \sn{(t,k)} & = & \cn{(t,k)} \dn{(t,k)}\\
\frac{d}{dt} \cn{(t,k)} & = & -\sn{(t,k)} \dn{(t,k)}\\
\frac{d}{dt} \dn{(t,k)} & = & -k^2 \sn{(t,k)} \cn{(t,k)}.
\end{eqnarray*}
The period of $\sn{(t,k)}$ and $\cn{(t,k)}$ is $4K(k)$, while
$\dn{(t,k)}$ has period $2K(k)$, where $K$ is the \textit{complete
elliptic integral of the first kind} defined by
\begin{equation}
\label{ellipK}
K(k) = \int_0^{\pi/2} \frac{dx}{\sqrt{1-k^2 \sin^2{x}}}.
\end{equation}
Similarly, the \textit{complete elliptic integral of the second kind} $E$
is defined by
\begin{equation}
\label{ellipE}
E(k) = \int_0^{\pi/2} {\sqrt{1-k^2 \sin^2{x}}} \ dx = \int_0^K \dn^2{t}\ dt.
\end{equation}
From the definitions of $K$ and $E$ it is clear that for $k\in [0,1),$
$K$ is a positive strictly increasing function of $k$ and
$E$ is a positive strictly decreasing
function of $k$.

The following expressions for the derivatives of $K$ and $E$ with respect to $k\in (0,1)$ are needed in \S \ref{areas},
\begin{eqnarray}
\label{dE:dk}
\frac{dE}{dk}  & = & \frac{1}{k}(E-K),\\
\label{dK:dk}
\frac{dK}{dk} & =  & \frac{1}{k{k'}^2}(E-{k'}^2K)
\end{eqnarray}
where ${k'}^2:= 1-k^2$. See \cite[\S 3.8]{lawden} for proofs of these formulae.

\section{Applications to desingularisations of singular slg submanifolds}
\label{analysis:spectral:geometry}


A fundamental question in special Lagrangian geometry is to understand when one can desingularize
a singular \slg variety.
This question is closely related to understanding compactifications of
moduli spaces of compact nonsingular \slg submanifolds.
When the singular \slg variety has isolated conical singularities this problem
can be approached by attempting to glue in appropriate pieces of nonsingular \slg submanifolds
which are asymptotic to \slg cones. Joyce has carried out such a programme in the recent series of
papers \cite{joyce:slgcones1,joyce:slgcones2,joyce:slgcones3,joyce:slgcones4,joyce:slgcones5}.

Joyce's programme has three main components:
the deformation theory of asymptotically conical SLG submanifolds,
the deformation theory of singular \slg varieties with isolated conical singularities,
and the desingularization theory of \slg varieties with isolated conical singularities.
The spectral geometry of \slg cones plays a key part in all three components.
In this appendix we review very briefly the most relevant parts of each component in order to explain how
results of the current paper can be applied to all three components of Joyce's programme.
For a more comprehensive review of Joyce's programme we refer the reader to the survey paper \cite{joyce:slgcones5}.

\subsection{Deformation theory of asymptotically conical \slg submanifolds in $\C^n$}
\label{acslg:section}
This section provides definitions and results that were used in Remark \ref{ac:dimension:unbounded}.

An asymptotically conical \slg (ACSLG) submanifold in $\C^n$ is
a nonsingular \slg submanifold which tends to a \slg cone $C$ at infinity.
If $L$ is an \acslg submanifold in $\C^n$
then $\lim_{t\ra0+}{tL}=C$ for some \slg cone $C$. Hence \acslg submanifolds give rise to
local models for how nonsingular \slg submanifolds can develop isolated singularities
based on the asymptotic cone $C$. Conversely, one can use \acslg submanifolds to desingularize
isolated conical singularities of \slg varieties.


We recall from \cite[\S 7]{joyce:slgcones1} the precise definition of an \acslg submanifold.
Let $L$ be a closed, nonsingular \slg submanifold in $\C^n$.
$L$ is \textit{\ac with rate $\lambda<2$ and cone $C$} if
for some compact subset $K \subset L$ and some $T>0$, there exists a diffeomorphism
$\phi: \Sigma \times (T,\infty) \ra L-K$ such that
\begin{equation}
\label{aclsg:rate}
|\nabla^k (\phi - \iota)| = O(r^{\lambda-1-k}) \textrm{\ as\ } r\ra \infty \textrm{\ for\ } k=0,1
\end{equation}
where $\iota(r,\sigma)=r\sigma$ and $\nabla$, $|\,.\,|$ are computed using the
cone metric on $(0,\infty)\times \Sigma$.
%
Natural function spaces for elliptic theory on an \acslg submanifold $L$ are
the weighted Sobolev spaces $L^p_{k,\beta}(L)$. When acting on these spaces, the Laplacian $\Delta$ has
good mapping properties which depend in a crucial way on $\dsig$ and $N_{\Sigma}$
\cite[\S 4]{marshall}, \cite[Thm 7.9]{joyce:slgcones1}.
The analytic properties of $\Delta$ acting on weighted Sobolev
spaces on an \acslg submanifold are treated in detail in \cite{marshall}, building
on \cite{lockhart,lockhart:mcowen}.

We can now give the asymptotically conical analogue of Prop \ref{slg:bdy:symmetry}.
A similar statement appears as Prop 10.2 in \cite{joyce:lectures}, where the proof is left
as an exercise for the reader. One suspects that part (iii) holds without the
assumption that $L$ has only one end. However, in this case the proof will be more involved.
\begin{prop}
\label{acslg:symmetry}
Let $L$ be an \acslg submanifold $L$ with rate $\lambda$ and regular \slg cone $C$.
Let $\text{G}$ be the identity component of the subgroup of $\text{SU(n)}$ preserving $C$. Then

(i) $\text{G}$ admits a moment map $\mu$ and $C\subset \mu^{-1}(0)$

(ii) If $\lambda<0$, then $L\subset \mu^{-1}(0)$ and $\text{G}$ also preserves $L$

(iii) If $\lambda =0$ and $L$ has one end, then $L\subset \mu^{-1}(c)$ for some $c\in Z(\mathfrak{g}^*)$ and $\text{G}$ preserves $L$.
\end{prop}
\begin{proof}
(i) \cite[Prop 4.2]{joyce:symmetries} implies that $\text{G}$ admits a moment map $\mu$, and $C\subset \mu^{-1}(c)$ for some
$c\in Z(\mathfrak{g}^*)$. But since $G\subset \text{SU(n)}$, $\mu$ is homogeneous of degree $2$, \textit{i.e.} $\mu(tz)=t^2\mu(z)$ for any
$t>0$. But since $C$ is dilation invariant we must have $c=0$.

(ii) Let $x_i$ and $\tilde{x_i}$ ($i=1, \ldots ,2n$) denote the coordinates of the imbeddings of $C$ and $L$ into $\C^n$.
The definition of asymptotically conical (\ref{aclsg:rate}) implies
\begin{equation}
\label{coordinate:asymptotics}
|x_i - \tilde{x_i}| = O(r^{\lambda-1}) \quad \textrm{and} \quad |x_i^2 - \tilde{x_i}^2| = O(r^{\lambda})
\end{equation}
for $i=1, \ldots , 2n$.

Hence when $\lambda<0$, the restriction of any quadratic function from $\C^n$ to $L$
approaches the same value as its restriction to $C$. In particular, since $\mu\equiv0$ on $C$,
we have $\lim_{r\ra\infty}\mu(v)=0$ for any $v\in \mathfrak{g}$. Also by Lemma \ref{harmonic:momentmap},
$\mu(v)$ is a harmonic function on $L$ for any $v\in \mathfrak{g}$. But on an \ac manifold
the only harmonic function which tends to zero at each end is the zero function itself.
Hence $\mu\equiv0$ on $L$ also, and arguing as in the proof of Theorem \ref{slg:bdy:symmetry}
we see that $L$ is also $\text{G}$-invariant.

(iii) When $\lambda=0$, (\ref{coordinate:asymptotics}) implies that $\mu(v)$ is a bounded harmonic function on $L$ for every $v\in \mathfrak{g}$.
If $L$ has only one end, then the only bounded harmonic functions on $L$ are the constants.
Hence $\mu=c$ on $L$ for some $c\in Z(\mathfrak{g}^*)$. Again
this implies $L$ is $\text{G}$-invariant.
\end{proof}
\begin{remark}If $L$ has more than one end, then there exist non-constant bounded
harmonic functions on $L$. In this case to conclude that a function $\mu(v)$ is constant on $L$
it would be enough to prove that it must take the same constant value on each end.
\end{remark}
The next lemma
associates to any \slg cone $C$ a $1$-parameter family of asymptotically conical \slg submanifolds
$L_t$, asymptotic with rate $\lambda=2-n$ to the union of two \slg cones $C\cup e^{\frac{i\pi}{n}}C$. The construction
was discovered independently by the author \cite{haskins:slgcones}, Joyce \cite{joyce:symmetries} and Castro-Urbano \cite{castro:urbano}.
Lemma \ref{rotated:cone} is used in Remark \ref{ac:dimension:unbounded}.
\begin{lemma}
\label{rotated:cone}
Let $C$ be any regular \slg cone in $\C^n$ with link $\Sigma$. For any $t\neq 0\in \R$ define $L_t$ as
$$ L_t= \left\{ (z\sigma \in \C^n : \sigma\in \Sigma, \ z\in \C, \textrm{\ with\ } \Imag(z^n)=t,\  \arg{z}\in (0,\frac{\pi}{n})\right\}.$$
Then $L_t$ is an asymptotically conical \slg submanifold with rate $\lambda=2-n$, asymptotic to $C\cup e^{\frac{i\pi}{n}}C$.
$L_t$ has the same symmetry group as $C$.
\end{lemma}
\begin{proof}
The fact that $L_t$ is a \slg submanifold asymptotic to $C\cup e^{\frac{i\pi}{n}}C$ is
the content of \cite[Thm A]{haskins:slgcones}. It is a simple calculation to verify that $L_t$
has rate $\lambda=2-n$. Hence by parts (i) and (ii) of Prop \ref{acslg:symmetry}, $\text{G}$ admits
a moment map $\mu$, $C\subset \mu^{-1}(0)$ and $L\subset \mu^{-1}(0)$ and $L$ is $\text{G}$-invariant.
Alternatively one can verify by direct computation that if $\mu(v)|_C=0$, then $\mu(v)|_{L_t}=0$ also.
\end{proof}
Applying Lemma \ref{rotated:cone} to any of the \slg $T^2$-cones gives a plethora of
examples of connected complete asymptotically conical \slg submanifolds in $\C^3$ with two ends.
This is different from the well-studied codimension one setting, where it is known that
any complete connected stable minimal (in particular volume-minimizing) hypersurface in any Euclidean space
has exactly one end \cite{cao:shen:zhu}.

Within the analytic framework of weighted Sobolev spaces on \ac submanifolds one can develop the deformation theory
of \acslg submanifolds.
Let $L$ be an \acslg submanifold in $\C^n$ with rate $\lambda<2$ and cone $C$. Define the
\textit{moduli space of deformations of $L$ with rate $\lambda$}, denoted $\mathcal{M}^{\lambda}_{L}$, to be
the set of \acslg submanifolds in $\C^n$ with rate $\lambda$ and cone $C$ that
are diffeomorphic to $L$ and isotopic to $L$ as an \ac submanifold of $\C^n$.
The main result in the deformation theory of \acslg submanifolds is the following theorem
(which we used in Remark \ref{ac:dimension:unbounded}).

\begin{theorem}[\mbox{Marshall \cite[Thm 6.2.15]{marshall}, Pacini \cite[Thms 2,3]{pacini}}]
\label{acslg:deformations}
Let $L$ be an \acslg submanifold in $\C^n$ with rate $\lambda<2$ and cone $C$.
Let $\mathcal{M}^{\lambda}_{L}$ be the set of \acslg submanifolds in $\C^n$ with rate $\lambda$ and cone $C$ that
are diffeomorphic to $L$ and isotopic to $L$ as an \ac submanifold of $\C^n$.
Then

(i) If $\lambda \in (0,2)-\dsig$ then $\mathcal{M}^{\lambda}_{L}$ is a manifold with
$$ \dim{\mathcal{M}^{\lambda}_{L}} = b^1(L) - b^0(L) + N_{\Sigma}(\lambda)$$

where $\Sigma$ is the link of $C$ and $\dsig$ and $N_{\Sigma}$ are defined in \S \ref{harm:fn:cone}.

(ii) If $\lambda \in (2-n,0)$, then $\mathcal{M}^{\lambda}_{L}$ is a manifold of dimension $b^{m-1}(L)=b^1_c(L)$.
\end{theorem}

That is, the dimension of $\mathcal{M}^{\lambda}_{L}$ is purely topological when we consider only
\acslg submanifolds with sufficiently fast decay at infinity. For \ac submanifolds with
slower decay at infinity the dimension of $\mathcal{M}^{\lambda}_{L}$ depends on the spectral geometry of the asymptotic
cone via $N_{\Sigma}(\lambda)$. If we choose the rate $\lambda$ sufficiently close to $2$,
\textit{i.e.} $\max{\{\dsig\cap (0,2)\}}< \lambda < 2$, then $N_{\Sigma}(\lambda) = b^0(\Sigma) + \lind(C)$.
In \S \ref{s1inv:indexbounds} we found a countably infinite family of \acslg \ submanifolds in $\C^3$ on
which $\lind$ is unbounded, and so that each member of the family is diffeomorphic to $T^2\times \R$. By
the discussion above this provides us with topologically simple \acslg submanifolds which still occur in moduli spaces of arbitrarily
large dimension.

\subsection{Deformation theory for \slg submanifolds with isolated conical singularities}
\label{conical:sing:section}
In \cite{joyce:slgcones1}, \slg submanifolds with isolated conical singularities are defined
and in \cite{joyce:slgcones2} moduli spaces of such submanifolds are considered.
A Riemannian manifold with isolated conical singularities is
a singular Riemannian $n$-manifold $X$ with a finite number of distinct singular points
$x_1, x_2, \dots , x_N$, which near any singular point, $x_i$, looks metrically like the metric cone $C_i$
over some smooth compact nonsingular link $\Sigma_i$.
Each singular point $x_i$ has an associated rate $\mu_i$ specifying how fast the
metric on $X$ approaches the cone metric as we tend to $x_i$.
For a precise definition of a \slg variety with isolated conical singularities see \cite{joyce:slgcones1}.

Let $X$ be a compact \slg variety with conical singularities $x_1, \ldots, x_N$, cones $C_1, \ldots, C_N$
and rates $\mu_1, \ldots \mu_N$ in an almost Calabi-Yau manifold $M$
and denote $X-\{x_1, \ldots x_N\}$ by $X'$.
Let $\mathcal{M}_X$ be the set of all compact \slg submanifolds in $M$ with isolated conical singularities
which are deformation equivalent to $X$. For a precise definition of
$\mathcal{M}_X$ see \cite[\S 5]{joyce:slgcones2}. The main result \cite[Thm 6.10]{joyce:slgcones2} of the deformation theory
developed by Joyce is a description of the structure of
$\mathcal{M}_X$ near $X$ in terms of the zeroes of a smooth map between an infinitesimal deformation
space $\mathcal{I}_{X'}$ and an obstruction space $\mathcal{O}_{X'}$ of dimension $\sum_{i=1}^n{\sind(C_i)}$,
where $\sind(C_i)$ is the stability index defined in \S \ref{harm:fn:cone}.
Since the results of the current paper provide many \slg $T^2$-cones $C$ for which $\sind(C)$ is arbitrarily large,
this implies that the deformation theory for a conically singular \slg variety $X$ is in general obstructed
and that the moduli space $\mathcal{M}_X$ can fail to be a smooth manifold.

$X$ is said to have \textit{strictly stable} conical singularities if each cone $C_1, \ldots, C_N$
is strictly stable as defined in (\ref{strictly:stable}).
If $X$ has strictly stable singularities then the obstruction space $\ox=\{0\}$,
and then $\mx$ is a smooth manifold of dimension equal to $\dim{\ix}$ \cite[Cor 6.11]{joyce:slgcones2}.
Hence it becomes an important question to identify the possible strictly stable
\slg cones. Theorem \ref{ct:is:only:stable:t2} in  \S \ref{conformal:area:section}
implies that the Clifford torus cone is the unique strictly stable \slg $T^2$-cone in $\C^3$.

More generally, there is an natural transversality condition for $X\in \mx$
which guarantees that near $X$, $\mx$ is a smooth manifold of dimension
$\dim{\ix}-\dim{\ox}$ \cite[Cor 6.13]{joyce:slgcones2}. It is  conjectured that
this transversality condition holds in any almost Calabi-Yau manifold
whose K\"ahler class is sufficiently generic \cite[Conj 9.5]{joyce:slgcones2}.


\subsection{Desingularization theory and singularities in generic families of \slg submanifolds in almost \cy manifolds}
\label{generic:singularities}
Given $X\in \mx$ one can try to desingularize $X$ as follows. At each singular point $x_i$ of $X$, we glue in an appropriately scaled
\acslg submanifold $L_i$ which is asymptotic to the cone $C_i$. This gives families of almost singular, almost \slg submanifolds $N^t$ for $t$ small.
The point is then to prove by an Implicit Function Theorem type argument that for small enough $t$
there exists a perturbation of $N^t$ which is still nonsingular, but which is exactly special Lagrangian.
In general, there are several types of obstructions to performing such desingularizations.
We refer the reader to \cite{joyce:slgcones3,joyce:slgcones4,joyce:slgcones5} for the whole story.

Let $(M,J,\omega, \Omega)$ be an almost \cy manifold, and $L$ be a compact nonsingular \slg submanifold in $M$.
Then $\ml$ the moduli space of deformations of $L$ in $M$ is a smooth manifold of dimension $b^1(L)$.
In general, $\ml$ is noncompact. One natural compactification of $\ml$ is to take the closure $\mlb$ of $\ml$
in the space of integral currents. Standard compactness results in geometric measure theory imply
that that $\mlb$ is compact. Define the \textit{boundary} $\bdyml$ to be $\mlb - \ml$.
The structure of $\bdyml$ is important but very poorly understood.

One possible class of points in $\bdyml$ are the moduli spaces $\mx$ of \slg submanifolds with
isolated singularities discussed in \S \ref{conical:sing:section}. Let $X\in \mx$, and suppose
that $\mx \subset \bdyml$ and that $\mx$ is a smooth manifold. Define the \textit{index} of the
singularities of $X$ to be $\ind(X) = \dim \ml - \dim \mx$, \textit{i.e.} the codimension
of $\mx$ in $\ml$. The most common and hence most important singularities are those with small index.
For example, Joyce has pointed out that to understand whether the weighted count of \slg homology $3$-spheres
in a given homology class that he proposed in \cite{joyce:counting} is really an invariant
of almost \cy $3$-folds it is sufficient to understand only index $1$ singularities.
When $X$ is transverse, Joyce gives a formula for $\ind(X)$.
\begin{theorem}[\mbox{\cite[Thm 8.10]{joyce:slgcones5}}]
\label{index:singularity:thm}
Let $X$ be a compact, transverse \slg submanifold with conical singularities at $x_1, \ldots ,x_N$ and cones $C_1, \ldots ,C_N$.
Let $L$ be a desingularization of $X$ constructed by gluing \acslg submanifolds $L_1, \ldots L_N$ in at $x_1, \ldots x_N$.
Then
\begin{equation}
\label{index:transverse}
 \ind(X) = \dim{\mathcal{Y}} + 1- q + \sum_{i=1}^Nb^1_\text{cs}(L_i) + \sum_{i=1}^N{\sind(C_i)}.
 \end{equation}
\end{theorem}
\no
Here $q$ is the number of connected components of $X-\{x_1, \ldots ,x_N\}$. $\mathcal{Y}$ is
a particular vector subspace in the product of the first cohomology groups
of the links $\Sigma_i$ of the cones $C_i$. For a definition of $\mathcal{Y}$ see \cite[Defn 8.2]{joyce:slgcones5}.

If $\omega$ is sufficiently generic in its \ka class then any compact \slg submanifold $X$
with conical singularities in $(M,J,\omega,\Omega)$ is expected to be transverse
and hence this formula for $\ind(X)$ will be valid.
(\ref{index:transverse}) makes it clear that any cone $C$ for which $\sind(C)$ is large will
occur only in singularities of high codimension, and hence will be less important to understand.

Using (\ref{index:transverse}), Joyce studied index $1$ singularities modelled on a rigid \slg cone $C\in \C^3$ \cite[\S 10.3]{joyce:lectures}.
He showed that either $C$ is the transverse intersection of two \slg $3$-planes or $C$ is a $T^2$-cone with
$\lind(C)=6$, \textit{i.e.} $C$ is a (Legendrian) stable $T^2$-cone. Joyce raised the question of classifying all such
$T^2$-cones and conjectured that the cone on the Clifford torus (Examples \ref{clifford:torus})
and perhaps also the cone $T_{1,2}$ constructed in Theorem \ref{u:alpha:0} are the only examples.
In Theorem \ref{ct:is:only:stable:t2} of \S \ref{conformal:area:section} we proved that the Clifford torus
is the unique \slg $T^2$-cone $C$ with $\lind(C)=6$. Hence Theorem \ref{ct:is:only:stable:t2} implies the following corollary.
\begin{corollary}
\label{index:1:classification}
In dimension $3$, the cone on the Clifford torus and the transverse intersection of two \slg $3$-planes are the only rigid
\slg cones which can occur in an index $1$ singularity.
\end{corollary}
\noindent
Corollary \ref{index:1:classification} should be important for subsequent work establishing invariance properties of the counting ``invariant'' proposed
in \cite{joyce:counting}.
\newpage

\bibliographystyle{abbrv}
\bibliography{paper}
\end{document}